%% file: main.tex
\documentclass[12pt]{amsart}

\usepackage{amsmath}
\usepackage{amsfonts}
\usepackage{amssymb}
\usepackage{amscd}
\usepackage{graphicx}
\usepackage[abbrev,alphabetic]{amsrefs}
\RequirePackage[dvipsnames,usenames]{color}
\usepackage{soul,xcolor}
\setstcolor{red}
\usepackage{stmaryrd}
\usepackage{mathtools}
\usepackage{booktabs}
\usepackage{multirow}
\usepackage{stmaryrd} % \mapsfrom := the opposite of \mapsto
\newtagform{tiny}{\tiny(}{)}

\usepackage{mathtools}
\usepackage{hyperref}
\usepackage[margin=1.25in]{geometry}

\usepackage{amsthm}
\usepackage{comment}
\usepackage[all,cmtip]{xy}
\usepackage{tikz-cd}
\usetikzlibrary{cd}
\usepackage{mathtools}
\usepackage[all]{xy}

\makeatletter
\def\@tocline#1#2#3#4#5#6#7{\relax
  \ifnum #1>\c@tocdepth % then omit
  \else
    \par \addpenalty\@secpenalty\addvspace{#2}%
    \begingroup \hyphenpenalty\@M
    \@ifempty{#4}{%
      \@tempdima\csname r@tocindent\number#1\endcsname\relax
    }{%
      \@tempdima#4\relax
    }%
    \parindent\z@ \leftskip#3\relax \advance\leftskip\@tempdima\relax
    \rightskip\@pnumwidth plus4em \parfillskip-\@pnumwidth
    #5\leavevmode\hskip-\@tempdima
      \ifcase #1
       \or\or \hskip 1em \or \hskip 2em \else \hskip 3em \fi%
      #6\nobreak\relax
    \hfill\hbox to\@pnumwidth{\@tocpagenum{#7}}\par% <---- \dotfill -> \hfill
    \nobreak
    \endgroup
  \fi}
\makeatother

\newsavebox{\pullback}
\sbox\pullback{%
\begin{tikzpicture}%
\draw (0,0) -- (1ex,0ex);%
\draw (1ex,0ex) -- (1ex,1ex);%
\end{tikzpicture}}

\newsavebox{\pullbackdl}
\sbox\pullbackdl{%
\begin{tikzpicture}%
\draw (-1ex,0ex) -- (0ex,0ex);%
\draw (0ex,-1ex) -- (0ex,0ex);%
\end{tikzpicture}}

\newsavebox{\pushoutdr}
\sbox\pushoutdr{%
\begin{tikzpicture}%
\draw (-1ex,-1ex) -- (-1ex,0ex);%
\draw (-1ex,0ex) -- (0ex,0ex);%
\end{tikzpicture}}

\newcommand{\rup}[1]{\lceil #1 \rceil}
\newcommand{\rdown}[1]{\lfloor #1 \rfloor}

\newcommand{\A}{\mathbb{A}}
\renewcommand{\P}{\mathbb{P}}
\newcommand{\Z}{\mathbb{Z}}
\newcommand{\Q}{\mathbb{Q}}

\newcommand{\R}{\mathbb{R}}
\newcommand{\F}{\mathbb{F}}

\newcommand{\cHom}{\mathcal{H}om}

\newcommand{\cF}{\mathcal{F}}

\newcommand{\cL}{\mathcal{L}}

\newcommand{\cM}{\mathcal{M}}
\newcommand{\cO}{\mathcal{O}}

\newcommand{\MO}{\mathcal{O}}
\newcommand{\sO}{\mathcal{O}}

\newcommand{\m}{\mathfrak{m}}

\newcommand{\univ}{\mathrm{univ}}

\newcommand{\Frac}{\mathrm{Frac}}
\newcommand{\Proj}{\mathrm{Proj}}
\renewcommand{\Im}{\mathrm{Im}}

\DeclareMathOperator{\Br}{Br}
\DeclareMathOperator{\Gal}{Gal}

\DeclareMathOperator{\gen}{gen}
\DeclareMathOperator{\pr}{pr}
\DeclareMathOperator{\Hilb}{Hilb}
\DeclareMathOperator{\Univ}{Univ}
\DeclareMathOperator{\PGL}{PGL}

\DeclareMathOperator{\Supp}{Supp}
\DeclareMathOperator{\Spec}{Spec}

\DeclareMathOperator{\Hom}{Hom}

\DeclareMathOperator{\Pic}{Pic}
\DeclareMathOperator{\NE}{NE}
\DeclareMathOperator{\Ex}{Ex}

\DeclareMathOperator{\Ker}{Ker}

\DeclareMathOperator{\Bl}{Bl}

\theoremstyle{plain}
\newtheorem{theorem}{Theorem}[section]
\newtheorem{thm}[theorem]{Theorem}

\newtheorem{prop}[theorem]{Proposition}

\newtheorem{lem}[theorem]{Lemma}

\newtheorem{cor}[theorem]{Corollary}

\newtheorem*{claim*}{Claim}
\newtheorem{step}{Step}

\theoremstyle{definition}

\newtheorem{dfn}[theorem]{Definition}
\newtheorem{definition}[theorem]{Definition}

\newtheorem{example}[theorem]{Example}

\newtheorem{nothing}[theorem]{}
\newtheorem*{setup*}{Setup}
\newtheorem{theo}{Theorem}

\theoremstyle{remark}
\newtheorem{rem}[theorem]{Remark}



\numberwithin{equation}{theorem}

\title[Liftability and vanishing theorems for Fano threefolds]{Liftability and vanishing theorems for Fano threefolds in positive characteristic II}

\author{Tatsuro Kawakami}
\address{Graduate School of Mathematical Sciences, 
The University of Tokyo, 
3-8-1 Komaba, Meguro-ku, Tokyo 153-8914, JAPAN}
\email{tatsurokawakami0@gmail.com}
\author{Hiromu Tanaka} 
\address{Department of Mathematics, 
Graduate School of Science, 
Kyoto University, 
Kyoto 606-8502, JAPAN} 
\email{tanaka.hiromu.7z@kyoto-u.ac.jp}

\begin{document}

\begin{abstract}
In our series of papers, we prove that smooth Fano threefolds in positive characteristic lift to the ring of Witt vectors. 
Moreover, we show that they 
satisfy Akizuki-Nakano vanishing, 
$E_1$-degeneration of the Hodge to de Rham spectral sequence, 
and torsion-freeness of Crystalline cohomologies. 
In this paper, we establish these results 
except when $|-K_X|$ is very ample and the Picard group is generated by $\omega_X$.
To this end, we show that an arbitrary smooth Fano threefold is quasi-$F$-split when the Picard number or the Fano index is larger than one.
\end{abstract}

\subjclass[2020]{14J45, 13A35, 14F17} 
\keywords{Fano threefolds, Liftability to characteristic zero, Vanishing theorems
(quasi-)$F$-split, Positive characteristic}
\maketitle

\setcounter{tocdepth}{2}

\tableofcontents

\input section1.tex

\input section2.tex

\input section3.tex

\input section4.tex
\input section5.tex

\input section6.tex
\input section7.tex

\input section8.tex

\bibliographystyle{skalpha}
\bibliography{bibliography.bib}

\end{document}

%% file: section1.tex
\section{Introduction}

This paper is a continuation of \cite{Kawakami-Tanaka(Lift1)}. 
In this paper,
we prove the following theorems except when
$\Pic X = \mathbb{Z}K_X$ and $-K_X$ is very ample.

\begin{theo}[\textup{\cite[Theorem A]{Kawakami-Tanaka(Lift1)}}]\label{Introthm:W(k)-lift}
Let $X$ be a smooth Fano threefold over an algebraically closed field $k$ of positive characteristic.
Then $X$ lifts to $W(k)$.
\end{theo}

\begin{theo}[\textup{\cite[Theorem B]{Kawakami-Tanaka(Lift1)}}]\label{Introthm:ANV}
Let $X$ be a smooth Fano threefold over an algebraically closed field $k$ of characteristic $p>0$. 
% Suppose that one of the following holds:
% \begin{enumerate}
%     \item $p>2$.
%     \item $(\rho(X), r_X, (-K_X)^3) \neq (1, 1, 10)$. %i.e., g =6
% \end{enumerate}
Then Akizuki-Nakano vanishing holds on $X$, that is, 
if $A$ is an ample  Cartier divisor $A$ on $X$, then we have
\[
H^j(X,\Omega^i_X \otimes \MO_X(-A))=0
\]
for all integers $i,j\geq 0$ satisfying  $i+j<3$.
\end{theo}

\begin{theo}[\textup{\cite[Theorem C]{Kawakami-Tanaka(Lift1)}}]
\label{Introthm:E_1-degeneration}
Let $X$ be a smooth Fano threefold over an algebraically closed field $k$ of characteristic $p>0$. 
% Suppose that one of the following holds:
% \begin{enumerate}
%     \item $p>2$.
%     \item $(\rho(X), r_X, (-K_X)^3) \neq (1, 1, 10)$.
% \end{enumerate}
Then the following hold. 
\begin{enumerate}
    \item The Hodge to de Rham spectral sequence 
\[
E_1^{i,j} %{\cyan E_1^{i,j}??} 
= H^j(X, \Omega_X^i) 
% {\cyan H^i(X, \Omega_X^j)??} 
\Rightarrow 
H^{i+j}(X, \Omega_X^{\bullet}) =E^{i+j}
\]
degenerate at $E_1$.
    \item Crystalline cohomology $H^i_{\mathrm{cris}}(X/W(k))$ is torsion-free for every $i\geq 0$.
\end{enumerate}
\end{theo}

\begin{theo}[\textup{\cite[Theorem D]{Kawakami-Tanaka(Lift1)}}]\label{Introthm:hoge number}
Let $X$ be a smooth Fano threefold over an algebraically closed field $k$ of positive characteristic. 
Take a lift $f\colon \mathcal{X}\to W(k)$ of $X$ to $W(k)$, whose existence is ensured by Theorem \ref{Introthm:W(k)-lift}. 
Let $X_{\overline{K}}$ be the geometric generic fibre over $W(k)$. 
%Suppose that $(\rho(X), r_X, (-K_X)^3) \neq (1, 1, 10)$.
Then all the Hodge numbers $h^{i,j}(X)\coloneqq \dim_{k} H^j(X,\Omega^i_X)$ of $X$ coincide with those of $X_{\overline{K}}$, that is, 
\[
h^j(X,\Omega^i_X)=h^j(X_{\overline{K}},\Omega^i_{X_{\overline{K}}})
\]
hold for all $i,j\geq 0$.
\end{theo}

\subsection{Quasi-$F$-splitting}

To prove most of the statements of the main theorems, we focus on quasi-$F$-splitting, a weaker notion than $F$-splitting introduced by Yobuko \cite{Yob19} (see also \cite{KTTWYY1} for some foundational results on quasi-$F$-splitting).

Recently, Petrov proved that a smooth projective quasi-$F$-split variety satisfies the Akizuki–Nakano vanishing theorem (and hence lifts to $W(k)$ if it is Fano) and the $E_1$-degeneration of the Hodge-to-de Rham spectral sequence \cite{Petrov}. Therefore, most of the statements of the main theorems could be settled if we could prove that smooth Fano threefolds are quasi-$F$-split.
Unfortunately, there exist smooth Fano threefolds that are not quasi-$F$-split when $\rho(X) = r_X = 1$ (see \cite[Example 7.2]{KTY} and Example \ref{e-p=5-nonQFS}).
Nevertheless, we can prove the following:

\begin{theo}[\textup{=Theorem \ref{thm:quasi-F-split}}]\label{Introthm:quasi-F-split}
    Let $X$ be a smooth Fano threefold over an algebraically closed field of positive characteristic.
    If $\rho(X)>1$ or $r_X>1$, then $X$ is quasi-$F$-split.
\end{theo}

Surprisingly, Theorem \ref{thm:quasi-F-split} asserts that Fano threefolds with wild conic bundle structures (i.e., conic bundles that are not generically smooth) are all quasi-$F$-split. These varieties are known to be non-$F$-split.
Thus quasi-$F$-splitting provides a framework that can be applied even to such pathological varieties in positive characteristic, 
enabling us to establish
useful vanishing theorems.
This highlights a significant advantage of quasi-$F$-splitting.

Another remarkable property of quasi-$F$-splitting is that every quasi-$F$-split smooth variety lifts to $W_2(k)$ together with an arbitrary effective Cartier divisor \cite{AZ21}, \cite[Section 7.2]{KTTWYY1}.
From this, we can deduce the logarithmic Akizuki–Nakano vanishing theorem when $p>2$.

\begin{theo}\label{Introthm:log ANV}
Let $X$ be a smooth Fano threefold over an algebraically closed field of characteristic $p>2$ such that  $\rho(X) >1$ or $r_X >1$. 
Take a  reduced divisor $E$ with simple normal crossing support and an ample $\Q$-divisor $A$ such that the support of the fractional part of $A$ is contained in $E$. 
Then 
\[
H^j(X, \Omega^i_X(\log E) \otimes \sO_X(-\lceil A \rceil))=0
\]
holds for every pair $(i, j)$ of integers $i$ and $j$ satisfying $i+j <3$.  
\end{theo}

\subsection{$F$-splitting}
It is well known that smooth del Pezzo surfaces are $F$-split when $p>5$.
It is then natural to ask when a smooth Fano threefold is $F$-split.
In the course of the proof of Theorem \ref{thm:quasi-F-split}, we establish the following theorem.

\begin{theo}\label{intro-F-split}
Every smooth Fano threefold over an algebraically closed field of characteristic $p>5$ such that  $\rho(X) >1$ or $r_X >1$
is $F$-split.
\end{theo}
\begin{rem}
    $F$-splitting of some smooth Fano threefolds has been proven by Totaro in a different way \cite[the proof of Lemma 1.5]{Totaro(Fano)}.
\end{rem}

\begin{rem}\,
\begin{enumerate}
    \item As mentioned above, 
    if $S$ is a non-$F$-split del Pezzo surface, then $X := S \times \P^1$ is a smooth Fano threefold which is not $F$-split. 
    Since this construction {is} %can be 
    applicable  for $p \in \{2, 3, 5\}$,  
    the assumption  $p>5$ in Theorem \ref{intro-F-split} is optimal.
    \item 
If $p=7$, then the Fermat quartic hypersurface 
$\{x_0^4+x_1^4+x_2^4+x_3^4+x_4^4=0\}\subset \mathbb{P}_k^4$ is not $F$-split 
by Fedder's criterion. 
\item 
If $p=11$, then 
we shall prove that there exists a smooth Fano threefold which is not $F$-split (Example \ref{e-p=11}). 
\item 
The authors do not know whether there exists a 
non-$F$-split smooth Fano threefold of characteristic $p \geq 13$. 
\end{enumerate}
\end{rem}

\subsection{Strategy of the proof of Theorem \ref{Introthm:quasi-F-split}}
We now overview how to show that a smooth Fano threefold 
with $\rho(X) \geq 2$ is quasi-$F$-split (Theorem \ref{thm:quasi-F-split}(1)). 
In order to show that $X$ is quasi-$F$-split, 
we shall use one of the following strategies. 
\begin{enumerate}
\item[(A)] Inversion of adjunction for $F$-splitting.  
\item[(B)] Inversion of adjunction for quasi-$F$-splitting.  
\item[(C)] Cartier operator criterion for quasi-$F$-splitting. 
\end{enumerate}

(A) 
For many cases (e.g., $3 \leq \rho(X) \leq  5$ except for No. 3-10), we can prove that $X$ is $F$-split. 
For example, let us consider the case when $X$ is  of No. 2-30, i.e., 
there is a blowup $X \to \P^3$ along a smooth conic $B$ on $\P^3$.  
 Take the plane $D (\simeq \P^2)$ containing $B$. 
For the the proper transform $D_X$ of $D$ on $X$, we have the following implications: 
\[
D \text{ is $F$-split} \overset{{\rm (i)}}{\Rightarrow} 
(\P^3, D) \text{ is $F$-split} \overset{{\rm (ii)}}{\Leftrightarrow} 
(X, D_X) \text{ is $F$-split}  \overset{{\rm (iii)}}{\Rightarrow}  X  \text{ is $F$-split}.  
\]
Of course, $D (\simeq \P^2)$ is $F$-split. 
The implication (i) follows from the inversion of adjunction for $F$-splitting, 
e.g., if $(Y, D)$ is $F$-split and $-(K_Y+D)$ is ample, 
then $(Y, D)$ is $F$-split. 
The equivalence (ii) is assured by $K_X + D_X = f^*(K_{\P^3}+ D)$. 
Finally, (iii) holds by definition.

(B) 
Even if the strategy (A) does not work, 
we can apply a quasi-$F$-split version of (A) for most of the remaining cases. 
The authors has proved that an inversion of adjunction for log Calabi-Yau pairs in \cite{Kawakami-Tanaka(dPvar)}.
This allows us to prove the following statement:
a smooth Fano threefold $X$ is quasi-$F$-split if there are smooth prime divisors $S$ and $S'$ such that $K_X+S+S' \sim 0$ and $S \cap S'$ is a smooth curve which is quasi-$F$-split (Corollary \ref{c-IOA-Fano3}). 

As other technical issues, we shall encounter the following obstructions: 
\begin{itemize}
\item $S$ is not necessarily smooth. For some cases, the generic member  is enough as a replacement (cf.\ Section \ref{ss-generic}). To this end, we shall need to treat algebraic varieties defined over an imperfect field. 
\item If the first prime divisor $S$ is a smooth weak del Pezzo surface, 
then it is often hard to find $S'$ such that $S \cap S'$ is smooth. 
For example, if $|-K_S|$ has no smooth member, then $S \cap S'$ can not be smooth. 
In order to avoid such a pathological phenomenon, 
we shall establish some properties on weak del Pezzo surfaces, 
e.g., if $V$ is a surface over a $C_1$-field $K$ of characteristic two and 
its base change $V \times_{\Spec K} \Spec \overline K$ is a  Langer surface (defined as the base change of the blowup of $\P^2_{\F_2}$ along all the $\F_2$-rational points), then $\rho(V) =8$ 
(Lemma \ref{l-Langer-rigid}).
\end{itemize}

(C) Except when $X$ is one of 2-2, 2-6, 2-8, and 3-10, 
we may apply one of (A) and (B).  
For these remaining cases, we shall apply a quasi-$F$-splitting criterion via Cartier operator, 
which has been established in \cite[Theorem F]{KTTWYY1} (cf. Proposition \ref{prop:criterion for qFs}). 
To this end, we need to prove $H^j(X, \Omega^i_X \otimes \MO_X(p^{\ell}K_X))=0$ for suitable triples $(i, j, \ell)$. 
Even if $X$ is explicitly given, 
it is often hard to compute such cohomologies directly. 
The main strategy is to embed $X$ into a (typically toric) fourfold $P$, 
and apply Bott vanishing for $P$. 
For example,  if $X$ is of No.~2-2, 
then we can find such an embedding with $P = \P_{\P^1 \times \P^2}(\MO_{\P^1 \times \P^2} \oplus \MO_{\P^1 \times \P^2}(-1, -2))$. 
Although $P$ is not necessaliry toric for the other cases,  
we shall find a closed embedding $X \hookrightarrow P$ to a fourfold $P$ which 
almost satisfies  Bott vanishing. 
For more details on (C), see Section \ref{sec:8}.

\medskip
\noindent {\bf Acknowledgements.}
The authors express their gratitude to Burt Totaro for valuable comments.
They also thank Teppei Takamatsu and Shou Yoshikawa for useful conversations.
Kawakami was supported by JSPS KAKENHI Grant number JP22KJ1771 and JP24K16897.
Tanaka was supported by JSPS KAKENHI Grant number JP22H01112 and JP23K03028.

%% file: section2.tex
\section{Preliminaries}

\subsection{Notation}
\label{ss:notation}

In this subsection, we summarise notation and basic definitions used in this article. 
\begin{enumerate}
\item Throughout the paper, $p$ denotes a prime number and we set $\F_p \coloneqq \Z/p\Z$. 
Unless otherwise specified, we work over an algebraically closed field $k$ of characteristic $p>0$. 
We denote by $F \colon X \to X$ the absolute Frobenius morphism on an $\F_p$-scheme $X$.
    \item We say that $X$ is a {\em variety} (over a field $\kappa$) if 
    $X$ is an integral scheme 
    that is separated and of finite type over $\kappa$. 
    We say that $X$ is a {\em curve} (resp. {\em surface}, resp. {\em threefold}) 
    if $X$ is a variety of dimension one (resp. two, resp. three). 
\item For a variety $X$, 
we define the {\em function field} $K(X)$ of $X$ 
as the stalk $\MO_{X, \xi}$ at the generic point $\xi$ of $X$. 
\item We say that an $\R$-divisor $D$ on a normal variety $X$ is \emph{simple normal crossing} if for every point $x \in \Supp D$, the local ring $\cO_{X,x}$ is regular and there exists a regular system of parameters $x_1,\ldots, x_d$ of the maximal ideal $\m$ of $\cO_{X,x}$ and $1 \leq r \leq d$ such that $\Supp (D|_{\Spec \MO_{X, x}}) = \Spec (\cO_{X,x}/(x_1 \cdots x_r))$. 
\item Given an integral normal Noetherian scheme $X$, a projective birational morphism $\pi \colon Y \to X$ is called \emph{a log resolution (of singularities) of $X$} if $Y$ is regular and $\mathrm{Exc}(f)$ is a simple normal crossing divisor.
\item We say that an $\F_p$-scheme $X$ is {\em $F$-finite} if the absolute Frobenius morphism $F: X \to X$ is a finite morphism. 
We say that an $\F_p$-algebra $R$ is {\em $F$-finite} if $\Spec R$ is $F$-finite. 
In particular, a field $\kappa$ is $F$-finite if and only if $[\kappa:\kappa^p]<\infty$. 
If $X$ is a variety over an $F$-finite field, then $X$ is $F$-finite. 
\item 
Given a normal variety $X$ and an $\R$-divisor $D$, 
we define the subsheaf $\MO_X(D)$ of the constant sheaf $K(X)$ on $X$ 
{by the following formula}
\[
\Gamma(U, \MO_X(D)) = 
\{ \varphi \in K(X) \mid 
\left({\rm div}(\varphi)+D\right)|_U \geq 0\}
\]
{for every} open subset $U$ of $X$. 
In particular, %we denote 
$\cO_X(\rdown{{D}}) = \cO_X({D})$.
\item 
Given a field $K$ and $K$-schemes $X$ and $Y$, 
we say that $X$ is {\em $K$-isomorphic} to $Y$ if there exists an isomorphism $\theta : X \to Y$ over $K$. 
\item 
Given a closed subscheme $X$ of $\P^n$, we set $\MO_X(a) := \MO_{\P^n}(a)|_X$ for $a \in \Z$ unless otherwise specified. 
Similarly, if $Y$ is a closed subscheme of $\P^n \times \P^m$, then we define $\MO_Y(a, b) \coloneqq \MO_{\P^n \times \P^m}(a, b)$ for $a, b \in \Z$. 
\item 
Given a coherent sheaf $\cF$ and a Cartier divisor $D$ on a variety $X$, 
we set $\cF(D) \coloneqq \cF \otimes \MO_X(D)$ unless otherwise specified. 
Note that $B_n\Omega_X^i(p^nD)$ (resp.\ $Z_n\Omega_X^i(p^nD)$) 
does not mean $B_n\Omega_X^i \otimes \MO_X(p^nD)$ (resp.\ $Z_n\Omega_X^i \otimes \MO_X(p^nD)$) even if $D$ is Cartier (cf.\ Subsection \ref{ss-Cartier-op}). 
\item Given two closed subschemes $Y$ and $Z$ on a scheme $X$, we denote by $Y \cap Z$ the scheme-theoretic intersection, i.e., 
$Y\cap Z \coloneqq Y \times_X Z$. 
\end{enumerate}

\subsection{Cartier operators}\label{ss-Cartier-op}

In this section, we recall the fundamental facts on the higher Cartier operators (\cite{Ill79}, \cite{KTTWYY1}).

Let $X$ be a smooth variety over a perfect field of characteristic $p>0$ and $D$ a Cartier divisor on $X$.
The Frobenius pushforward of the de Rham complex
\[
F_{*}\Omega^{\bullet}_X\colon  F_{*}\sO_X \xrightarrow{F_{*}d} F_{*}\Omega_X \xrightarrow{F_{*}d} \cdots
\]
is a complex of $\sO_X$-modules.
Tensoring with $\sO_X(D)$, we obtain a complex 
\[
F_{*}\Omega^{\bullet}_X\colon  F_{*}\sO_X(pD) \xrightarrow{F_{*}d\otimes \sO_X(D)} F_{*}\Omega_X(pD) \xrightarrow{F_{*}d\otimes \sO_X(D)} \cdots
\]

We define locally free $\sO_X$-modules as follows.
\[
\begin{array}{rl}
&B^1\Omega^i_X(pD)\coloneqq\Im(F_{*}d : F_{*}\Omega^{i-1}_X(pD)\to F_{*}\Omega^i_X(pD)),\\
&Z_1\Omega^i_X(pD)\coloneqq\Ker(F_{*}d : F_{*}\Omega^{i}_X(pD) \to F_{*}\Omega^{i+1}_X(pD)).\\
\end{array}
\]

We have an isomorphism
\[
Z_1\Omega_X^i(pD)/B_1\Omega_X^i(pD) \overset{C(D)}{\simeq} \Omega_X^i(D).
\]
resulting from the Cartier isomorphism.
In fact, tensoring with $\sO_X(D)$ with the usual Cartier isomorphism
\[
Z_1\Omega_X^i/B_1\Omega_X^i \overset{C}{\simeq} \Omega_X^i,
\]
we obtain the above isomorphism.

Taking the Frobenius pushforward, we obtain
\[
F_{*}Z_1\Omega_X^i(p^2D) \to F_{*}Z_1\Omega_X^i(p^2D)/F_{*}B_1\Omega_X^i(p^2D) \overset{F_{*}C(D)}{\simeq} F_{*}\Omega_X^i(pD).
\]
We denote by $B_2\Omega_X^i(p^2D)$ and $Z_2\Omega_X^i(p^2D)$ the preimages of $B_1\Omega_X^i(pD)\subset F_{*}\Omega_X^i(pD)$ and $Z_1\Omega_X^i(pD)\subset F_{*}\Omega_X^i(pD)$ by the above map.
Inductively, we define locally $\sO_X$-module $B_n\Omega_X^i(p^nD)$ and $Z_n\Omega_X^i(p^nD)$ for all $n\geq 0$.
Moreover, we set $B_0\Omega_X^i(D)=0$ and $Z_0\Omega_X^i(pD)=\Omega_X^i(pD)$.

\begin{lem}
    Then we have the following exact sequences
\begin{align}
0 \to B_n\Omega^i_X(p^nD) \to  Z_n\Omega^i_X(p^nD) \to \Omega^i_X(D) \to 0. \label{exact:B}
\end{align}
\begin{align}
0 \to Z_n\Omega_X^{i}(p^nD) \to F_{*}Z_{n-1}\Omega_X^i(p^nD) \to B_1\Omega_X^{i+1}(pD)\to 0. \label{exact:Z}
\end{align}
for all $i\geq 0$ and all $n\geq 1$.
\end{lem}
\begin{proof}
The assertion follows from \cite[(5.7.1) and Lemma 5.8]{KTTWYY1}.
\end{proof}
\begin{rem}
Taking $n=1$, we have the following exact sequence:
\begin{align}
0 \to B_1\Omega^i_X(pD) \to  Z_1\Omega^i_X(pD) \xrightarrow{C(D)}\Omega^i_X(D) \to 0. \label{exact:B_1}
\end{align}
\begin{align}
0 \to Z_1\Omega_X^{i}(pD) \to F_{*}\Omega_X^i(pD) \xrightarrow{F_{*}d\otimes \sO_X(D)} B_1\Omega_X^{i+1}(pD)\to 0. \label{exact:Z_1}
\end{align}
for all $i\geq 0$.
\end{rem}
\begin{rem}\label{rem:tensoring Z-divs}
Taking $D=0$, we have short exact sequences
\begin{align}
0 \to B_n\Omega^i_X \to  Z_n\Omega^i_X \xrightarrow{C^{n}} \Omega^i_X\to 0,\label{exact:B,D=0}
\end{align}
\begin{align}
0 \to Z_n\Omega_X^{i} \to F_{*}Z_{n-1}\Omega_X^i\xrightarrow{F_{*}d\circ F_{*}C^{n-1}} B_1\Omega_X^{i+1}\to 0, \label{exact:Z,D=0}
\end{align}
which coincides with \cite[(2.15.1) and Lemma 2.16]{KTTWYY1} respectively. 
Then we can confirm that \[
\eqref{exact:B}=\eqref{exact:B,D=0}\otimes \sO_X(D)\,\,\text{and}\,\,\eqref{exact:Z}=\eqref{exact:Z,D=0}\otimes \sO_X(D)
\]
holds for all $i\geq 0$.
In particular,
\begin{align*}
    &B_n\Omega_X^{i}(p^nD)= B_n\Omega_X^{i}\otimes \sO_X(D)\,\,\text{and}\\
    &Z_n\Omega_X^{i}(p^nD)= Z_n\Omega_X^{i}\otimes \sO_X(D)
\end{align*}
hold for all $n\geq 0$.
\end{rem}

\subsection{Generic members}\label{ss-generic}

Let $X$ be a regular projective variety $X$ over a field $k$ 
and let $D$ be a Cartier divisor on $X$ satisfying $h^0(X, \MO_X(D)) \geq 2$. 
For a base point free linear system $\Lambda \subset |D|$ and the corresponding linear subspace $V_{\Lambda} \subset H^0(X, \MO_X(D))$, 
the {\em generic member}  $X^{\gen}_{\Lambda}$ of $\Lambda$ is defined by the following diagram: 
\[
\begin{tikzcd}
	X^{\gen}_{\Lambda}  & X^{\univ}_{\Lambda} \\
	X \times_k \kappa & X \times_k (\mathbb P^n_k)^* & X\\
	\Spec \kappa & (\mathbb P^n_k)^* & \Spec\,k
	\arrow[from=3-2, to=3-3]
	\arrow[from=3-1, to=3-2, "\Theta"]
	\arrow[from=2-3, to=3-3]
	\arrow[from=2-2, to=3-2, "{\rm pr}_2"]
	\arrow[from=2-1, to=3-1]
	\arrow[from=2-1, to=2-2]
	\arrow[hook, from=1-2, to=2-2]
	\arrow[hook, from=1-1, to=2-1]
	\arrow[from=1-1, to=1-2]
	\arrow[from=1-2, to=2-3]
	\arrow[from=2-2, to=2-3, "{\rm pr}_1"]
\end{tikzcd}\qquad \kappa := K((\mathbb P^n_k)^*)
\]
where 
\begin{enumerate}
\item $X^{\univ}_{\Lambda}$ denotes the universal family that parametrises all the members of $\Lambda$, 
%effective Cartier divisors on $X$ linearly equivalent to 
\item $\kappa = K((\mathbb P^n_k)^*)$ is the function field 
of the projective space $(\mathbb P^n_k)^*$ 
and 
$\Theta: \Spec K((\mathbb P^n_k)^*) \to (\mathbb P^n_k)^*$ is the induced morphism. 
\end{enumerate}
Then the following hold. 
\begin{enumerate}
\setcounter{enumi}{2}
    \item $\kappa / k$ is a purely transcendental extension of finite transcendence degree. 
    \item $X \times_k \kappa$ is a regular projective variety. 
    \item $X^{\gen}_{\Lambda}$ is a regular prime divisor \cite[Theorem 4.9(4)(12)]{Tan-Bertini}. 
\end{enumerate}
For more details, we refer to \cite{Tan-Bertini}. 
By abuse of notation, also $(X^{\gen}_{\Lambda}) \times_{\kappa} \kappa'$ is called the generic member when $\kappa'/\kappa$ is a purely transcendental extension, 
because we shall encounter the situation as in the following remark. 

\begin{rem}\label{r-generic-snc}
We now consider the case when we have two base point free linear systems $\Lambda_1$ and $\Lambda_2$ on $X$. 
As above, we obtain two generic members $X^{\gen}_{\Lambda_1}$ on $X \times_{\kappa} \kappa_1$ 
and $X^{\gen}_{\Lambda_2}$ on $X \times_{\kappa} \kappa_2$. 
For $\kappa_1 = k(s_1, \ldots, s_a)$ and $\kappa_2 =k(t_1,\ldots, t_b)$, i.e., 
each of $\{s_i\}$ and $\{t_j\}$ is a  transcendental basis, 
we set 
\[
\kappa := \Frac (k(s_1, \ldots, s_a) \otimes_k k(t_1, \ldots, t_b)) =k(s_1,\ldots, s_a, t_1,\ldots, t_b). 
\]
For 
$(X^{\gen}_{\Lambda_1})_{\kappa} := X^{\gen}_{\Lambda_1} \times_{\kappa_1} \kappa$ and 
$(X^{\gen}_{\Lambda_2})_{\kappa} := X^{\gen}_{\Lambda_2} \times_{\kappa_2} \kappa$, 
\begin{enumerate}
\item[$(\star)$] the sum $(X^{\gen}_{\Lambda_1})_{\kappa} + (X^{\gen}_{\Lambda_2})_{\kappa}$ is 
a simple normal crossing divisor. 
\end{enumerate}
%anaka, Bertini, Proposition 5.10??]. 
The similar statement holds even if we start with finitely many base point free linear systems 
$\Lambda_1,\ldots, \Lambda_r$ on $X$, i.e., the sum 
\[
(X^{\gen}_{\Lambda_1})_{\kappa} + \cdots + (X^{\gen}_{\Lambda_r})_{\kappa} 
\]
of the generic members 
$(X^{\gen}_{\Lambda_1})_{\kappa}, \ldots, (X^{\gen}_{\Lambda_r})_{\kappa}$ is simple normal crossing, 
where $\kappa := \Frac(\kappa_1 \otimes_k \cdots \otimes_k \kappa_r)$ and $(-)_{\kappa}$ denotes the base change to $\kappa$. 
\end{rem}

\begin{proof}[Proof of $(\star)$]
Both $(X^{\gen}_{\Lambda_1})_{\kappa}$ and $(X^{\gen}_{\Lambda_2})_{\kappa}$ are clearly regular prime divisors. 
It suffices to show that 
the scheme-theoretic intersection $(X^{\gen}_{\Lambda_1})_{\kappa} \cap (X^{\gen}_{\Lambda_2})_{\kappa}$ is regular. 
For each $i \in \{1, 2\}$, 
let $D_i$ be the Cartier divisor with $\Lambda_i \subset |D_i|$ 
and let $V_i \subset H^0(X, \MO_X(D_i))$ be the $k$-vector subspace corresponding to $\Lambda_i$. 
Consider the restriction map: 
\[
\rho : H^0(X \times_k \kappa_1, \MO_{X \times_k \kappa_1}(D_2 \times_k  \kappa_1)) 
\to 
H^0(X^{\gen}_{\Lambda_1}, \MO_{X \times_k \kappa_1}(D_2 \times_k  \kappa_1)|_{X^{\gen}_{\Lambda_1}}). 
\]
Then the generic member $(X^{\gen}_{\Lambda_2})_{\kappa}$ coincides with the base change 
of the generic member of $V_2 \times_k \kappa \subset H^0(X \times_k \kappa_1, \MO_{X \times_k \kappa_1}(D_2 \times_k  \kappa_1))$. 
Therefore, the intersection $(X^{\gen}_{\Lambda_1})_{\kappa} \cap (X^{\gen}_{\Lambda_2})_{\kappa}$ 
coincides with the generic member of ${\rm Image}(\rho)$ by \cite[Proposition 5.10(2)]{Tan-Bertini}, 
which is regular \cite[Theorem 4.9(4)]{Tan-Bertini}. 
\end{proof}

\subsection{$F$-splitting criteria}

\begin{dfn}
Let $X$ be a normal variety and let $\Delta$ be an effective $\Q$-divisor on $X$. 
\begin{enumerate}
\item We say that $(X, \Delta)$ is {\em $F$-split} if 
\[
\MO_X \xrightarrow{F^e} F_*^e\MO_X \hookrightarrow F_*^e\MO_X( \rdown{(p^e-1)\Delta})
\]
splits as an $\MO_X$-module homomorphism for every $e \in \Z_{>0}$. 
\item We say that $(X, \Delta)$ is {\em sharply $F$-split} if 
\[
\MO_X \xrightarrow{F^e} F_*^e\MO_X \hookrightarrow F_*^e\MO_X( \rup{(p^e-1)\Delta})
\]
splits as an $\MO_X$-module homomorphism for some $e \in \Z_{>0}$. 
\item 
We say that $(X, \Delta)$ is {\em globally $F$-regular} 
if, given an effective $\Z$-divisor $E$, there exists $e \in \Z_{>0}$ such that 
\[
\MO_X \xrightarrow{F^e} F_*^e\MO_X \hookrightarrow 
F_*^e\MO_X( \rup{(p^e-1)\Delta} +E)
\]
splits as an $\MO_X$-module homomorphism. 
\end{enumerate}
We say that $X$ is {\em $F$-split} (resp. {\em globally $F$-regular}) 
if so is $(X, 0)$. 
\end{dfn}

\begin{rem}
We have the following implications: 
\[
\text{globally $F$-regular}\Longrightarrow
\text{sharply $F$-split}\Longrightarrow
\text{$F$-split} 
 \]
 where the former implication is clear and the latter one holds by 
 the same argument as in  \cite{Sch08g}*{Proposition 3.3}. 
 Moreover, 
 if the condition ($\star$) holds, then 
 $(X, \Delta)$ is sharply F-split if and only if $(X, \Delta)$ is $F$-split. 
 \begin{enumerate}
    \item[($\star$)] $(p^e-1)\Delta$ is a $\Z$-divisor for some $e \in \Z_{>0}$
 \end{enumerate}
In particular, $X$ is $F$-split if and only if $F\colon\MO_X \to F_*\MO_X$ splits as an $\MO_X$-module homomorphism.  
 For more foundational properties, we refer to \cite{SS10}. 
\end{rem}

In what follows, we summarise some $F$-splitting criteria, which are well known to experts. 

\begin{prop}\label{p-F-split-blowup}
Let $f\colon X \to Y$ be a birational morphism of normal projective varieties. 
Take an effective $\Q$-divisor $\Delta_Y$ on $Y$ such that 
$(p^e-1)(K_Y+\Delta_Y)$ is Cartier for some $e \in \Z_{>0}$. 
Assume that the $\Q$-divisor $\Delta$ defined by $K_X+\Delta = f^*(K_Y+\Delta_Y)$ is effective. 
Then $(X, \Delta)$ is  $F$-split if and only if $(Y, \Delta_Y)$ is $F$-split.  
\end{prop}

\begin{proof}
If $(X, \Delta)$ is $F$-split, then so is $(Y, \Delta_Y)$ (take the pushforward). 
As for the opposite implication, the same argument as \cite{HX15}*{the first paragraph of the proof of Proposition 2.11} works. 
\end{proof}

\begin{prop}\label{p-F-split-IOA}
Let $\kappa$ be an $F$-finite field of characteristic $p>0$. 
Let $X$ be a normal Gorenstein projective variety over $\kappa$.  
Take a normal prime Cartier divisor $S$ 
and an effective $\Q$-Cartier $\Q$-divisor $B$ on $X$ such that $S \not\subset \Supp\,B$. 
Assume that 
\begin{enumerate}
\item $(S, B|_S)$ is $F$-split, and 
\item there exists $e \in \Z_{>0}$ such that 
$(p^e-1)(K_X+S+B)$ is Cartier and 
\[
H^1(X, \MO_X(-S-(p^e-1)(K_X+S+B)))=0. 
\]
\end{enumerate}
Then $(X, S+B)$ is $F$-split. 
\end{prop}

\begin{proof}
The same argument as in \cite[Lemma 2.7]{CTW17} works. 
\end{proof}

\begin{cor}\label{c-F-split-IOA}
Let $\kappa$ be an $F$-finite field of characteristic $p>0$. 
Let $X$ be a normal {Gorenstein} projective variety over $\kappa$. 
Take a normal prime Cartier divisor $S$ such that $S$ is $F$-split and $-(K_X+S)$ is ample. 
Then $(X, S)$ is $F$-split. 
\end{cor}

\begin{proof}
By applying Proposition \ref{p-F-split-IOA} with $B=0$, it suffices to show that $H^1(X, \sO_X(-S-(p^e-1)(K_X+S))))=0$ for $e \gg 0$, 
which follows from Serre vanishing.
\end{proof}

\begin{rem}
    In the setting of Corollary \ref{c-F-split-IOA}, even if $S$ is quasi-$F$-split, $X$ is not necessarily quasi-$F$-split \cite[Example 7.7]{KTY}.
\end{rem}

\begin{cor}\label{c-P1-2pts-Fsplit}
Let $X$ be a  {smooth} Fano threefold over $k$. 
Assume that there exist a field extension $k \subset \kappa$ and 
effective divisors $D_1, D_2, D_3$ on $X_{\kappa} := X \times_k \kappa$ such that 
\begin{enumerate}
\item $\kappa$ is an $F$-finite field, 
\item $K_{X_{\kappa}}+D_1+D_2+D_3 \sim 0$, 
\item $D_1$ is a normal prime divisor, 
\item $D_1 \cap D_2$ is one-dimensional, smooth, and
\item $D_1 \cap D_2 \cap D_3$ is non-empty, zero-dimensional, and smooth over $\kappa$, and 
\item $H^1(X_{\kappa}, \MO_{X_{\kappa}}(-D_1)) = H^1(X_{\kappa}, \MO_{X_{\kappa}}(-D_2)) = H^1(X_{\kappa}, \MO_{X_{\kappa}}(-D_3))=0$ (cf.\ Lemma \ref{l-P1-2pts-Fsplit}).
\end{enumerate}
Then $(X \times_k \kappa, D_1+D_2+D_3)$ is $F$-split. In particular, $X$ is $F$-split. 
\end{cor}

\begin{proof}
First of all, we  show that 
\begin{equation}\label{e1-P1-2pts-Fsplit}
H^1(D_1, \MO_{D_1}(-D_2|_{D_1}))=0. 
\end{equation}
To this end, it suffices to prove $H^1(X_{\kappa}, \MO_X(-D_2))=0$ and $H^2(X_{\kappa}, \MO_{X_{\kappa}}(-D_2-D_1))=0$. 
The former one follows from (6). 
The latter one holds by  
\[
h^2(X_{\kappa}, \MO_{X_{\kappa}}(-D_2-D_1)) = h^1(X_{\kappa},  \MO_{X_{\kappa}}(K_X + D_1+D_2)) 
\overset{{\rm (2)}}{=}h^1(X_{\kappa},  \MO_{X_{\kappa}}(-D_3))\overset{{\rm (6)}}{=}0.
\]
This completes the proof of (\ref{e1-P1-2pts-Fsplit}).

By $H^1(X_{\kappa}, \MO_{X_{\kappa}} (-D_1))\overset{{\rm (6)}}{=}0$ and (\ref{e1-P1-2pts-Fsplit}), 
we obtain 
\[
\kappa = H^0(X_{\kappa}, \MO_{X_{\kappa}}) \xrightarrow{\simeq} H^0(D_1, \MO_{D_1}) \xrightarrow{\simeq} 
H^0(D_1 \cap D_2, \MO_{D_1 \cap D_2}). 
\]
Therefore, $C := D_1 \cap D_2$ is a smooth projective curve. 
By (4) and (5), $C$ and $C \cap D_3$ are smooth over $\kappa$. 
In particular, $(C \times_{\kappa} \overline \kappa, (C \cap D_3) \times_{\kappa} \overline{\kappa}) \simeq (\P^1_{\overline \kappa}, P+Q)$ for 
the algebraic closure $\overline \kappa$ of $\kappa$ and some distinct points $P$ and $Q$. 
Then $(C, C \cap D_3)$ is $F$-split.

By Proposition \ref{p-F-split-IOA}, $(D_1, (D_2+D_3)|_{D_1})$ is $F$-split, 
where  Proposition \ref{p-F-split-IOA} is applicable by (\ref{e1-P1-2pts-Fsplit}). 
Again by Proposition \ref{p-F-split-IOA} and (6), $(X \times_k \kappa, D_1+D_2+D_3)$ is $F$-split.  
\end{proof}

\begin{lem}\label{l-P1-2pts-Fsplit}
Let $X$ be a smooth Fano threefold over $k$ and let $k \subset \kappa$ be a field extension. 
Take a divisor on $X_{\kappa} \coloneqq X \times_k \kappa$. 
Assume that one of the following conditions hold. 
\begin{enumerate}
\item $D$ is nef and $\nu(X_{\kappa}, D) \geq 2$. 
\item There exists a morphism $\pi \colon X \to \P^1_{\kappa}$ such that 
$\pi_*\MO_X = \MO_{\P^1_{\kappa}}$ and $\MO_{X_{\kappa}} (D) \simeq \pi^*\MO_{\P^1}(1)$. 
\end{enumerate}
Then $H^1(X_{\kappa}, \MO_{X_{\kappa}}(-D))=0$. 
\end{lem}

\begin{proof}
Taking the base change to the algebraic closure $\overline{\kappa}$ of $\kappa$, 
we may assume that $\kappa$ is algebraically closed. 
Then the problem is reduced to the case when $k=\kappa$. 
If (1) holds, then we may apply \cite[Corollary 3.6]{Kaw1}. 
If (2) holds, then we may assume that $D$ is a general fibre of $\pi \colon X \to \P^1$, 
which is a prime divisor, and hence $H^1(X, \MO_X)=0$ implies $H^1(X, \MO_{X}(-D))=0$. 
\end{proof}

\begin{rem}\label{r-g0-to-P1}
Let $\kappa$ be an $F$-finite field and let $C$ be a Gorenstein projective curve over $\kappa$. 
If $-K_C$ is ample and there exists a Cartier divisor $D$ satisfying $\deg D =1$, 
then $C \simeq \P^1_{\kappa}$ \cite[Lemma 10.6]{Kol13}, and hence $C$ is $F$-split. 
\end{rem}

\begin{example}\label{e-toric}
If $X$ is a normal toric variety and $D$ is a torus-invariant reduced divisor, 
then $(X, D)$ is $F$-split. 
This  follows essentially from \cite[2.6]{Fuj07} (cf.\ \cite[Proposition 2.17]{Tan-toric}). 
\end{example}

\begin{prop}\label{prop:criterion for F-split}
    Let $X$ be a smooth Fano threefold.
    Suppose that the following condition.
    \begin{enumerate}
        \item $H^0(X, \Omega_X^2(pK_X))=0$.
        \item $H^1(X, \Omega_X^1(K_X))=0$.
        \item $H^2(X, \Omega_X^1(-pK_X))=0$.
        \item $H^1(X, \Omega_X^2(-pK_X))=0$.
    \end{enumerate}
    Then $X$ is $F$-split.
\end{prop}
\begin{proof}
    Firstly, the global $F$-splitting of $X$ (i.e., splitting of $F\colon \MO_X \to F_*\MO_X$) is equivalent to the surjectivity of the evaluation map
    \[
    \mathrm{Hom}_{\sO_X}(F_{*}\sO_X, \sO_X) \xrightarrow{F^*}  
    \mathrm{Hom}_{\sO_X}(\sO_X, \sO_X) (\simeq H^0(X, \sO_X)).
    \]
    By Grothendieck duality, they are equivalent to the injectivity of 
    \[
    F: H^3(X, \omega_X) \to H^3(X, F_{*}\sO_X\otimes\omega_X).
    \]
    Thus it suffices to show that $H^2(X, B_1\Omega_X^1(pK_X))=0$.
    By \eqref{exact:B_1} and %the condition 
    (2), this vanishing can be reduced to $H^2(X, Z_1\Omega^1_X(pK_X))$.
    By (\ref{exact:Z_1}), it is enough to prove %the vanishing is reduced to 
    \[
    H^1(X, B_1\Omega^2_X(pK_X))=0\,\,\,\text{and}\,\,\,
    H^2(X, \Omega^1_X(pK_X))=0.
    \]
    The latter one follows from (4) and Serre duality.
    It suffices to show the first one. 
    By \eqref{exact:B_1} and the condition (1), this vanishing is reduced to 
    \[
    H^1(X, Z_1\Omega_X^2(pK_X))=0.
    \]
    By \ref{exact:Z_1}, it suffices to show
    \[
    H^0(X, B_1\Omega^3_X(pK_X))=0\,\,\,\text{and}\,\,\,
    H^1(X, \Omega^2_X(pK_X))=0.
    \]
    The first vanishing holds by 
     $B_1\Omega_X^3(pK_X)\subset F_{*}\omega^{p+1}_X$ and the ampleness of $\omega^{-1}_X$. 
    The latter one follows from (3) and Serre duality. 
\end{proof}

\subsection{Quasi-$F$-splitting criteria}

We recall that definition of $F$-splitting and quasi-$F$-splitting.
We refer to \cite[Section 3]{Kawakami-Tanaka(dPvar)} for details.

\begin{dfn}
Let $X$ be a normal variety.
We define a $W_n\MO_X$-module $Q^{e}_{X, n}$ and a $W_n\MO_X$-module homomorphism $\Phi^{\Delta,e}_{X, \Delta, n}$ by the following pushout diagram: 
\[
\begin{tikzcd}
W_n\sO_X \arrow[r, "F^e"] \arrow[d, "R^{n-1}"'] & F_*^eW_n\sO_X \arrow[d]\\
\sO_X \arrow[r, "\Phi^{e}_{X,n}"] & Q^{e}_{X,n}. 
\end{tikzcd}
\]
Applying $(-)^* \coloneqq \cHom_{W_n\MO_X}(-, W_n\omega_X(-K_X))$ to $\Phi^{e}_{X, n}$, we have a $W_n\MO_X$-module homomorphism
\[
(\Phi^{e}_{X, n})^* \colon (Q^{e}_{X,n})^* \to \sO_X. 
\]
We say that $X$ is {\em quasi-$F^e$-split} if 
\[
H^0(X, (\Phi^{e}_{X, n})^*) \colon
H^0(X, (Q^{e}_{X, n})^*) \to H^0(X, \sO_X) 
\]
is surjective.
\end{dfn}
\begin{rem}
Let $X$ be a normal variety.
 Then $X$ is $F$-split if and only if 
    \[
H^0(X, F^*) \colon
H^0(X, (F_{*}\sO_X)^*) \to H^0(X, \sO_X) 
\]
is surjective for every $e>0$, where $(-)^* \coloneqq \cHom_{\sO_X}(-, W_1\omega_X(-K_X))=\cHom_{\sO_X}(-, \sO_X)$.
In particular, if $X$ is $F$-split, then $X$ is quasi-$F^e$-split for every $e>0$.
\end{rem}

\begin{cor}\label{c-IOA-Fano3}
Let $k$ be an algebraically closed field of characteristic $p>0$ and let $k \subset \kappa$ be a field extension, where $\kappa$ is $F$-finite. 
Let $X$ be a smooth Fano threefold over $k$. %regular projective variety over $\kappa$. 
Assume that there exist a prime divisor $D$ and a reduced divisor $D'$ on $X \times_k \kappa$ such that 
\begin{enumerate}
\renewcommand{\labelenumi}{(\alph{enumi})}
\item $K_{X \times_k \kappa}+D+D' \sim 0$, 
\item $D$ is regular, $D \cap D'$ is smooth over $\kappa$, $D \not\subset \Supp D'$, 
\item $H^1(X \times_k \kappa, \MO_{X \times_k \kappa}(-D))=0$, and 
\item $H^1(X \times_k \kappa, \MO_{X \times_k \kappa}(-D'))=0$. 
\end{enumerate}
Then $X$ is $2$-quasi-$F^e$-split for all $e>0$. 
\end{cor}

\begin{proof}
Set $X_{\kappa} \coloneqq X \times_k \kappa$. 
We show that $(X \times_k \kappa, D_1+\cdots+ D_{n-1})$ is weakly $2$-quasi-$F$-split (see \cite[Definition 5.13]{Kawakami-Tanaka(dPvar)} for the definition of weak quasi-$F$-splitting).

Since $D+D'$ is ample, $D + D'$ is connected, and hence $D'|_{D} = D \cap D' \neq \emptyset$. 
By \cite[Corollary 5.17]{Kawakami-Tanaka(dPvar)}, it is enough to show that $(D, D'|_D)$ is weakly $2$-quasi-$F$-split. 
Again by \cite[Corollary 5.17]{Kawakami-Tanaka(dPvar)}, 
it suffices to prove that 
\begin{enumerate}
\renewcommand{\labelenumi}{(\roman{enumi})}
    \item $D'|_D$ is a prime divisor,
    \item $H^1(D, \MO_D(-D'|_D))=0$, and 
    \item $D'|_D$ is $2$-quasi-$F$-split. 
\end{enumerate}
Consider an exact sequence
\[
0 \to \MO_{X_{\kappa}}(-D-D') \to \MO_{X_{\kappa}}(-D') \to \MO_{D}(-D'|_D) \to 0, 
\]
which induces the following one: 
\[
0\overset{{\rm (d)}}{=}
H^1(X_{\kappa}, \MO_{X_{\kappa}}(-D')) \to H^1(D,  \MO_{D}(-D'|_D)) \to 
H^2(X_{\kappa}, \MO_{X_{\kappa}}(-D-D')).  
\]
Then (ii) follows from 
\[
H^2(X_{\kappa}, \MO_{X_{\kappa}}(-D-D')) \simeq 
H^2(X_{\kappa}, \MO_{X_{\kappa}}(K_{X_{\kappa}})) \simeq H^2(X, \MO_X(K_X)) \otimes_k \kappa 
%\simeq H^1(X, \MO_X) \otimes_k \kappa 
=0. 
\]
By the induced exact sequence 
\[
H^0(D, \MO_D) \to H^0(D'|_D, \MO_{D'|_D}) \to H^1(D, \MO_D(-D'|_D)) \overset{{\rm (ii)}}{=}0, 
\]
$H^0(D'|_D, \MO_{D'|_D})$ is a field. {Since $D'|_D = D \cap D'$ is smooth}, we obtain (i).
Let us show (iii). 
Since $D'|_D$ is a smooth projective curve over $\kappa$ with $K_{D'|_D} \sim 0$, 
the base change $(D'|_D) \times_{\kappa} \overline \kappa$ 
to the algebraic closure $\overline \kappa$ of $\kappa$ is an elliptic curve, which is 2-quasi-$F$-split (\cite[Remark 2.11]{KTTWYY1}). %, 
%$D'|_D$ is smooth by [PW], i.e., $D'|_D$ is a smooth elliptic curve.  
Thus so is $D'|_D$ \cite[Proposition 2.12]{KTY}, and (iii) holds. 

Now, we show that $X$ is quasi-$F$-split.
Since $(X \times_k \kappa, D_1+\cdots+ D_{n-1})$ is weakly $2$-quasi-$F$-split
it follows that $X \times_k \kappa$ is $2$-quasi-$F$-split.
By \cite[Proposition 2.12]{KTY}, $X$ is $2$-quasi-$F$-split.
\end{proof}

\begin{rem}\label{r-IOA-Fano3}
The assumptions (c) and (d) automatically hold when $k=\kappa$ and $D'$ is connected. 
Indeed, we have an exact sequence $H^0(X, \MO_X) \xrightarrow{\simeq}H^0(E, \MO_E) \to H^1(X, \MO_X(-E)) \to H^1(X, \MO_X)=0$ for any reduced connected divisor $E$ on $X$.  
\end{rem}

Although the following result is contained in \cite[Theorem F]{KTTWYY1}, 
we include a proof for the reader's convenience, as its proof is quite short.

\begin{prop}\label{prop:criterion for qFs}
    Let $X$ be a smooth Fano threefold. % 3-fold.
    Suppose that the following condition.
    \begin{enumerate}
        \item $H^0(X, \Omega_X^2(p^iK_X))=0$ for all $i>0$.
        \item $H^1(X, \Omega_X^1(K_X))=0$.
        \item $H^2(X, \Omega_X^1(-p^iK_X))=0$ for all $i>0$.
    \end{enumerate}
    Then $X$ is quasi-$F$-split.
\end{prop}
\begin{proof}
It is enough to show that $H^2(X, B_n\Omega_X^1(p^nK_X))=0$ for some $n>0$.
    By \eqref{exact:B} and (2), it suffices to prove $H^2(X, Z_n\Omega^1_X(p^nK_X))=0$ for some $n>0$. 
    By (\ref{exact:Z}), this vanishing is reduced to 
    \[
    H^1(X, B_1\Omega^2_X(pK_X))=0\,\,\,\text{and}\,\,\,
    H^2(X, Z_{n-1}\Omega^1_X(p^nK_X))=0.
    \]
    Repeating this procedure, the vanishing $H^2(X, Z_n\Omega^1_X(p^nK_X))=0$ is reduced to  
    \[
    H^1(X, B_1\Omega^2_X(p^lK_X))=0\,\,\,\text{and}\,\,\,
    H^2(X, \Omega^1_X(p^nK_X))=0.
    \]
    for every $l\in\{1,\ldots, n\}$.
    Taking $n\gg0$, we may assume $H^2(X, \Omega^1_X(p^nK_X))=0$ by Serre vanishing. 
    By \eqref{exact:B} and (1), it suffices to show that $H^1(X, Z_1\Omega^2_X(p^lK_X))=0$ for every $l\in\{1,\ldots, n\}$.
    By  Serre duality, we have 
    \[
    H^1(X, \Omega_X^2(p^lK_X))\simeq H^2(X, \Omega_X^1(-p^lK_X)) \overset{{\rm (3)}}{=} 0
    \]
    for every $l>0$. 
    By (\ref{exact:Z}), the problem is finally reduced to $H^0(X, B_1\Omega_X^3(p^lK_X))=0$ for all $l\in\{1,\ldots, n\}$. 
    This holds because 
    $B_1\Omega_X^3(p^lK_X)\subset F_{*}\omega^{p^l+1}_X$ and $\omega^{-1}_X$ is ample. 
\end{proof}

\subsection{Weak del Pezzo surfaces}
In this subsection, we recall when canonical weal del Pezzo surfaces are $F$-split.
\begin{dfn}
Let $S$ be a normal Gorenstein projective surface. 
\begin{enumerate}
\item We say that $S$ is {\em weak del Pezzo} if $-K_S$ is nef and big. 
\item We say that $S$ is {\em del Pezzo} if $-K_S$ is ample. 
\end{enumerate}
\end{dfn}

\begin{thm}\label{t-dP-F-split}
Let $S$ be a canonical weak del Pezzo surface. 
Then $S$ is $F$-split if $p>5$ or $K_S^2 \geq 5$. 
\end{thm}
\begin{proof}
See \cite[Theorem A]{KT}. 
\end{proof}

\subsection{Akizuki-Nakano vanishing for hhypersurfaces} 

We will need the following variant of \cite[Lemma 2.8]{Kawakami-Tanaka(Lift1)}
when we deal with Fano threefolds $V_1$ and $V_2$ 
(i.e., the case of index two satisfying $(-K_X/2)^3  \in \{1, 2\}$).

\begin{lem}\label{lem:ANV for hypersurfaces}
Let $P$ be an lci projective variety of $\dim\,P=d+1$ over an algebraically closed field. 
Let $X$ be an ample effective Cartier divisor on $P$ 
which is a smooth variety. 
%    Let $X$ be a $d$-dimensional smooth complete intersection variety of ample Cartier divisors $H_1,\ldots, H_n$. 
%Let $i,j\in\Z_{\geq 0}$ be 
%Take non-negative integers $i$ and $j$. 
Assume that the equality 
    \begin{equation}\label{e1 ANV for CI}
    H^j(P,\Omega^{[i]}_P(-H))=0
    \end{equation}
 holds    if  $i+j<d+1$ and $H$ is an ample Cartier divisor on $P$. 
    Then the following hold. 
    \begin{enumerate}
        \item $H^j(X,\Omega^{i}_{X}(-H))=0$ if $i+j<d$ and $H$ is an ample Cartier divisor on $P$.
        \item $H^{j}(X, \Omega^{i}_{X}) \simeq H^j(P, \Omega^{i}_{P})$ 
        if   $i+j<d$.
%        \item $H^{j}(X, \Omega^{i}_{X})\simeq H^{j+n}(P, \Omega^{i+n}_{P})$ if $i+j>d$.
    \end{enumerate}
\end{lem}
\begin{proof}
    Note that $X$ is contained the smooth locus of $P$, 
    because $X$ is a smooth effective Cartier divisor on $P$. Then the assertion holds by 
    the same argument as in \cite[Lemma 2.8]{Kawakami-Tanaka(Lift1)}
    after replacing $\Omega^{i}_P$ by $\Omega^{[i]}_P$. 
\end{proof}

%% file: section3.tex
\section{Fano threefolds with $\rho \geq 4$}\label{sec:5}

\begin{prop}\label{p-rho6}
Let $X$ be a smooth Fano threefold with $\rho(X) \geq 6$. 
Then $X$ is quasi-$F$-split. 
\end{prop}

\begin{proof}
In this case, $X \simeq S \times \P^1$ for a smooth del Pezzo surface \cite[Section 7.6]{FanoIV}. 
Since $S$ is quasi-$F$-split and $\P^1$ is $F$-split, 
$S \times \P^1$ is quasi-$F$-split \cite[Proposition 6.7]{KTY}. 
\end{proof}

\begin{prop}\label{p-rho5}
Let $X$ be a smooth Fano threefold with $\rho(X)=5$. 
Then $X$ is $F$-split. 
\end{prop}

\begin{proof}
Recall that $X$ is of No.~5-1, 5-2, or 5-3 \cite[Section 7.5]{FanoIV}. 
If $X$ is 5-3, then $X \simeq S \times \P^1$ for a smooth del Pezzo surface $S$ with $\rho(S)=4$, and hence $F$-split (e.g.~$S$ is toric, and hence so is $X$). 

Assume that $X$ is 5-1. 
Then $X = \Bl_{B_1 \amalg B_2 \amalg B_3} Y$, 
where $Y \coloneqq \Bl_C Q$ is a blowup of $Q$ along a conic $C$ and 
$B_1, B_2, B_3$ are mutually distinct one-dimensional fibres of 
$\sigma\colon Y \coloneqq \Bl_C Q \to Q$ \cite[Section 7.5]{FanoIV}. 
Since the smallest linear subvariety $\langle C \rangle$ of $\P^4$ containing $C$ is a plane, 
we obtain $C \subset \langle C \rangle =\overline{H}_1 \cap \overline{H}_2$ for 
suitable hyperplanes $\overline{H}_1$ and $\overline{H}_2$ on $\P^4$. 
Set $H_i \coloneqq Q \cap \overline{H}_i$ for each $i \in \{1, 2\}$. 
We then get a scheme-theoretic equality $C = H_1 \cap H_2$. 
Note that each $H_i$ is a (possibly singular) quadric surface in $\overline{H}_i = \P^3$, which is smooth along $C = H_1 \cap H_2$. 
It is easy to see that $\Delta$ is effective for the divisor $\Delta$ defined by $h^*(K_Q+H_1+H_2) = K_X +\Delta$, where $h \colon X \to Q$ denotes the induced birational morphism. 
It is enough to show that $(Q, H_1 + H_2)$ is $F$-split (Proposition \ref{p-F-split-blowup}). 
Take a general hyperplane section $H_3$. 
Since $K_Q+H_1+H_2+H_3 \sim 0$ and 
all $H_3, H_3 \cap H_2$, and $H_3 \cap H_2 \cap H_1$ are smooth, 
$(Q, H_1+H_2+H_3)$ is $F$-split (Corollary \ref{c-P1-2pts-Fsplit}), and hence $(Q, H_1+H_2)$ is $F$-split.

\medskip

Assume that $X$ is No.~5-2. 
Then $X = \Bl_{B \amalg B'} Y$, 
where $Y \coloneqq \Bl_{L_1 \amalg L_2} \P^3$ is a blowup of $\P^3$ along a disjoint union of lines $L_1$ and $L_2$, 
and $B$ and $B'$ are mutually distinct one-dimensional fibres of 
$\sigma\colon Y =  \Bl_{L_1 \amalg L_2} \P^3 \to \P^3$ lying over $L_1$ \cite[Section 7.5]{FanoIV}. 
Take two planes $H_1$ and $H_2$ on $\P^3$ such that $H_1 \cap H_2 =L_1$. 
Note that each $H_i \cap L_2$ is a smooth point. 
Pick a plane $H_3$ containing $L_2$. 
Then it is easy to see that $\Delta$ is effective for the divisor $\Delta$ defined by $h^*(K_{\P^3}+H_1+H_2+H_3) = K_X +\Delta$. 
It is enough to show that $(\P^3, H_1 + H_2+H_3)$ is $F$-split (Proposition \ref{p-F-split-blowup}). 
Pick a general hyperplane $H_4$. 
Apply Corollary \ref{c-P1-2pts-Fsplit} by setting $D_1 \coloneqq H_1, D_2 \coloneqq H_2$, and $D_3 \coloneqq H_3 + H_4$.  
Then $(X, D_1+D_2+D_3)$ is $F$-split, and hence  $(\P^3, H_1 + H_2 +H_3)$ is $F$-split. 
\end{proof}

\begin{prop}\label{p-rho4}
Let $X$ be a smooth Fano threefold with $\rho(X)=4$. 
Then $X$ is $F$-split. 
\end{prop}

\begin{proof}
We treat the following six cases separately: 
\begin{enumerate}
\item 4-4, 4-10, 4-12. 
\item  4-3, 4-6, 4-8, 4-13. 
\item 4-5, 4-7, 4-9, 4-11. 
\item 4-11. 
\item 4-2. 
\item 4-1. 
\end{enumerate}

(1) If $X$ is 4-4, then there is a smooth curve $B$ on $X$ such that the blowup $\widetilde X$ of $X$ along $B$ is Fano \cite[Proposition 5.31]{FanoIV}. 
Since $\widetilde X$ is $F$-split (Proposition \ref{p-rho5}), so is $X$. 
If $X$ is 4-10, then we can write $X \simeq S \times \P^1$ for a smooth del Pezzo surface $S$ 
with $K_S^2 = 7$ \cite[Section 7.4]{FanoIV}. 
In this case,  $X$ is clearly $F$-split. 
Assume that $X$ is 4-12. 
Then we have $X = \Bl_{B \amalg B'} Y_{\text{2-33}}$, 
where $Y_{\text{2-33}} = \Bl_L \P^3$ is the blowup of $\P^3$ along a line $L$, and 
$B$ and $B'$ are mutually disjoint one-dimensional fibres of the induced blowup 
$\rho : Y_{\text{2-33}}= \Bl_L \P^3 \to \P^3$. 
In this case, we can apply a similar argument to that of 5-2 in the proof of Proposition \ref{p-rho5}.

(2) Assume that $X$ is one of 4-3, 4-6, 4-8, 4-13. 
In this case, 
there is a blowup $h: X \to \P^1_1 \times \P^1_2 \times \P^1_3$ along a curve  $B$ of tridegree $(1, 1, c)$ for some $c \geq 0$ \cite[Section 7.4]{FanoIV}. 
Let $B' \subset \P^1_1 \times \P^1_2$ be the image of $B$, which is a curve of bidegree $(1, 1)$. 
Then the induced morphism $B' \to \P^1$ to the first direct product factor is an isomorphism. 
Hence $B' \simeq \P^1$. 
Set $D \subset \P^1 \times \P^1 \times \P^1$ to be the inverse image of $B'$, 
which satisfies $D \in |\MO_{\P^1 \times \P^1 \times \P^1}(1, 1, 0)|$ and $D \simeq \P^1 \times \P^1$.  
Then it is easy to see that $\Delta$ is effective for the divisor $\Delta$ defined by $h^*(K_{\P^1 \times \P^1 \times \P^1}+D) = K_X +\Delta$. 
It is enough to show that $(\P^1 \times \P^1 \times \P^1, D)$ is $F$-split (Proposition \ref{p-F-split-blowup}), which follows from the fact that $D$ is $F$-split 
and $-(K_{\P^1 \times \P^1 \times \P^1} +D)$ is ample (Corollary \ref{c-F-split-IOA}).

(3) Assume that $X$ is one of  4-5,  4-7, 4-9. 
Let $\tau : \F_1 \to \P^2$ be the blowdown of the $(-1)$-curve 
$\Gamma$ and set $P :=\tau(\Gamma)$. 
Then there is a blowup $f: X \to Y = \F_1 \times \P^1$ 
along a smooth curve $B$ lying over $B' \subset \F_1$, 
where $B'$ is disjoint from the $(-1)$-curve $\Gamma$ and 
$B'$ is the inverse image of a line $L \subset \P^2$ \cite[Section 7.4]{FanoIV}. 
Then the induced composite  morphism 
\[
h : X \to \F_1 \times \P^1 \to \P^2 \times \P^1
\]
is the blowup along $(P \times \P^1) \amalg \overline{B}$, 
where $\overline{B} \subset \P^2 \times \P^1$ denotes the image of $B \subset \F_1 \times \P^1$. 
In particular, we get $(P \times \P^1) \amalg \overline{B} \subset 
L' \times \P^1 \cup L \times \P^1$ for a line $L'$  passing through $P$. 
Then it is easy to see that $\Delta$ is effective for the divisor $\Delta$ defined by $h^*(K_{\P^2 \times \P^1}+L \times \P^1 + L' \times \P^1) = K_X +\Delta$. 
It is enough to show that $(\P^2 \times \P^1, L \times \P^1 + L' \times \P^1)$ is $F$-split (Proposition \ref{p-F-split-blowup}). 
This holds, because we may assume that $L \times \P^1 + L' \times \P^1$ is a torus-invariant reduced divisor on a toric variety $\P^2 \times \P^1$.

(4) 
Assume that $X$ is 4-11. 
Then there is a blowup $f \colon X \to \F_1 \times \P^1$ along $C$, 
where $C = \Gamma \times \{ t\}$ for the $(-1)$-curve $\Gamma$ on $\F_1$ and a closed point $t$ of $\P^1$ \cite[Section 7.4]{FanoIV}. 
For the blowup $\tau \times {\rm id} \colon \F_1 \times \P^1 \to \P^2 \times \P^1$, consider the composite birational morphism: 
\[
h \colon X \xrightarrow{f} \F_1 \times \P^1 \xrightarrow{\tau \times {\rm id}} \P^2 \times \P^1. 
\]
For the blowup centre $P \coloneqq \tau(\Gamma)$ of $\tau \colon \F_1 \to \P^2$, 
pick two lines $L$ and $L'$ on $\P^2$ passing through $P$. 
Then we have \[K_Y + (\Gamma \times \P^1) + L_Y+L'_Y= (\tau \times {\rm id})^*(K_{\P^2 \times \P^1} + (L \times \P^1) + (L' \times \P^1)),\] 
where $L_Y$ and $L'_Y$ denote the proper transforms of 
$L \times \P^1$ and $L' \times \P^1$, respectively. 
Since $C=\Gamma \times \{ t\}$ is contained in $\Gamma \times \P^1$, we see that
the divisor $\Delta$ 
defined by $K_X + \Delta = h^*(K_{\P^2 \times \P^1} + (L \times \P^1) + (L' \times \P^1))$ is effective. 
Then $X$ is $F$-split, because 
$(\P^2 \times \P^1, L \times \P^1 + L' \times \P^1)$ is $F$-split (Proposition \ref{p-F-split-blowup}).

(5) 
Assume that $X$ is 4-2. 
Then $X$ is a blowup of $Y=\P_{\P^1 \times \P^1}(\MO \oplus \MO(1, 1))$ 
along an elliptic curve $B$ on a section $S$ of the $\P^1$-bundle $\pi:Y \to \P^1 \times \P^1$ disjoint from the negative section $S'$ of $\pi$ \cite[Proposition 5.28]{FanoIV}. 
Note that $S \simeq S' \simeq \P^1 \times \P^1$. 
%For the induced projection $\pi : Y =\P_{\P^1 \times \P^1}(\MO \oplus \MO(1, 1)) \to \P^1 \times \P^1$, 
We have that $K_Y+S+S' \sim \pi^*K_{\P^1 \times \P^1}$. 
{Indeed, since $K_Y+S+S'$ is $\pi$-numerically trivial, we can write $K_Y+S+S'\sim \pi^{*}D$ for some Cartier divisor $D$ on $\P^1 \times \P^1$. By restricting to $S$, we obtain $D\sim K_{\P^1 \times \P^1}$.} 

Then, the $\Q$-divisor 
\[
-(K_Y+S+(1-\epsilon)S') = \epsilon S' +\pi^*K_{\P^1 \times \P^1}
\]
is ample for $0< \epsilon \ll 1$. 
After perturbing $\epsilon$, 
the problem is reduced to the case when $(p^e-1)(K_Y+S+(1-\epsilon)S')$ is Cartier for some $e>0$. 
Replacing $e$ by some $e' \in e\Z_{>0}$,  
we may assume that $H^1(X, \MO_X(-S-(p^e-1)(K_Y+S+(1-\epsilon)S')))=0$. 
Since $(S, (1-\epsilon)S'|_S)  = (\P^1 \times \P^1, 0)$ is sharply $F$-split, 
$(X, S+(1-\epsilon)S')$ is sharply $F$-split (Proposition \ref{p-F-split-IOA}). 
Hence $X$ is $F$-split.

(6) Assume that $X$ is 4-1. 
Then $X$ is a prime divisor on $\P^1_1 \times \P^1_2 \times \P^1_3 \times \P^1_4$
of multi-degree $(1, 1, 1, 1)$ \cite[Section 7.4]{FanoIV}. 
For each $i \in \{1, 2, 3, 4\}$, we set $H_i := \pi_i^*\MO_{\P^1_i}(1)$ and $H'_i := \pr_i^*\MO_{\P^1_i}(1)$, where $\pi_i$ and $\pr_i$ denote the the induced morphisms: 
\[
\pi_i :X \hookrightarrow \P^1_1 \times \P^1_2 \times \P^1_3 \times \P^1_4 \xrightarrow{\pr_i} \P^1_i. 
\]
It holds that $-K_X \sim H_1 + H_2 + H_3 + H_4$. 
Note that 
\[
H_1 \cdot H_2 \cdot H_3 = H_1' \cdot H'_2 \cdot H'_3 \cdot (H'_1 + H'_2 + H'_3 + H'_4) = H'_1 \cdot H'_2 \cdot H'_3 \cdot H'_4 =1. 
\]
%By symmetry, $H_1 \cdot H_2 \cdot H_4=1$. 
Take the generic members $H_1^{\gen}, H_2^{\gen}, H_3^{\gen}, H_4^{\gen}$ of $H_1, H_2, H_3, H_4$, where each $H_i^{\gen}$ is an effective Cartier divisor on $X \times_k \kappa$ for suitable 
purely transcendental field extension $\kappa/k$ (Remark \ref{r-generic-snc}). 
Set $D_1 := H_1^{\gen}, D_2 := H_2^{\gen}, D_3 := H_3^{\gen} +H_4^{\gen}$. 
Then $D_1^{\gen} \cap D_2^{\gen}$ is smooth, because $\deg (H_3^{\gen}|_{D_1 \cap D_2})=1$ (Remark \ref{r-g0-to-P1}). 
%$(D_1 \cap D_2, (D_3+D_4)|_{D_1 \cap D_2})$
By Corollary \ref{c-P1-2pts-Fsplit}, $X$ is $F$-split. 
\end{proof}

%% file: section4.tex
\section{Fano threefolds with $\rho =3$ (except for 3-10)}

The purpose of this subsection is to prove that 
an arbitrary smooth Fano threefold $X$ with $\rho(X)=3$ is $F$-split 
except for 3-10 (Proposition \ref{p-rho3}). 
We start with some complicated cases: 3-1, 3-3 and 3-4.

\begin{lem}\label{l-3-1}
Let $X$ be a smooth Fano threefold of No.~3-1. 
Then the following hold. 
\begin{enumerate}
    \item 
Let $\varphi \colon X \to \P^1 \times \P^1 \times \P^1$ be a finite double cover 
such that $(\varphi_*\MO_X/\MO_Y)^{-1} \simeq \MO_{\P^1 \times \P^1 \times \P^1}(1, 1, 1)$. 
Then $\varphi$ is separable, i.e., the induced field extension $K(X)/K(Y)$ is separable. 
    \item 
    $X$ is $F$-split. 
\end{enumerate}
\end{lem}

\begin{proof}
Let us show (1). 
Set $Y\coloneqq \P^1 \times \P^1 \times \P^1$ and  $\cL \coloneqq \MO_{\P^1 \times \P^1 \times \P^1}(1, 1, 1)$. 
By \cite[Proposition 0.1.2]{CD89} and \cite[Lemma 3.2]{Ful98}, 
it suffices to show that 
\[
\deg c_3(\Omega_Y^1 \otimes \cL^{\otimes 2}) \neq 0.
\]
Set $\cM_i \coloneqq \pr_i^*\Omega_{\P^1}^1 \otimes \cL^{\otimes 2}$ for every $i \in \{1, 2, 3\}$, i.e., 
\[
\cM_1 \coloneqq \MO(0, 2, 2), \qquad 
\cM_2 \coloneqq \MO(2, 0, 2), \qquad 
\cM_3 \coloneqq \MO(2, 2, 0). \qquad 
\]
We have $\Omega_Y^1 \otimes \cL^{\otimes 2}= \cM_1 \oplus \cM_2 \oplus \cM_3$. 
Then the following holds (cf. \cite[Appendix A, Section 3]{Har77}): 
\begin{eqnarray*}
\sum_{i=0}^{\infty} c_i(\Omega_Y^1 \otimes \cL^{\otimes 2})  t^i
&=& c_t(\Omega_Y^1 \otimes \cL^{\otimes 2}) \\
&=& c_t(\cM_1) c_t(\cM_2)c_t(\cM_3)\\
&=&(1+c_1(\cM_1)t) (1+c_1(\cM_2)t) (1+c_1(\cM_3)t). 
\end{eqnarray*}
Therefore, we get 
\[
\deg c_3(\Omega_Y^1 \otimes \cL^{\otimes 2}) =\deg(c_1(\cM_1)c_1(\cM_2)c_1(\cM_3)) = \cM_1 \cdot \cM_2 \cdot \cM_3 >0, 
%=(2H_2+3H_3)\cdot (2H_1+2H_3) \cdot (2H_1+2H_2) >0,  
\]
as required. Thus (1) holds. 

Let us show (2). 
There is a finite separable double cover $\varphi \colon X \to Y = \P^1 \times \P^1 \times \P^1$ as in (1) \cite[Theorem 6.7]{FanoIII}. 
Moreover, for each $i \in \{1, 2, 3\}$, 
the composition 
$\varphi_i \colon X \xrightarrow{\varphi} \P^1 \times \P^1 \times \P^1 \xrightarrow{\pr_i} \P^1$ is the contraction of an extremal ray \cite[Remark 6.8]{FanoIII}. 
Since $\varphi$ is separable, if $H_1, H_2, H_3$ are general members, 
then the scheme-theoretic intersection $H_1 \cap H_2 \cap H_3$ is reduced two points. 
Let $H_1^{\gen}, H_2^{\gen}, H_3^{\gen}$ be their  generic members. 
Then the regular curve $H_1^{\gen} \cap H_2^{\gen}$ is automatically smooth (Remark \ref{r-g0-to-P1}). 
Then $X$ {is $F$-split} (Corollary \ref{c-P1-2pts-Fsplit}, Lemma \ref{l-P1-2pts-Fsplit}).
%and hence so is $X$. 
\end{proof}

\begin{lem}\label{l-3-4-Y}
Let $Y$ be a smooth Fano threefold of No.~2-18. 
Let $f_1 \colon Y \to \P^1$ and $f_2 \colon Y \to \P^2$ be the contractions of the extremal rays. 
Take a general point $P \in \P^1$ and a general line $L$ on $\P^2$. 
Then $H^0(\Gamma, \MO_{\Gamma}) = k$ 
for the scheme-theoretic intersection $\Gamma \coloneqq f_1^{-1}(P) \cap f_2^{-1}(L)$. 
\end{lem}

\begin{proof} 
Set $S_1 \coloneqq f_1^{-1}(P)$ and $S_2 \coloneqq f_2^{-1}(L)$. 
By \cite[Theorem 15.2]{FS20} and \cite[Theorem 3.3]{BT22}, 
$S_1$ is a canonical del Pezzo surface. 
In particular, $S_1$ is a rational surface. 
Since $-K_Y \sim S_1 + 2S_2$ \cite[Proposition 5.9 and Section 9.2]{FanoIII}, the divisor 
$-K_Y|_{Y_K}=2S_2|_{Y_K}$ is ample 
for the generic fibre $Y_K$ of $f_1 \colon Y \to \P^1$, 
where $K \coloneqq K(\P^1)$. 
Hence $S_2|_{S_1} (=\Gamma)$ is ample, 
as $S_1$ is chosen to be a general fibre of $f_1$. 
Therefore, $H^1(S_1, \MO_{S_1}(-\Gamma))=0$ by \cite[Proposition 3.3]{CT19} (or \cite[Theorem 3]{Muk13}), 
and hence $k= H^0(S_1, \MO_{S_1}) \xrightarrow{\simeq} H^0(\Gamma, \MO_{\Gamma})$. 
\end{proof}

\begin{lem}\label{l-3-3-3-4}
Let $X$ be a smooth Fano threefold of No.~3-3 or 3-4. 
Then $X$ is $F$-split. 
\end{lem}

\begin{proof}
For each case, $X$ has exactly three extremal rays and 
there is a conic bundle structure $f\colon X \to \P^1 \times \P^1$ \cite[Section 7.3]{FanoIV}. 
In what follows, we shall use their properties obtained in \cite[Propositions 4.33 and 4.35]{FanoIV}. 
For each $i \in \{1, 2\}$, let 
\[
\varphi_i \colon X \xrightarrow{f} \P^1 \times \P^1 \xrightarrow{\pr_i} \P^1 =:Z_i
\]
be the induced composite contraction. 
Note that each $\varphi_i$ corresponds to %a contraction of 
a two-dimensional extremal face of $\NE(X)$. 
Let $\varphi_3 \colon X \to Z_3 = \P^2$ be the contraction  of the remaining two-dimensional extremal face of $\NE(X)$. 
For each $i \in \{1, 2, 3\}$, let $H_i$ be the pullback of the ample generator on $Z_i$. 
We recall $-K_X \sim H_1 +H_2 +H_3$ (\cite[Propositions 4.33 and 4.35]{FanoIV}). 
In particular, 
\[
2 = -K_X \cdot H_1 \cdot H_2 =(H_1 + H_2 + H_3) \cdot H_1 \cdot H_2 = H_1 \cdot H_2 \cdot H_3. 
\]

For each $i \in \{1, 2, 3\}$, $H'_i$ denotes the generic member of $H_i$, 
which is a regular prime divisor on $X' \coloneqq X \times_k k'$ for a suitable purely transcendental extension $k'/k$ (Remark \ref{r-generic-snc}). 
Note that 
$H'_1 \cap H'_2$ is a smooth curve, because $f$ is not wild \cite[Corollary 8 and Remark 9]{MS03}.

We now finish the proof by assuming that 
\begin{enumerate}
\item[(i)] $H^1(X', K_{X'}+H'_1 + 2H'_3) = H^1(X', K_{X'}+H'_2 + 2H'_3)=0$, and 
\item[(ii)] $H^2(X', K_{X'} +2H'_3) = 0$. 
\end{enumerate}
By Corollary \ref{c-P1-2pts-Fsplit} and Lemma \ref{l-P1-2pts-Fsplit}, 
it suffices to find $H''_3 \in |H'_3|$ such that 
$H'_1 \cap H'_2 \cap H''_3$ is smooth and zero-dimensional. 
Since $k'$ is an infinite field and $H'_1 \cap H'_2$ is isomorphic to a smooth conic on $\P^2_{k'}$, 
the Bertini theorem enables us to find a smooth zero-dimensional effective Cartier divisor 
$D$ on $H'_1 \cap H'_2$ such that $H'_3|_{H'_1 \cap H'_2} \sim D$.
Therefore, it is enough to show that the restriction maps 
\[
H^0(X', \MO_X(H'_3)) \xrightarrow{\alpha} H^0(H'_1, \MO_X(H'_3)|_{H'_1}) 
\xrightarrow{\beta} H^0(H'_1 \cap H'_2, \MO_X(H'_3)|_{H'_1 \cap H'_2}) 
\]
are surjective. 
The restriction map $\alpha$
is surjective, because 
\[
H^1(X', H'_3-H'_1) \simeq H^1(X', K_{X'}+H'_2 +2H'_3) \overset{{\rm (i)}}{=} 0. 
\]
The problem is reduced to the surjectivity of $\beta$. 
To this end, it suffices to prove 
$H^1(H'_1, \MO_{X'}(-H'_2 +H'_3)|_{H'_1}) =0$. 
We have an exact sequence 
\[
H^1(X', -H'_2+H'_3) \to H^1(H'_1, \MO_{X'}(-H'_2 +H'_3)|_{H'_1}) \to H^2(X', -H'_1-H'_2 +H'_3). 
\]
By 
\[
-H'_2 + H'_3  \sim K_{X'}+H'_1 +2H'_3 \qquad {\rm and} \qquad 
-H'_1 -H'_2 + H'_3 \sim K_{X'} + 2H'_3, 
\]
(i) and (ii) imply $H^1(H'_1, (- H'_2+H'_3)|_{H'_1}) =0$. 
Therefore, it is enough to prove (i) and (ii). 
\begin{claim*}
The following hold. 
\begin{enumerate}
\item $H'_i$ is a regular weak del Pezzo surface for every $i \in \{1, 2, 3\}$. 
\item $H^0(H'_i, \MO_{H'_i}) = k'$ for every $i \in \{1, 2, 3\}$. 
\item $H^2(X', K_{X'}+H'_i)=0$ for every $i \in \{1, 2, 3\}$. 
\item $H^2(X', K_{X'}+H'_i+H'_j)= H^2(X', K_{X'}+2H'_3)=0$ for $1 \leq i < j \leq 3$. 
%\item $H^2(H_i, $
\item $H^1(H'_i, \MO_{H'_i}) =0$ and $H^1(H'_i, K_{H'_i})=0$ for every $i \in \{1, 2, 3\}$. 
\item $H^0(H'_1 \cap H'_3, \MO_{H'_1 \cap H'_3}) = k'$ and $H^0(H'_2 \cap H'_3, \MO_{H'_2 \cap H'_3}) = k'$. 
\item $H^1(H'_3, K_{H'_3}+(H'_1+H'_3)|_{H'_3})=0$
\end{enumerate}
\end{claim*}

\begin{proof}[Proof of Claim]
Since $H'_i$ is the generic member of a base point free linear system $|H_i|$, 
$H'_i$ is a regular prime divisor on $X'$. 
If $i \in \{1, 2\}$, then $-K_{H'_i}$ is ample, because $H'_i$ is the generic fibre of 
$\varphi_i \colon X \to \P^1$. 
We see that $-K_{H'_3}$ is nef and big by $-K_{H'_3} \sim (H'_1+H'_2)|_{H'_3}$ and 
$(H'_1+H'_2)^2 \cdot H'_3 = 2H'_1 \cdot H'_2 \cdot H'_3 >0$. 
Thus (1) holds. 
If $i \in \{1, 2\}$, then (2) holds by the fact that $H_i$ is 
(a base change of) the generic fibre of a contraction $\varphi_i \colon X \to \P^1$. 
We have $H^0(H'_3, \MO_{H'_3})=k'$, because general members of the complete linear system $|H_3|$ are geometrically integral \cite[Proposition 2.10]{FanoI}. Thus, (2) holds. 

Let us show (3). 
Consider an exact sequence 
\begin{multline*}
    0 = H^2(X', K_{X'}) \to H^2(X', K_{X'}+H'_i) \to H^2(H'_i, K_{H'_i})\\  \to H^3(X', K_{X'}) \to H^3(X', K_{X'}+H'_i)=0.
\end{multline*}
By (2) and Serre duality, we obtain 
$h^2(H'_i, K_{H'_i})=1$ and $h^3(X', K_{X'}) =1$. 
Therefore, $H^2(X', K_{X'}+H'_i)=0$. 
Thus (3) holds. 

Let us show (4). 
If $i\neq j$, then $H^2(X', K_{X'}+H'_i+H'_j)=0$ by an exact sequence 
\[
0 = H^2(X', K_{X'}+H'_i) \to H^2(X', K_{X'}+H'_i+H'_j) \to H^2(H'_j, K_{H'_j} +H'_i)=0,  
\]
where $H^2(H'_j, K_{H'_j} +H'_i)=0$ follows from Serre duality and $H^0(H'_j, \MO_{X'}(-H'_i)|_{H'_j})=0$. 
Similarly, $H^2(X', K_{X'}+2H'_3)=0$  by $H^0(H'_3, \MO_{X'}(-H'_3)|_{H'_3})=0$. 
Thus (4) holds.

Let us show (5). By an exact sequence 
\[
0= H^1(X', \MO_{X'}) \to H^1(H'_i, \MO_{H'_i}) \to H^2(X', \MO_{X'}(-H'_i)),   
\]
it suffices to prove $H^2(X', \MO_{X'}(-H'_i))=0$, 
which follows from (4) by using $-H'_i \sim K_{X'} +H'_1+H'_2+H'_3 -H'_i$. 
Thus (5) holds.

Let us show (6). 
Since $X$ has exactly three extremal rays, 
there is the extremal ray $R$ such that $R$ is the intersection of the extremal faces corresponding to $X \to \P^1$ and $X \to \P^2$. 
Let $f\colon X \to Y$ be the contraction of $R$. 
%the contraction $f\colon X \to Y$ of 
If $X$ is 3-3, then $Y=\P^1 \times \P^2$ and $f\colon X \to Y \coloneqq \P^1 \times \P^2$ and $H'_1 \cap H'_3$ is isomorphic to the corresponding intersection $H_1'^Y \cap H_3'^Y$ on $Y$, because $H_1'^Y \cap H_3'^Y$ is disjoint from $f(\Ex(f))$. 
Therefore, $H'_1 \cap H'_3$ is geometrically integral. 
If $X$ is 3-4, then $Y = \F_1$ or $Y$ is a smooth Fano threefold of No.~2-18. 
If $Y = \F_1$, then $H'_1\cap H'_3$ is a smooth fibre of  $f\colon X \to \F_1$. 
The other case follows from  Lemma \ref{l-3-4-Y} by using the upper semi-continuity \cite[Ch. III, Theorem 12.8]{Har77}. 
%, which is a smooth rational curve. 
Thus (6) holds. 

Let us show (7). 
We have 
\[
H_1 \cdot H^2_3 = H_2 \cdot H^2_3 =1, 
\]
because $2 = -K_X \cdot H^2_3 = (H_1 + H_2) \cdot H_3^2$ and a fibre $\zeta \equiv H_3^2$ of $\varphi_3$ is not contracted by $\varphi_i$ for each $i \in \{1, 2\}$.
%In particular, $H'_1 \cap H'_3 \cap H''_3 = k'$. Then $H^1(H_3, -D)=0$ by CT-proof 
Fix a general member $H''_3$ of $|H'_3|$. 
By Serre duality, it is enough to show $H^1(H'_3, \MO_{H'_3}(-D))=0$
for an effective Cartier divisor $D \coloneqq (H'_1 + H''_3)|_{H'_3}$ on $H'_3$. 
By (2), the problem is reduced to $H^0(D, \MO_D)=k'$. 
Clearly, $D$ is nef. 
It holds that $D^2 = (H_1 +H_3)^2 \cdot H_3 \geq H_1 \cdot H^2_3 =1>0$.  
Hence $D$ is nef and big. 
Then $H^0(D, \MO_D)$ is a field \cite[Corollary 3.17]{Eno}. 
Since $\Supp D$ contains $H'_1 \cap H'_3 \cap H''_3$ which is a $k'$-rational point 
by $H'_1 \cdot H'_3 \cdot H''_3 = H_1 \cdot H_3^2 =1$, we obtain field extensions 
\[
k' = H^0(H'_3, \MO_{H'_3}) \hookrightarrow H^0(D, \MO_D) \hookrightarrow 
H^0(H'_1 \cap H'_3 \cap H''_3, \MO_{H'_1 \cap H'_3 \cap H''_3}) =k', 
\]
which implies $H^0(D, \MO_D) =k'$, as required. 
This completes the proof of Claim.  
\end{proof}

It is enough to show (i) and (ii). 
As (ii) has been settled by Claim(4), 
let us  show (i). 
By $H^0(H'_3, \MO_{H'_3})=k'$ and $H^0(H'_1 \cap H'_3, \MO_{H'_1 \cap H'_3})=k'$ (Claim(2)(6)), 
we get  $H^3(H'_3, \MO_{H'_3}(-H'_1|_{H'_3}))=0$ by the following exact sequence:
\[
H^0(H'_3, \MO_{H'_3}) \xrightarrow{\simeq} H^0(H'_1 \cap H'_3, \MO_{H'_1 \cap H'_3}) \to 
H^1(H'_3, \MO_{H'_3}(-H'_1|_{H'_3})) \to 
H^1(H'_3, \MO_{H'_3})=0. 
\]
By Serre duality, we get  
\[
H^1(H'_3, K_{H'_3} + H'_1|_{H'_3})=0. 
\]
We have the following exact sequences:
\[
0=H^1(X', K_{X'}) \to H^1(X', K_{X'}+H'_1) \to H^1(H'_1, K_{H'_1})=0, 
\]
\[
0= H^1(X', K_{X'}+H'_1) \to H^1(X', K_{X'}+H'_1+H'_3) \to H^1(H'_3, K_{H'_3} + H'_1|_{H'_3})=0
\]
\begin{multline*}
    0= H^1(X', K_{X'}+H'_1+H'_3) \to H^1(X', K_{X'}+H'_1+2H'_3)\\ \to H^1(H'_3, K_{H'_3}+(H'_1+H'_3)|_{H'_3})=0, 
\end{multline*}
where the last equality follows from Claim (7). 
Thus $H^1(X', K_{X'}+H'_1 +2H'_3)=0$. Similarly, we obtain $H^1(X', K_{X'}+H'_2+2H'_3)=0$. 
Thus (i) holds. 
\end{proof}

\begin{prop}\label{p-rho3}
Let $X$ be a smooth Fano threefold with $\rho(X)=3$. 
\begin{enumerate}
    \item If $X$ is not 3-10, then $X$ is $F$-split. 
    \item Assume that $X$ is 3-10. Then $X$ is $F$-split if and only if 
    $X$ has no wild conic bunlde structure. 
\end{enumerate}
\end{prop}

\begin{proof}
We may assume that $X$ is none of 3-1, 3-3, and 3-4 (Lemma \ref{l-3-1}, Lemma \ref{l-3-3-3-4}). 
Note that if $X$ has a wild conic bundle structure, then $X$ is not $F$-split 
\cite[Theorem 2.1 or Corollary 2.5]{GLP15}. 
In what follows, we assume that any conic bundle structure from $X$ is generically smooth. 
Under this additional assumption, it suffices to show that $X$ is $F$-split. 
We divide the proof into the following four cases. 

\begin{enumerate}
\item 3-27, 3-28, 3-31. 
\item 3-5, 3-8, 3-12, 3-13, 3-15, 3-16, 3-17, 3-19,  
3-20, 3-21, 3-22, 3-23, 3-24, 3-26, 3-29. 
\item  3-6, 3-10, 3-18, 3-25. 
\item 3-2, 3-7, 3-9, 3-11, 3-14, 3-30. 
\end{enumerate}

(1) In this case, $X$ is toric \cite[Subsection 7.3]{FanoIV}, and hence $F$-split. 

(2) In this case, 
there exist a smooth Fano threefold $Y$ with $\rho(Y)=2$, 
a $\P^1$-bundle $g\colon Y \to \P^2$, and a subsection $B \subset Y$ of $g$ such that 
$X \simeq \Bl_B Y$, 
$B_{\P^2} \coloneqq g(B) \simeq \P^1$, and 
$Y$ is one of 2-32, 2-34, 2-35 \cite[Subsection 7.3, cf.~Theorem 4.23]{FanoIV}. 
Set $S \coloneqq g^{-1}(B)$, which is a $\P^1$-bundle over $\P^1$, and hence 
$F$-split. 
{Let $T$ be the pullback of the ample generator by the contraction of the other extremal ray $R$. 
By \cite[Proposition 5.9(3)]{FanoIII}, we can write $-K_Y \sim aS+2T$, where $a$ is the length of $R$.
Then we can check that $a>1$ in each case \cite[Subsection 7.2]{FanoIV}.
Thus $-(K_Y+S) \sim (a-1)S+bT$ is ample by Kleiman's criterion.} 
Therefore, $(Y, S)$ is $F$-split (Proposition \ref{p-F-split-IOA}), and hence $X$ is $F$-split (Proposition \ref{p-F-split-blowup}). 

(3) 
In this case, we can write $X \simeq \Bl_{C \amalg C'} \P^3$ or $X \simeq \Bl_{C \amalg C'} Q$ 
for a disjoint union of smooth curves $C$ and $C'$ on $\P^3$ or $Q$ \cite[Subsection 7.3]{FanoIV}. 
We only treat the case when $X$ is 3-6, as the other cases are similar. 
In this case, {$X \simeq \Bl_{C \amalg C'} \P^3$}, $C$ is a line, and 
we can write $C' = S_1 \cap S_2$ for some quadric surfaces $S_1$ and $S_2$.
Take a general plane $H$ containing the line $C$ and 
a general quadric surface $S$ containing $C'$. 
Let $H_X$ and $S_X$ be the proper transforms  of $H$ and $S$, respectively. 
Although $H_X \to H$ and $S_X \to S$ are not necessarily isomorphisms, 
these birational morphisms are isomorphic over $H \cap S$. 
Therefore, we obtain $H_X \cap S_X \xrightarrow{\simeq} H \cap S$. 
This is nothing but a general fibre of the contraction $X \to \P^1 \times \P^1$ \cite[Proposition 4.37]{FanoIV}. 
Since this is not a wild conic bundle, 
we get $H \cap S \simeq H_X \cap S_X \simeq \P^1$. 
Then $(\P^3, H+S)$ is $F$-split by applying Proposition \ref{p-F-split-IOA} twice. 
Hence $X$ is $F$-split (Proposition \ref{p-F-split-blowup}). 

(4) In what follows, we treat the remaining  cases separately.

\underline{3-2}:
%cf the proof of \ref{p-rho3}. 
We use the same notation as in \cite[Proposition 4.32]{FanoIV}. 
We have $-K_X \sim 2H_1 + H_2 +D$ and a conic bundle $f:X \to \P^1 \times \P^1$, 
where $D \simeq \P^1 \times \P^1$, and 
each $H_i$ is the pullback of $\MO_{\P^1}(1)$ by $X \xrightarrow{f} \P^1 \times \P^1 \xrightarrow{\pr_i} \P^1$. 
Moreover $f|_D \colon D =\P^1 \times \P^1 \to \P^1 \times \P^1$ 
is a finite double cover which can be written 
as ${\rm id} \times \psi$ for some double cover $\psi\colon\P^1 \to \P^1$. 
For the generic member $H^{\gen}_2$ of $|H_2|$, the intersection
$C \coloneqq H_2^{\gen} \cap D$ is a regular curve of genus zero
with $K_C+2H_1|_C \sim 0$, {because $H_1\cdot C=H_1\cdot H_2 \cdot D>0$}. 
By $\deg (H_1|_C)=1$, we get $C \simeq \P^1_{\kappa}$. 
Hence $X$ is $F$-split (Corollary \ref{c-P1-2pts-Fsplit}).

\underline{3-7}: 
In this case, 
$f\colon X \to W$ is a blowup along an elliptic curve $B = S_1 \cap S_2$ with $S_1, S_2 \in |-(1/2)K_W|$ \cite[Subsection 7.3]{FanoIV}.  
Let $S$ be a general member of $|-(1/2)K_W|$  containing $B$. 
Since $B$ is an ample effective Cartier divisor on $S$, it follows that
$S$ is smooth along $B$ and $S$ is normal. 
Note that the proper transform $S_X$ of $S$ on $X$ is a fibre of a contraction $\pi \colon X \to \P^1$. 
Therefore, the geometric generic fibre $X_{\overline K}$ of $\pi$ 
is normal, where $K \coloneqq K(\P^1)$ and $\overline K$ denotes the algebraic closure of $K$. 
Since $X_K$ is a regular del Pezzo surface, 
$X_{\overline K}$ has at worst canonical singularities \cite[Theorem 3.3]{BT22}. 
Therefore, a general fibre $S_X$ of $\pi \colon X \to \P^1$ is a  canonical del Pezzo surface. 
%Then a general member $S$ is canonical. 
By $K_S^2 =6$, Theorem \ref{t-dP-F-split} shows that $S (\simeq S_X)$ is $F$-split. 
Since $-(K_W+S)$ is ample, we have $(W, S)$ is $F$-split (Proposition \ref{p-F-split-IOA}), and hence $X$ is $F$-split (Proposition \ref{p-F-split-blowup}).

\underline{3-9}: 
By \cite[Proposition 4.42]{FanoIV},   there is a blowup
$f\colon X \to Y=Y_{\text{2-36}} = \P_{\P^2}(\MO \oplus \MO(2))$ 
along a smooth curve $B$ such that 
\begin{itemize}
\item $B$ is contained in a section $S$ of 
the $\P^1$-bundle structure $\pi \colon \P_{\P^2}(\MO \oplus \MO(2)) \to \P^2$, and 
\item $S$ is disjoint from another section $T$ of $\pi \colon\P_{\P^2}(\MO \oplus \MO(2)) \to \P^2$. 
%which can be contracted to one point. 
\end{itemize}
Since $\pi$ is $\P^1$-bundle, we have $K_Y+S+T = \pi^*K_{\P^2}$. 
Since $-(K_Y+S+ (1-\epsilon)T) = -\pi^*K_{\P^2} +\epsilon T$ is ample 
for some $0 < \epsilon \ll 1$ and $(S, 0)$ is $F$-split, 
$(Y, S+(1-\epsilon)T)$ is $F$-split (Proposition \ref{p-F-split-IOA}). 
Hence $X$ is $F$-split (Proposition \ref{p-F-split-blowup}). 

\underline{3-11}: 
By \cite[Subsection 7.3]{FanoIV}, there exists a blowup $f \colon X \to V_7$ 
along an elliptic curve $B = S_1 \cap S_2$ with $S_1, S_2 \in |-(1/2)K_{V_7}|$. 
By the same argument as in that of 3-7, a general member $S$ 
of $|(-1/2)K_{V_7}|$ is a canonical del Pezzo surface with $K_S^2 =7$. 
Then $S$ is $F$-split  (Theorem \ref{t-dP-F-split}). 
Since $-(K_{V_7}+S)$ is ample, $(V_7, S)$ is $F$-split (Proposition \ref{p-F-split-IOA}), and hence $X$ is $F$-split (Proposition \ref{p-F-split-blowup}).

\underline{3-14}: 
By \cite[Subsection 7.3]{FanoIV}, 
we have $X = \Bl_{P \amalg C} \P^3$, 
where $C$ is a smooth cubic curve contained in a plane $H$ and $P$ is a point satisfying $P \not\in H$. 
Take two {distinct} planes $H'$ and $H''$ containing $P$. 
Then $(\P^3, H+H'+H'')$ is $F$-split, and hence so is $X$ (Proposition \ref{p-F-split-blowup}). 

\underline{3-30}: 
By \cite[Subsection 7.3]{FanoIV}, there exist blowups 
\[ 
X \to V_7 \to \P^3,
\]
where $V_7 \to \P^3$ is a blowup at a point $P \in \P^3$ and 
the blowup centre of $X \to V_7$ is the proper transform of a line $L$ passing through $P$. 
Take two planes $H, H'$ containing $L$. 
Then $(\P^3, H+H')$ is F-split, and hence $X$ is F-split (Proposition \ref{p-F-split-blowup}). 
\end{proof}

%% file: section5.tex
\section{Fano threefolds with $\rho =2$ (except for 2-2, 2-6, 2-8)}

\subsection{Quasi-$F$-splitting (imprimitive case)}\label{ss-rho2-imprim}

\subsubsection{$D+E$ (imprimitive)}

\begin{prop}\label{p-DE}
Let $X$ be a smooth Fano threefold with $\rho(X)=2$. 
If the types of the extremal rays are $D+E$, then 
$X$ is quasi-$F$-split. 
\end{prop}

\begin{proof}
In this case, $X$ is imprimitive \cite[Subsection 7.2]{FanoIV}, i.e., 
the types of the extremal rays are $D+E_1$. 
Let $f: X \to Y$ (resp. $\pi : X \to \P^1$) 
be the contraction of type $E_1$ (resp. type $D$). 
Let $B$ be the smooth curve on $Y$ that is the blowup centre of $f: X \to Y$, i.e., $X = \Bl_B Y$.  
By \cite[Subsection 7.2]{FanoIV}, we may assume that one of the following holds. 
\begin{enumerate}
\item $Y =\P^3$ (2-4, 2-25, 2-33). 
\item $Y =Q$, where $Q$ is a smooth quadric hypersurface on $\P^4$ 
(2-7, 2-29). 
\item $Y =V_d$ with $1 \leq d \leq 5$, where 
$V_d$ is a smooth Fano threefold of index two satisfying $(-K_{V_d})^3 = 8d$
(2-1, 2-3, 2-5, 2-10, 2-14). 
%Since $Y$ is a Fano threefold of index $\geq 2$
\end{enumerate}

\begin{claim*}
A general fibre $D_X$ of $\pi : X \to \P^1$ is a canonical del Pezzo surface. 
\end{claim*}

\begin{proof}[Proof of Claim]
By \cite[Theorem 15.2]{FS20}, the generic fibre $X_K$ of $\pi$ is geometrically normal, 
where $K$ denotes the function field of $\P^1$. 
Then  its base change $X_{\overline K} := X_{K} \times_K \overline{K}$ 
to the algebraic closure $\overline K$ has at worst canonical singularities \cite[Theorem 3.3]{BT22}. 
Hence a general fibre $D_X$ of $\pi$ is normal and has at worst canonical singularities. 
This completes the proof of Claim. 
\end{proof}

(1) Assume that $Y=\P^3$. 
In this case, $B =D \cap D'$ for some surfaces $D$ of degree $1 \leq e \leq 3$, 
i.e., $Y$ is 2-4 ($e=3$),  2-25 ($e=2$), or 2-33 ($e=1$). %{\cred add the numbers}. 
Although $D$ might be singular, 
$D$ is a normal prime divisor on $X$, 
because $D$ is smooth along an effective ample Cartier divisor $B=D'|_D$. 
After replacing $D$ by a general member of the pencil generated by $D$ and $D'$, 
we may assume that $D$ is a canonical del Pezzo surface (Claim). 
If $e \in \{1, 2\}$, then $D$ is $F$-split  (Theorem \ref{t-dP-F-split}). 
Then $(Y, D)$ is $F$-split (Corollary \ref{c-F-split-IOA}), which implies that $X$ is $F$-split (Proposition \ref{p-F-split-blowup}). 
We may assume that $D$ is a cubic surface, i.e., $X$ is 2-4. 
Recall that $-K_X \sim D_X + H$, where $D_X$ denotes the proper transform of $D$ and $H\coloneqq f^*\MO_{\P^3}(1)$. 
Replacing $D_X$ and $H$ by  general members of $|D_X|$ and $|H|$,  
we obtain $D_X \cap H \simeq D \cap {f_{*}}H$. 
Since $D \cap H$ is a general hyperplane section of a normal cubic surface $D$, 
we have $D \cap {f_{*}}H \simeq D_X \cap H$ is an elliptic curve. 
Hence $X$ is quasi-$F$-split (Corollary \ref{c-IOA-Fano3}, Remark \ref{r-IOA-Fano3}).

(2) Assume that $Y=Q$. 
In this case, $B =D \cap D'$ for some surfaces $D \in |\MO_Q(e)|$ with $1 \leq e \leq 2$, i.e., $Y$ is 2-29 ($e=1$), or 2-7 $(e=2)$.
By Claim, $D$ is a canonical del Pezzo surface.
{We have 
\[
K_D^2 =(K_Q+D)^2 \cdot D=
2e(3-e)^2. 
\]}
If $e=1$, then {$K_D^2=8$}, and thus $D$ is $F$-split  (Theorem \ref{t-dP-F-split}). 
Then $(Y, D)$ is $F$-split (Corollary \ref{c-F-split-IOA}), which implies that $X$ is $F$-split (Proposition \ref{p-F-split-blowup}). 
We may assume that $e=2$, i.e., $X$ is 2-7. 
%Replacing $D_X$ and $H$ by  general members of $|D_X|$ and $|H|$, we obtain $D_X \cap H \simeq D \cap H$. 
Set $H\coloneqq f^{*}\sO_{Q}(1)$. Replacing $D$ and $f_{*}H$ by general members of $|D_X|$ and $|\sO_{Q}(1)|$,  
we obtain $D_X \cap H \simeq D \cap f_{*}H$.
Since $D \cap {f_*}H$ is a general hyperplane section of 
{a canonical del Pezzo surface $D$ of degree $4$}  with $-K_D \sim 
{f_{*}}H|_D = D \cap {f_*}H$, 
it follows that {$D\cap f_{*}H \simeq D_X \cap H$} is an elliptic curve. 
Hence $X$ is quasi-$F$-split (Corollary \ref{c-IOA-Fano3}, Remark \ref{r-IOA-Fano3}).

(3) Assume that $Y=V_d$. 
In this case,  $B$ is an elliptic curve 
which is a complete intersection of two prime divisors $D \in |-(1/2)K_Y|$ 
and $D'\in  |-(1/2)K_Y|$, 
i.e., $Y$ is 2-1, 2-3, 2-5, 2-10, or 2-14. 
We then get $K_X +D_X+D'_X +E \sim 0$, 
where $E := \Ex(f)$ and $D_X$ and $D'_X$  are the proper transforms 
of $D$ and $D'$, respectively. 
%two general fibres $D_X, D'_X$ and the exc $E$. 
Fix general members {$D$ and} $D'$ of $|-(1/2)K_Y|$ containing $B$. 
Let $D^{\gen}$ be the generic member of the pencil generated by $D$ and $D'$.
Let $\kappa$ be the function field $\kappa$ of this pencil. 
For every $k$-scheme $Z$, we set $Z_{\kappa} := Z \times_k \kappa$. 
Then $D^{\gen} \cap D'_{\kappa} = B_{\kappa}$ and 
$X_{\kappa} \to Y_{\kappa}$ is the blowup along $B_{\kappa}$. 
Note that the proper transform $D^{\gen}_X$ is regular \cite[Theorem 4.9]{Tan-Bertini}. 
Since $D^{\gen}_X \cap ((D'_X)_{\kappa}+E_{\kappa}) = D^{\gen}_X \cap E_{\kappa} \simeq B_{\kappa}$ 
is smooth over $\kappa$, 
we conclude that $X$ is quasi-$F$-split (Corollary \ref{c-IOA-Fano3}, Remark \ref{r-IOA-Fano3}).
\end{proof}

\subsubsection{$E+E$ (imprimitive)}

\begin{prop}\label{p-EE1}
Let $X$ be a smooth Fano threefold with $\rho(X)=2$. 
Assume that  the types of the extremal rays are $E_1$ and $E$. 
Then the following hold. 
\begin{enumerate}
    \item $X$ is quasi-$F$-split. 
    \item If $X$ is not 2-12, then $X$ is $F$-split. 
\end{enumerate}
\end{prop}

\begin{proof}
Let $f \colon X \to Y$ be a contraction of type $E_1$ and 
let $f'\colon X \to Y'$ be the contraction of the other extremal ray, which is of type $E$. 
Let $H_Y$ (resp.~$H_{Y'}$) be the ample Cartier divisor that generates $\Pic\,Y$ (resp.~$\Pic\,Y'$). 
Set $H \coloneqq f^*H_Y$ and $H' \coloneqq f'^*H_{Y'}$. 
The list of such smooth Fano threefolds is as follows \cite[Subseciton 7.2]{FanoIV}: 
\[
\text{2-12, 2-15, 2-17, 2-19, 2-21, 2-22, 2-23, 2-26, 2-28, 2-30}. 
\]
In particular, $Y$ is $\P^3$, $Q$, or $V_d$ with $3 \leq d \leq 5$. 
Then $|H_Y|$ is very ample, and hence we may assume that $H_Y$ is a smooth prime divisor on $Y$. 
Note that we have $-K_X \sim H+H'$ except when $X$ is 2-30 \cite[Remark 3.4 and Proposition 5.9]{FanoIII}.

\setcounter{step}{0}
\begin{step}\label{s1-EE1}
If $X$ is 2-15, 2-28, or  2-30, then $X$ is $F$-split. 
\end{step}

\begin{proof}[Proof of Step \ref{s1-EE1}]
In this case, there is a blowup $f: X \to  Y = \P^3$ along a smooth curve $B$ 
such that $B$ is contained in a prime divisor $D$ on $\P^3$ of degree $\leq 2$ 
(\cite[Proposition 9.3]{FanoIII}, \cite[Subsection 7.2]{FanoIV}). 
Since $D$ is $F$-split and $-(K_{\P^3}+D)$ is ample, it follows that $(\P^3, D)$ is $F$-split (Corollary \ref{c-F-split-IOA}), which implies that $X$ is $F$-split (Proposition \ref{p-F-split-blowup}). 
This completes the proof of Step \ref{s1-EE1}. 
\end{proof}

\begin{step}\label{s2-EE1}
If $X$ is 2-17 2-19, 2-21, 2-22, 2-23, or 2-26, then $X$ is $F$-split. 
\end{step}

\begin{proof}[Proof of Step \ref{s2-EE1}]
Replace $H$ by a general member of $|H|$. 
We now prove that $H$ is a smooth prime divisor that is $F$-split. 
We first treat the case when $X$ is 2-17, 2-19, or 2-22. 
In this case, there is a blowup $f\colon X \to Y = \P^3$ along a smooth curve $B$ of degree $\leq 5$. 
For a general plane $H_Y$, its pullback $H=f^{*}H_Y$ is nothing but the blowup of $H_Y$ along $H_Y \cap B$. 
Since $H_Y \cap B$ is smooth and $-K_H = -(K_X+H)|_H = H'|_H$ is nef and big, 
it follows that $H$ is a smooth weak del Pezzo surface. 
By $K_H^2 \geq  K_{H_Y}^2 -5 =4$, we have $H$ is $F$-split \cite[Proposition 3.6]{KT}. 
For the the remaining case (i.e., $X$ is 2-21, 2-23, or 2-26), 
we can apply the same argument, because there is a blowup $f\colon X \to Y = Q$ along a smooth curve $B$ of degree $\leq 4$, and hence we may apply \cite[Proposition 3.6]{KT}. 
Therefore, $H$ is a smooth prime divisor which is $F$-split.

In order to prove that $(X, H)$ is $F$-split, 
it is enough to show that 
\[
H^1(X,-H-(p^e-1)(K_X+H))=H^1(X, K_X +p^{e} H')=0 
\]
for some $e > 0$ by $-K_X \sim H+H'$ and Proposition \ref{p-F-split-IOA}. 
Set $E' \coloneqq \Ex(f')$. 
By the Fujita vanishing theorem \cite[Theorem 3.8.1]{Fuj17}, we can find $e_0 \in \Z_{>0}$ and $s_0 \in \Z_{>0}$   such that 
\[
H^1(X, K_X +p^e H' - s_0E' ) =0  
\]
for every integer $e \geq e_0$. 
Indeed, we can find $a_0,t_0 \in \Z_{>0}$ such that $a_0H'-t_0E$ is ample.
Then, by Fujita vanishing, there exists $m\gg0$ such that
\[
H^1(X, K_X +m(a_0 H' - t_0E')+N ) =0  
\]
for any nef Cartier divisor $N$ on $X$.
Then we take $s_0=mt_0$ and $e_0 \in \Z_{>0}$ satisfying $p^e\geq ma_0$.

Recall that $\MO_X(mH'-E')|_{E'}$ is ample for some $m \gg 0$ 
and 
$E'$ satisfies the Kodaira vanishing theorem {(for Cartier divisors), because $E'$ is a smooth projective surface with negative Kodaira dimension or a singular quadric surface \cite[Definition 3.3]{FanoIII}}.
%Again by the Fujita vanishing theorem, 
Hence, for each $0 \leq s \leq s_0$, we can find $e(s) \in \Z_{>0}$ such that 
\[
 H^1(E', \MO_X(K_X+ p^eH' -sE')|_{E'})=H^1(E', \MO_{E'}(K_{E'}+p^eH'-(s+1)E'))=0 
 \]
for every $e \geq e(s)$.
Set 
\[
e \coloneqq \max \{e_0, e(0), e(1), ..., e(s_0)\}. 
\]
It suffices to prove 
\[
H^1(X, K_X+p^eH' -s E')=0
\]
for every $0 \leq s \leq s_0$. 
By descending induction on $s$, this follows from 
\begin{multline*}
    H^1(X, K_X+p^eH' -(s+1) E') \to H^1(X, K_X+p^eH'-sE')\\ \to H^1(S, \MO_X(K_X+ p^eH' -sE'))=0. 
\end{multline*}
This completes the proof of Step \ref{s2-EE1}. 
\end{proof}

\begin{step}\label{s3-EE1}
If $X$ is 2-12, then $X$ is quasi-$F$-split. 
\end{step}

\begin{proof}[Proof of Step \ref{s3-EE1}]
In this case, 
 the contraction of each extrmeal ray is a blowup $X$ of $\P^3$  along a smooth curve $B$ of degree $6$ \cite[Subsection 7.2]{FanoIV}. 
{As in the argument in Step 2,} replacing $H$ and $H'$ by general members of $|H|$ and $|H'|$ respectively, 
we may assume that 
$H$ and $H'$ are smooth weak del Pezzo surfaces with $K_H^2=3$. 

We now finish the proof by assuming that the restriction map 
\[
\rho \colon 
H^0(X, \MO_X(-K_X-H)) \to H^0(H, \MO_H(-K_H))
\]
is surjective. 
Since a general member $C$ of $|-K_H|$ is an elliptic curve \cite[Theorem 1.4]{KN} and $H'$ is a general member in $|H'| = |-K_X-H|$, 
it follows that $H \cap H'$ is an elliptic curve.
Therefore, $X$ is quasi-$F$-split by Corollary \ref{c-IOA-Fano3}. 

It suffices to show that $\rho$ is surjective. 
By an exact sequence 
\[
0 \to \MO_X(-K_X-2H) \to \MO_X(-K_X-H) \to \MO_H(-K_H) \to 0, 
\]
it is enough to prove that $H^1(X, -K_X-2H)=0$. 
Note that 
\[
-K_X -2H \sim K_X+2H'.  
\]

Since $H'$ is a smooth rational surface and $H'|_{H'}$ is nef and big, we have $H^1(H', K_{H'}+sH')=0$ for $s \geq 0$
by \cite[ Theorem 3]{Muk13}.
Thus, we have a surjection 
\[
H^1(X, K_X+sH') \to H^1(X, K_X+(s+1)H') \to H^1(H', K_{H'}+sH')=0 
\]
for every integer $s \geq 0$.
Since we have $H^1(X, K_X)=0$, using the above surjectivity for $s\in\{0,1\}$, we obtain $H^1(X, K_X+2H')=0$. 
This completes the proof of Step \ref{s3-EE1}. 
\end{proof}
Step \ref{s1-EE1}, Step \ref{s2-EE1}, and Step \ref{s3-EE1} complete the proof of 
Proposition \ref{p-EE1}. 
\end{proof}

\subsubsection{Langer surface}

\begin{dfn}\label{d-Langer}\,
\begin{enumerate}
\item 
For all the $\F_2$-points $P_1, \ldots, P_7 \in  \P^2_{\F_2}$, 
we set 
\[
V_{L, \F_2} \coloneqq \Bl_{P_1 \amalg \cdots \amalg P_7} \P^2_{\F_2}.
\]
For a field $K$ of characteristic two, we set $V_{L, K} := V_{L, \F_2} \times_{\F_2} K$, 
which is called the {\em Langer surface} over $K$. 
\item 
For a field of characteristic two and a zero-dimensional closed subscheme $Z$ of $\P^2_K$, 
we say that $Z$ is a {\em Langer configuration} if $\Bl_Z \P^2_K$ is 
$K$-isomorphic to the Langer surface $V_{L, K}$ over $K$. 
\end{enumerate}
\end{dfn}

\begin{lem}\label{l-Langer-conf}
Let $K$ be an algebraically closed field of characteristic two. 
Take a Langer configuration $Z \subset \P^2_K$. 
Then there exists a $K$-automorphism $\sigma : \P^2_K \xrightarrow{\simeq} \P^2_K$ 
such that 
$\sigma(Z) = Z_0$, 
where 
\[
Z_0 \coloneqq \{ [1:0:0], [0:1:0], [0:0:1], [1:1:0], [1:0:1], [0:1:1], [1:1:1]\}. 
\]
\end{lem}

\begin{proof}
Fix a $K$-isomorphism $\theta \colon \Bl_Z \P^2_K \xrightarrow{\simeq} \Bl_{Z_0} \P^2_K = V_{L, K}$. 
We have two birational contractions 
\[
\varphi : \Bl_Z \P^2_K \to \P^2_K, \qquad \varphi_0 : \Bl_{Z_0} \P^2_K \to \P^2_K, 
\]
where $\varphi$ (resp. $\varphi_0$) is the blowup along $Z$ (resp. $Z_0$). 
Recall that the Langer surface $V_{L, K}$ has exactly $7$ $(-1)$-curves 
(\cite[Theorem 5.4]{CT18} or \cite[Lemma 4.5(4)]{KN}). 
Then  both $\varphi$ and $\varphi_0$ contracts all the $(-1)$-curves on 
$\Bl_Z \P^2_K$ and $\Bl_{Z_0} \P^2_K$. 
Therefore, we obtain a $K$-automorphism $\sigma : \P^2_K \to \P^2_K$ 
which completes the following commutative diagram: 
\[
\begin{tikzcd}
\Bl_Z \P^2_K \arrow[r, "{\theta, \simeq}"] \arrow[d, "\varphi"] &  \Bl_{Z_0} \P^2_K\arrow[d, "\varphi_0"]\\
\P^2_K \arrow[r, "{\sigma, \simeq}"] & \P^2_K. 
\end{tikzcd}
\]
This diagram shows that $\sigma(Z) = Z_0$. 
\end{proof}

\begin{lem}\label{l-Langer-rigid}
Let $K$ be a $C_1$-field of characteristic two %Set $K:=K(B)$ and 
and take its algebraic closure $\overline K$. 
Let $V$ be a smooth projective surface over $K$ whose base change $V_{\overline K} := V \times_{\Spec K} \Spec {\overline K}$ is $\overline K$-isomorphic to the Langer surface over $\overline K$. 
Then the following hold. 
\begin{enumerate}
\item If $K$ is perfect, then $V$ is $K$-isomorphic to the Langer surface over $K$. 
\item $\rho(V)=8$. 
\end{enumerate}
\end{lem}

\begin{proof}
We now show the implication  (1)$\Rightarrow$ (2). 
Set $K'$ to be the purely inseparable closure of $K$ in $\overline K$, i.e., 
\[
K' \coloneqq\bigcup_{e =0}^{\infty} K^{1/2^e}, \qquad K^{1/2^e} := \{ a \in \overline K\,|\, a^{2^e} \in K\}. 
\] 
Note that $K'$ is a $C_1$-field, because being  $C_1$ is stable under algebraic extensions \cite[Definition 6.2.1 and Lemma 6.2.4]{GS17}. 
Therefore, (1) is applicable to 
the perfect $C_1$-field $K'$ and the base change $V_{K'} := V \times_K K'$. 
Therefore, $V_{K'}$ is the Langer surface over $K'$, and hence $\rho(V_{K'})=8$. 
Since the field extension $K \subset K'$ is purely inseparable, 
it holds that $\rho(V) = \rho(V_{K'})=8$ \cite[Proposition 2.4]{Tan18b}. 
This completes the proof of the implication (1) $\Rightarrow$ (2). 

\medskip

It suffices to show  (1). 
Assume that $K$ is perfect. 
Recall that $V_{\overline K}$ contains the exactly $7$ $(-1)$-curves $E_1, ..., E_7$. 
Set $\overline{\Gamma} \coloneqq E_1 + \cdots +E_7$ and $\overline{\cL} := \MO_{V_{\overline K}}(\overline{\Gamma})$. 
We now show  that 
\begin{enumerate}
\item[(i)] there is an invertible sheaf $\cL$ on $V$ such that $\alpha^*\cL \simeq \overline{\cL}$, 
{where $\alpha\colon V_{\overline{K}}\to V$ is the natural morphism} and 
\item[(ii)] there exists an effective Cartier divisor $\Gamma$ on $V$ such that $\MO_X(\Gamma) \simeq \cL$ and 
the equality $\alpha^*\Gamma = \overline{\Gamma}$ of Weil divisors holds. 
\end{enumerate}
Let us show (i).
By $H^0(V, \MO_V) = K$ and $\Br(K)=0$ 
\cite[Proposition 6.2.3]{GS17}, 
we obtain $\Pic\,V \xrightarrow{\simeq} \Pic\,(V_{\overline K})^{\Gal(\overline K/K)}$ \cite[Proposition 2.3]{FanoII} (essentially due to \cite[Proposition 5.4.2]{CTS21}). 
Then it holds that 
\[
\sigma^*\overline{\cL} \simeq  \MO_{V_{\overline K}}(\sigma^*\overline{\Gamma}) \simeq \MO_{V_{\overline K}}(\overline \Gamma) \simeq \overline{\cL}
\]
for every $\sigma \in \Gal(\overline K/K)$, and hence $\overline{\cL} \simeq \alpha^* \cL$ for some $\cL \in \Pic\,V$. 
Thus (i) holds. 
Let us show (ii). 
By the flat base change theorem, we obtain 
\[
H^0(V, \cL) \otimes_K \overline{K}  \simeq H^0(V_{\overline K}, \overline L), 
\]
which implies  $\dim_K H^0(V, \cL) = \dim_{\overline K} H^0(V_{\overline K}, \overline{\cL}) 
=\dim_{\overline K} H^0(V_{\overline K}, \MO_{V_{\overline K}}(\overline{\Gamma})) =1$. 
In particular, there exists an effective Cartier divisor $\Gamma$ on $V$ such that $\MO_V(\Gamma) \simeq \cL$ and $\alpha^* \Gamma = \overline{\Gamma}$. 
Thus (ii) holds. 

\medskip

Since $\overline{\Gamma} = \sum_{i=1}^7 E_i$ is smooth, so is $\Gamma$ (note that $\Gamma$ might be irreducible, although $\overline \Gamma$ is not). 
Let $\varphi :V \to W$ be the contraction of $\Gamma$, 
where $W$ is a smooth projective surface over $K$. 
Then its base change ${\varphi}_{\overline K} : V_{\overline K} \to W_{\overline K}$ to $\overline K$ is the birational contraction of $\overline{\Gamma}$, i.e., $W_{\overline K} \simeq \P^2_{\overline K}$. 
By $\Br(K)=0$, we get 
a $K$-isomorphism $W \simeq \P^2_K$ 
\cite[Proposition 7.1.6]{CTS21}. 
Via this isomorphism, we identify $W$ and $\P^2_K$ (resp. $W_{\overline K}$ and $\P^2_{\overline K}$): 
\[
\varphi : V \to W = \P^2_K, \qquad 
{\varphi}_{\overline K} : V_{\overline K} \to W_{\overline K}=\P^2_{\overline K}. 
\]
Recall that $\varphi$ is the blowup along some closed subscheme $Z$ on $W = \P^2_K$ \cite[Ch. II, Theorem 7.17]{Har77}. 
Since blowups commutes with flat base changes \cite[Section 8, Proposition 1.12(c)]{Liu02}, 
$\varphi_{\overline K}$ is the blowup along 
the base change $Z_{\overline K}$. 
Since $Z_{\overline K}$ is a zero-dimensional reduced scheme consisting of $7$ points, 
$Z$ is a smooth zero-dimensional closed subscheme of $\P^2_K$ satisfying $h^0(Z, \MO_Z) = 7$. 

Set $H := \Hilb^7_{\P^2_{\overline{\F_2}/\overline{\F}_2}}$, 
which is the Hilbert scheme of $\P^2_{\overline \F_2} \to \Spec \overline{\F}_2$ that parametrises 
the zero-dimensional closed subschemes $W$ satisfying $h^0(W, \MO_W) =7$. 
Let $U := \Univ^7_{\P^2_{\overline{\F}_2}/\overline{\F}_2}$ be its universal family 
(cf.~\cite[Section 5]{FAG}): 
$U \subset \P^2_{\overline{\F}_2} \times_{\overline{\F}_2} H \to H$.  
Let $H_L(\overline{\F}_2) \subset H(\overline{\F}_2)$ be the subset consisting of 
the Langer configurations over $\overline{\F}_2$ 
(Definition \ref{d-Langer}(2)). 
Set $G := \PGL_{3, \overline{\F}_2}$, which is an algebraic group over $\overline{\F}_2$. 
By Lemma \ref{l-Langer-conf}, 
$H_L(\overline{\F}_2)$ is equal to the $G(\overline{\F}_2)$-orbit  of $[Z_0] \in H(\overline{\F}_2)$, 
where 
\[
Z_0 := \{ [1:0:0], [0:1:0], [0:0:1], [1:1:0], [1:0:1], [0:1:1], [1:1:1]\} \subset \P^2_{\overline{\F}_2}. 
\]
Recall that the $G(\overline{\F}_2)$-orbit $H_L(\overline{\F}_2)$ is a locally closed subset of $H(\overline{\F}_2)$ \cite[Proposition 1.65(b)]{Mil16}. 
Since $G=\PGL_3(\overline{\F}_2)$ is irreducible, so is $H_L(\overline{\F}_2)$. 
There exists an integral locally closed subscheme $H_L$ of $H$ 
whose set of the $\overline{\F}_2$-valued points coincides 
with $H_L(\overline{\F}_2)$. 
Set $U_L \coloneqq U \times_H H_L$. 
We have the Langer configuration as the fibre 
of $\pi \colon U \to H$ over $[Z_0]$: 
\[
\{P_1^U,\ldots, P_7^U\}  = \pi^{-1}([Z_0]) \subset U. 
\]
Since $G=\PGL_{3, \overline{\F}_2}$ equivariantly acts on $\pi : U \to H$, 
we have the orbits \[O_G(P_1^U), \ldots, O_G(P_7^U)\] of the above 
$7$ points $P_1^U, \ldots, P_7^U$. 
Since each $O_G(P_i^U)$ is a locally closed subset, this is a subvariety (integral scheme). 
Since any fibre of $U_L \to H_L$ is geometrically reduced, 
it follows that $U_L$ is reduced. 
Therefore, we get a scheme-theoretic equality
\[
U_L = O_G(P_1^U) \amalg \cdots \amalg O_G(P_7^U), 
\]
because we have the corresponding set-theoretic equality.

Note that $Z \subset \P^2_K$ corresponds to a $K$-rational point $[Z] \in H_K(K)$ such that 
the corresponding $\overline K$-rational point 
$[Z_{\overline K}] \in H_{\overline K}(\overline K)$ is contained 
in $(H_L \times_{\overline{\F}_2} \overline K)(\overline K)$. 
Then the image 
\[
[Z] \hookrightarrow H(K) \to H
\]
is contained in the Langer locus $H_L$ (over $\overline{\F}_2$). 
Therefore, $Z \to [Z]$ is obtained by a base change of 
\[
U_L \to H_L, 
\]
and hence $Z$ must be split up into $7$ distinct points. 
\end{proof}

\subsubsection{$C+E$ (imprimitive)}

\begin{prop}\label{p-CE1}
Let $X$ be a smooth Fano threefold such that $\rho(X)=2$ and the types of the extremal rays are $C+E_1$. 
Then the following hold. 
\begin{enumerate}
\item $X$ is quasi-$F$-split. 
\item If $X$ is not 2-9, then  $X$ is $F$-split. 
\end{enumerate}
\end{prop}

\begin{proof}
Let $\pi : X \to \P^2$ (resp. $g : X\to Y$) be the contraction 
of the extremal ray of type $C$ (resp. $E_1$). 
Let $S$ be a general member of $|\pi^*\MO_{\P^2}(1)|$ and let $H$ be a general member of $|g^*\MO_Y(1)|$, 
where $\MO_Y(1)$ denotes the ample generator of $\Pic Y$. 
Since $Y =\P^3, Q$, or $V_d$ with $3 \leq d \leq 5$ \cite[Subsection 7.2]{FanoIV}, the complete linear system $|\MO_Y(1)|$ is very ample. 
Note that $H$ is the blowup of a smooth surface $\overline{H}\in |\MO_Y(1)|$ along 
the zero-dimensional smooth closed subscheme $B \cap \overline{H}$, 
and hence $H$ is a smooth prime divisor. We have that 
\[
-K_X \sim S + \mu H 
\]
for the length $\mu$ of the extremal ray of type $C$ \cite[Proposition 5.9]{FanoIII}, i.e., 
if $\pi$ is of type $C_i$ with $i \in \{1, 2\}$, then $\mu = i$. 
Note that $-K_H=-(K_X+H)|_H=(S+(\mu-1)H)|_H$ is nef and $S^2 \cdot H =\frac{2}{i}$ \cite[Lemma 5.3 and Proposition 5.9(2)]{FanoIII}.

\setcounter{step}{0}

\begin{step}\label{s1-CE1}
If $\pi$ is of type $C_2$, then $X$ is $F$-split. 
\end{step}

\begin{proof}[Proof of Step \ref{s1-CE1}]
Assume that $\pi$ is of type $C_2$. 
In this case, 
$X$ is 2-27 or 2-31 \cite[Subsection 7.2]{FanoIV}. 
If $X$ is 2-27 (resp. 2-31), then 
$X = \Bl_C\,Y$  for $Y=\P^3$ (resp. $Y = Q$), 
where $C$ is a smooth curve of degree $3$ (resp. $1$). 
Recall that $-K_X \sim S+2H$. 
We then get $S \cdot H^2 = (-K_X) \cdot H^2 -2H^3 = (-K_Y) \cdot 
\sO_Y(1)^2 -2\sO_Y(1)^3 =2$ in both cases.
Since we have
\[
K_H^2=(S+H)^2 \cdot H=S^2\cdot H+ 2(S\cdot H^2)+H^3=1+4+H^3,
\]
it follows that $H$ is a smooth del Pezzo surface with $K_H^2 = 6$ (resp. $K_H^2=7$).
Since $-(K_X+H) \sim S+H$ is ample, so is $-K_H$. 
Thus $H$ is $F$-split  (Theorem \ref{t-dP-F-split}). 
Therefore, $X$ is $F$-split (Corollary \ref{c-F-split-IOA}). 
This completes the proof of Step \ref{s1-CE1}. 
\end{proof}

\begin{step}\label{s2-CE1}
Assume that $\pi$ is of type $C_1$. 
Then the following hold. 
\begin{enumerate}
\renewcommand{\labelenumi}{(\roman{enumi})}
\item $K_H^2 =2$ and $H$ is a smooth weak del Pezzo surface. 
\item The induced composite morphism 
\[
\pi_H \colon H \hookrightarrow X \xrightarrow{\pi} \P^2
\]
coincides with the morphism induced by the complete linear system $|-K_H|$. 
Moreover, $\pi_H$ is a generically finite morphism of degree two. 
\item 
If there exists a smooth prime divisor $C$ on $H$ satisfying $C \sim -K_H$, 
then $X$ is quasi-$F$-split. 
\end{enumerate}
\end{step}

\begin{proof}[Proof of Step \ref{s2-CE1}]
It holds that 
\[
K_H^2 = (K_X+H)^2 \cdot H = S^2 \cdot H %=2 \ell \cdot H 
=2. 
\]
Thus (i) holds. 

Let us show (ii). 
We have that $\pi^*\MO_{\P^2}(1)|_H \sim S|_H \sim  -(K_X+H)|_H \sim  -K_H$. 
The Riemann--Roch theorem, together with (i), implies $h^0(H, -K_H) = K_H^2 +1 =3$. 
By $h^0(H, -K_H) = 3 = h^0(\P^2, \MO_{\P^2}(1))$, 
the 
composition $\pi_H \colon H  \hookrightarrow X \xrightarrow{\pi} \P^2$ 
coincides with the morphism induced by the complete linear system $|-K_H|$. 
In particular, $\pi_H$ is a generically finite morphism of degree two. 
Thus (ii) holds. 

Let us show (iii). 
By (ii), we have the induced isomorphism: 
\[
H^0(H, -K_H) \xleftarrow{\pi_H^*, \simeq} H^0(\P^2, \MO_{\P^2}(1))
\]
via $\pi_H \colon H \to \P^2$, i.e., $\pi_H^* : |\MO_{\P^2}(1)| \to |-K_H|$ is bijective. 
By our assumption, there exists  a line $L \in |\MO_{\P^2}(1)|$ 
such that $C := \pi^*_H(L) \in |-K_H|$ is a smooth prime divisor. 
This property holds even after replacing $L$ by a general member of $|\MO_{\P^2}(1)|$, 
and hence we may assume that $S = \pi^*L$. 
We then obtain $S \cap H =\pi^*L \cap H = \pi^*_H(L)=C$, which is a smooth elliptic curve. 
By $K_X+S+H \sim 0$, $X$ is quasi-$F$-split (Corollary \ref{c-IOA-Fano3}). 
Thus (iii) holds. 
This completes the proof of Step \ref{s2-CE1}. 
\end{proof}

\begin{step}\label{s3-CE1}
Assume that $\pi$ is of type $C_1$. 
Then there exists a smooth prime divisor $C$ on $H$ satisfying $C \sim -K_H$. 
\end{step}

\begin{proof}[Proof of Step \ref{s3-CE1}]
Recall that $H$ is a general member of $|g^*\MO_Y(1)|$. 
Suppose that any member of $|-K_H|$ is singular. 
By Step \ref{s2-CE1}(iii), it is enough to derive a contradiction. 
Note that, for every line $L$ on $\P^2$ and every smooth member $H \in |g^*\MO_Y(1)|$, 
its intersection $H \cap \pi^*L$ is not smooth (Step \ref{s2-CE1}(ii)). 

We now show that every smooth member of $|H|$ is isomorphic to the Langer surface over $k$ (Definition \ref{d-Langer}). 
The Stein factorisation $H'$ of the composition $\pi_H : H \hookrightarrow X \xrightarrow{\pi} \P^2$ is the anti-canonical model of $H$ (Step \ref{s2-CE1}(ii)). 
Note that each fibre of $\pi_H : H \hookrightarrow X \xrightarrow{\pi} \P^2$ 
is contained in a fibre of $\pi : X \to \P^2$, which is a conic. 
In particular, any singularity on $H'$ is either $A_1$ or $A_2$. 
By $K_H^2=2$ (Step \ref{s2-CE1}(i)) and  \cite[Theorem 1.4]{KN}, 
it holds that $p=2$ and $H$ is isomorphic to the Langer surface.

Fix two general members $H_1$ and $H_2$ of $|H|$. 
In particular, $\Gamma \coloneqq H_1 \cap H_2$ is a smooth curve and 
each of $H_1$ and $H_2$ is isomorphic to the Langer surface. 
Let $\sigma: Y \to X$ be the blowup along $\Gamma  = H_1 \cap H_2$, 
and hence we get the morphism $\alpha \colon Y \to \P^1$ induced by the pencil generated by 
$H_1^Y$ and $H_2^Y$, where $H^Y_i \coloneqq \sigma^*H_i - \Ex(\sigma)$, which coincides with the proper transform of $H_i$ on $Y$.
By construction, every general fibre of $\alpha$ is isomorphic to the Langer surface 
(as otherwise we could find a smooth member of $|H|$ which is not the Langer surface). 
Set $V \coloneqq Y \times_{\P^1} \Spec K$ to be the generic fibre of $\alpha : Y \to \P^1$, 
where $K \coloneqq \Frac\,\P^1$.

For the algebraic closure $\overline{K}$ of $K$, 
it is enough to show, by Lemma \ref{l-Langer-rigid}, that
the base change $V_{\overline K} \coloneqq V \times_{\Spec K} \Spec \overline K$ is isomorphic to the Langer surface over $\overline K$. 
In fact, this implies $8=\rho(V)=\rho(Y \times_{\P^1} \Spec K)=\rho(Y)-\rho(\P^1)=1$, 
which is a contradiction.
Fix a general fibre $V'$ of $\alpha : Y \to \P^1$. 
Let $F_1, \ldots, F_n$ be all the $(-2)$-curves on $V_{\overline K}$ 
($F$ is called a {\em $(-2)$-curve} if $F^2 =-2$ and $F \simeq \P^1$). 
Since $F_1, \ldots, F_n$ can be defined around the generic point of the base and $V'$ is general,
we obtain the corresponding $(-2)$-curves $F'_1,\ldots, F'_n$ on $V'$. 
By the invariance of intersection numbers for flat families, 
we see that $F'_i \cdot F'_j = F_i \cdot F_j$ for every $i, j$. 
Since $V'$ is the Langer surface over $k$, 
there are exactly $7$ $(-2)$-curves and they are mutually disjoint \cite[Theorem 5.4]{CT18}. 
Then we see that $n \leq 7$ and $F'_i \cdot F'_j = 0$ (i.e., $F'_i \cap F'_j = \emptyset$) for every $1 \leq i < j \leq n$. 
As $|-K_{V'}|$ has no smooth member, neither does $|-K_{V_{\overline K}}|$.  
By \cite[Theorem 1.4]{KN}, we get $n=7$, i.e., $V_{\overline K}$ is 
isomorphic to the Langer surface over $\overline K$. 
This completes the proof of Step \ref{s3-CE1}. 
\end{proof}
Step \ref{s1-CE1}, 
Step \ref{s2-CE1},  and 
Step \ref{s3-CE1} complete the proof of Proposition \ref{p-CE1}. 
\end{proof}

\subsection{Quasi-$F$-splitting (primitive case)}\label{ss-rho2-prim}

In Section \ref{ss-rho2-imprim}, 
the quasi-$F$-splitting for smooth Fano threefolds with $\rho=2$ 
has been settled for the imprimitive case. 
Hence the remaining cases are as follows \cite[Subsection 7.2]{FanoIV}: 
\[
\text{2-2, 2-6, 2-8, 2-18, 2-24, 2-32, 2-34, 2-35, 2-36}. 
\]
In what follows, we shall settle the cases except for 2-2, 2-6, and 2-8. 
These cases will be treated in Section \ref{sec:8}. 

\begin{lem}\label{l-prim-toric}
If $X$ is a smooth Fano threefold of No.~2-32, 2-34, 2-35, or 2-36, 
then $X$ is  $F$-split. 
\end{lem}

\begin{proof}
It is well known that $X$ is toric except when it is 2-32. 
Assume that $X$ is 2-32, i.e., $X$ is a smooth hypersurface on $\P^2 \times \P^2$ of bidegree $(1, 1)$. 
Note that $f_1 \colon X \hookrightarrow \P^2 \times \P^2 \xrightarrow{\pr_1} \P^2$ 
is a $\P^1$-bundle \cite[Subsection 7.2]{FanoIV}.  
Take a general member $D \in |f_1^*\MO_{\P^2}(1)|$, 
which is a $\P^1$-bundle over $\P^1$. 
Hence $D$ is $F$-split. 
As in the proof of Proposition \ref{p-rho3}(2), we can see
$-(K_X+D)$ is ample by \cite[Proposition 5.9(3)]{FanoIII} and \cite[Section 7.2]{FanoIV}.
Therefore, $X$ is $F$-split (Corollary \ref{c-F-split-IOA}). 
\end{proof}

\begin{lem}\label{l-2-18}
Let $Y$ be a smooth Fano threefold of No.~2-18. 
Then the following hold. 
\begin{enumerate}
\item $Y$ is $F$-split. 
\item 
Let $B$ be a smooth fibre of the contraction $g \colon Y \to \P^2$ of type $C_1$. Then the blowup $X \coloneqq \Bl_B Y$ of $Y$ along $B$ 
is a smooth Fano threefold of No.~3-4. 
\end{enumerate}
\end{lem}

The following argument is almost identical to that of \cite[Proposition 4.35]{FanoIV}.

\begin{proof}
By Proposition \ref{p-rho3}, (2) implies (1). 
Let us show (2). 
It is enough to show that $-K_X$ is ample \cite[Subsection 7.2]{FanoIV}. 
By construction, we have the following commutative diagram except for $f_1, g_{11}, g_{12}$. 
Since $\varphi_1 \times \varphi_2 \colon X \to Z_1 \times Z_2 = \P^1 \times \P^1$ is not a finite morphism,  
its Stein factorisation $f_1 \colon X \to S$ of $\varphi_1 \times \varphi_2$ is not an isomorphism. 
We get $\dim S=2$, because we have $\dim S \leq \dim (\P^1 \times \P^1)=2$, 
and $S \to \P^1 \times \P^1 \xrightarrow{\pr_1} \P^1$ is not an isomorphism. 
Then $f_1 \colon X \to S$ is a contraction of a (possibly non $K_X$-negative) extremal ray. 
Since this extramal ray is contained in the two-dimensional extremal faces corresponding to $\varphi_1$ and $\varphi_2$, 
$\NE(X)$ is generated by three extremal rays, which are corresponding to 
$f_1, f_2, f_3$. 
\[
\begin{tikzpicture}[commutative diagrams/every diagram,
    declare function={R=3;Rs=R*cos(60);}]
 \path 
  (0,0)  node(X) {$X$} 
  (90:R) node (Y1) {$S$}
  (210:R) node (Y2) {$\F_1$}
  (-30:R) node (Y3) {$Y$}  
  (30:Rs) node(Z1) {$\P^1$} 
  (150:Rs) node(Z2) {$\P^1$} 
  (270:Rs) node(Z3) {$\P^2$};
 \path[commutative diagrams/.cd, every arrow, every label]
 (X) edge[swap] node {$f_1$} (Y1)
 (X) edge[swap] node {$f_2$} (Y2)
 (X) edge node {$f_3$} (Y3)
 (X) edge node {$\varphi_1$} (Z1)
 (X) edge[swap] node {$\varphi_2$} (Z2)
 (X) edge[swap] node {$\varphi_3$} (Z3)
 (Y1) edge node {$g_{11}$} (Z1)
 (Y1) edge[swap] node {$g_{12}$} (Z2)
 (Y2) edge node {$g_{22}$} (Z2)
 (Y2) edge[swap] node {$g_{23} =\tau$} (Z3)
 (Y3) edge node {$g_{33}$} (Z3)
 (Y3) edge[swap] node {$g_{31}$} (Z1);
\end{tikzpicture}
\]
By the same argument as in \cite[Proposition 4.35]{FanoIV}, we obtain 
\[
-K_X \sim H_1 + H_2 +H_3, 
\]
where each $H_i$ is the pullback of the ample generator by $\varphi_i$. 
Since $\varphi_1 \times \varphi_2 \times \varphi_3 \colon X \to \P^1 \times \P^1 \times \P^2$ is a finite morphism, 
$-K_X \sim H_1 + H_2 +H_3$ is ample, as required. 
\end{proof}

\begin{lem}\label{l-2-24-gene-sm}
Let $Y$ be a smooth Fano threefold of No.~2-24. 
Assume that the 
contraction $\psi \colon Y \to \P^2$  of type $C_1$ is generically smooth. 
Then the following hold. 
\begin{enumerate}
\item $Y$ is $F$-split. 
\item
Let $B$ be a smooth fibre of $\psi \colon Y \to \P^2$. 
Then the blowup $X := \Bl_B Y$ of $Y$ along $B$ is a Fano threefold of No.~3-4. 
\end{enumerate}
\end{lem}

The following argument is almost identical to that of \cite[Proposition 4.40]{FanoIV}. 

\begin{proof}
Since (2) implies (1), it is enough to show (2). 
Let us show (2). It is enough to show that $-K_X$ is ample \cite[Subsection 7.2]{FanoIV}. 
By construction, we get the following commutative diagram 
except for $f_3, g_{31}, g_{33}$. 
By the same argument as in \cite[Proposition 4.40]{FanoIV}, we see that 
\begin{itemize}
\item     $\varphi_3 \times \varphi_1: X \to \P^2 \times \P^1$ is not a finite morphism, and 
\item $-K_X \sim H_1 + H_2 +H_3$, 
\end{itemize}
where each $H_i$ is the pullback of the ample generator by $\varphi_i$. 
Then $X$ has exactly three extremal rays and $-K_X$ is ample.  
\[
\begin{tikzpicture}[commutative diagrams/every diagram,
    declare function={R=3;Rs=R*cos(60);}]
 \path 
  (0,0)  node(X) {$X$} 
  (90:R) node (Y1) {$\F_1$}
  (210:R) node (Y2) {$Y$}
  (-30:R) node (Y3) {$\P^2 \times \P^1$}  
  (30:Rs) node(Z1) {$\P^1$} 
  (150:Rs) node(Z2) {$\P^2$} 
  (270:Rs) node(Z3) {$\P^2$};
 \path[commutative diagrams/.cd, every arrow, every label]
 (X) edge[swap] node {$f_1$} (Y1)
 (X) edge[swap] node {$f_2$} (Y2)
 (X) edge node {$f_3$} (Y3)
 (X) edge node {$\varphi_1$} (Z1)
 (X) edge[swap] node {$\varphi_2$} (Z2)
 (X) edge[swap] node {$\varphi_3$} (Z3)
 (Y1) edge node {$g_{11}$} (Z1)
 (Y1) edge[swap] node {$g_{12} =\tau$} (Z2)
 (Y2) edge node {$g_{22}$} (Z2)
 (Y2) edge[swap] node {$g_{23}$} (Z3)
 (Y3) edge node {$g_{33}$} (Z3)
 (Y3) edge[swap] node {$g_{31}$} (Z1);
\end{tikzpicture}
\]    
\end{proof}

\begin{lem}\label{l-2-24-wild}
Let $Y$ be a smooth Fano threefold of No.~2-24. 
Assume that the contraction $\psi \colon Y \to \P^2$  of type $C_1$ is not generically smooth. 
Then the following hold. 
\begin{enumerate}
\item $X$ is isomorphic to $\{ x_0y_0^2 + x_1y_1^2+x_2 y_2^2 =0\} \subset \P^2_x \times \P^2_y$, 
where $\P^2_x \coloneqq \Proj\,k[x_0, x_1, x_2] (\simeq \P^2)$ 
and $\P^2_y \coloneqq \Proj\,k[y_0, y_1, y_2] (\simeq \P^2)$. 
\item $X$ is quasi-$F$-split. 
\end{enumerate}
\end{lem}

\begin{proof}
Since (1) implies (2) \cite[Example 7.13]{KTY}, 
it is enough to show (1). 
Recall that $\psi$ can be written as follows: 
\[
\psi: X \hookrightarrow \P^2_x \times \P^2_y \xrightarrow{\pr_1} \P^2_x. 
\]
Since every fibre of $\psi$ is a non-reduced conic, we can write 
\[
X = \{ f_0(x) y_0^2 + f_1(x)y_1^2 + f_2(x) y_2^2 =0\}, 
\]
where each $f_i(x) \in k[x_0, x_1, x_2]$ is a homogeneous polynomial of degree $1$. 
We see that 
\begin{enumerate}
\item[$(*)$] none of $f_0(x), f_1(x), f_2(x)$ is zero.  
\end{enumerate}
Indeed, if $f_2(x) =0$, then the affine open subset of $X$ defined by $\{ x_0 \neq 0\} \cap \{ y_2 \neq 0\}$ contains a singular point. 
After applying a suitable coordinate change, we may assume that $f_0(x)=x_0$. 
Therefore, we can write 
\[
X = \{ x_0 y_0^2 + f_1(x)y_1^2 + f_2(x) y_2^2 =0\}.  
\]
%If $f_1(x) = f_2(x)=0$, then $X$ would be singular. Hence we may assume that 
We can write $f_1(x) =a x_0 +bx_1+cx_2$. 
By applying $y_0 \mapsto y_0 + \sqrt{a}y_1$, we may assume that $a=0$. 
Note that $(*)$ implies $f_1(x) \neq 0$. 
Replacing  $f_1(x) = bx_1 +cx_2 (\neq 0)$ by $x_1$, we may assume that $f_1(x)=x_1$, i.e., 
\[
X = \{ x_0 y_0^2 + x_1y_1^2 + f_2(x) y_2^2 =0\}.  
\]
Similarly, by applying a coordinate change $y_0 \mapsto y_0 +dy_2, y_1 \mapsto y_1 +e d_3$ for some $d, e \in k$, 
the problem is reduced to the case when $f_2(x) = \alpha x_2$ for some $\alpha \in k$. 
By $(*)$, we may assume that  $\alpha=1$. Thus (1) holds. 
\end{proof}

\begin{rem}
    Let $X$ be a smooth Fano threefold of No.~2-24.
    Combining Lemma \ref{l-2-24-gene-sm}, Lemma \ref{l-2-24-wild}, and \cite[Example 7.13]{KTY},  
$X$ is $2$-quasi-$F$-split. 
    Moreover, 
    the following hold. 
    \begin{enumerate}
    \item The following are equivalent. 
    \begin{enumerate}
        \item The quasi-$F$-split height is $1$, i.e., $X$ is $F$-split. 
        \item The contraction $\psi : X \to \P^2$ of type $C_1$ is generically smooth. 
    \end{enumerate}
    \item The following are equivalent. 
    \begin{enumerate}
        \item[(a)] The quasi-$F$-split height is $2$. 
        \item[(b)] The contraction $\psi : X \to \P^2$ of type $C_1$ is not generically smooth (called wild). 
        \item[(c)] $X \simeq \{ x_0 y_0^2 + x_1y_1^2 + x_2 y_2^2 =0\} \subset \P^2_x \times \P^2_y$.
    \end{enumerate}
    \end{enumerate}
\end{rem}

\begin{prop}\label{p-rho2}
Let $X$ be a smooth Fano threefold with $\rho(X)=2$. 
Then $X$ is quasi-$F$-split if $X$ is none of 2-2, 2-6, and 2-8. 
\end{prop}

\begin{proof}
If $X$ is imprimitive (resp.~primitive), 
then the assertion follows from Proposition \ref{p-DE}, Proposition \ref{p-EE1}, and 
Proposition \ref{p-CE1} (resp.~Section \ref{ss-rho2-prim}). 
\end{proof}

\subsection{$F$-splitting}

\begin{thm}\label{t-GFS}
Assume $p>5$. 
Let $X$ be a smooth Fano threefold with $\rho(X) \geq 2$.  
If $X$ is neither 2-2 nor 2-6, then $X$ is $F$-split. 
\end{thm}

\begin{proof}
By 
Proposition \ref{p-rho6}, 
Proposition \ref{p-rho5}, 
Proposition \ref{p-rho4}, and 
Proposition \ref{p-rho3}, 
we may assume that $\rho(X)=2$. 
The types $R_1+R_2$ of the extremal rays $R_1$ and $R_2$ are as follows, because the case $D+D$ does not occur \cite[Subsection 7.2]{FanoIV}: 
\begin{enumerate}
\item[(I)] $C+E$ or $D+E$. 
\item[(II)] $E+E$. 
\item[(III)] $C+D$. 
\item[(IV)] $C+C$. 
\end{enumerate}
For each $i \in \{1, 2\}$, let $f_i : X \to Y_i$ be the contraction of $R_i$, $\mu_i$ denotes the length of $R_i$, and 
set $H_i$ to be the pullback of the ample generator of $\Pic\,Y_i (\simeq \Z)$. 
Recall that we can write $-K_X \sim \mu_2 H_1 + \mu_1 H_2$ \cite[Proposition 5.9]{FanoIII}. 
In what follows, we treat the above four cases separately.

\medskip

(I) 
Since $R_1+R_2$ is $C+E$ or $D+E$, 
we have $Y_1 =\P^1$ or $Y_1 = \P^2$. 
Pick a general member $S$ of $|H_1|$. 

\begin{claim*}
$S$ is a canonical weak del Pezzo surface. 
\end{claim*}

\begin{proof}[Proof of Claim]
We first treat the case when $Y_1 = \P^1$. 
Let $S'$ be the geometric generic fibre of $f_1 \colon X \to Y_1 = \P^1$. 
Then $S'$ is normal and $-K_{S'}$ is ample \cite[Theorem 15.2]{FS20}, 
which implies that $S'$ is canonical \cite[Theorem 3.3]{BT22}. 
Hence $S$ is a canonical del Pezzo surface for the case when $Y_1 =\P^1$. 

It is enough to settle the case when  $Y = \P^2$. 
In this case, $S$ is the inverse image of a general line $L$ on $\P^2$. 
Since the discriminant scheme $\Delta_{f_1}$ of $f_1 \colon X \to S= \P^2$ is a reduced divisor \cite[Proposition 7.2]{Tan-conic}, 
the scheme theoretic intersection $L \cap \Delta_{f_1}$ is a zero-dimensional smooth scheme. 
For the resulting conic bundle $g\colon S = X \times_{Y_1} L \to L$, 
we have $\Delta_g = L \cap \Delta_{f_1}$ \cite[Remark 3.4]{Tan-conic}. 
Since $L$ and $\Delta_g$ are smooth, also $S$ is smooth \cite[Theorem 4.4]{Tan-conic}. 
We have that $-K_X \sim S + H$ for $H\coloneqq \mu_1 H_2 + (\mu_2-1)H_1$. 
In particular, $H$ and $H|_S$ are nef and big. 
By the adjunction formula: $K_S \sim (K_X+S)|_S \sim -H|_S$, we have
$S$ is a smooth weak del Pezzo surface. 
This completes the proof of Claim. 
\end{proof}

By Claim, 
$S$ is $F$-split (Theorem \ref{t-dP-F-split}). 
We have  $-K_X \sim S+H$ for $H\coloneqq \mu_1 H_2 + (\mu_2-1)H_1$. 
By $H^1(X, -S+(p^e-1)(K_X+S))=  H^1(X, K_X+p^e H)$ and  Proposition \ref{p-F-split-IOA}, it is enough to find $e \in \Z_{>0}$ such that 
\[ 
H^1(X, K_X+p^e H)=0\] by. 
Fix $e_1 \in \Z_{>0}$ such that $p^{e_1}H -D$ is ample for $D := \Ex(f_2)$. 
By the Serre vanishing theorem, there is $e_2 \in \Z_{>0}$ such that 
$H^1(X, K_X+ p^{e_2}(p^{e_1}H -D))=0$. 
Set $e:=e_1 + e_2$. 
It suffices to prove 
\[
H^1(X, K_X + p^e H -sD)=0
\]
for every $0 \leq s \leq p^{e_2}$ by descending induction on $s$. 
The base case $s = p^{e_2}$ has been checked already. 
Fix an integer $s$ satisfying $0 \leq s < p^{e_2}$. 
By the induction hypothesis, %$H^1(X, K_X +p^eH -D) =0$ {\cred ref}, 
we have  the following exact sequence 
\begin{multline*}
    0 = H^1(X, K_X +p^eH -(s+1)D) \to H^1(X, K_X+p^eH -sD )\\
    \to H^1(D, K_D +(p^eH-(s+1)D)|_D). 
\end{multline*}
It suffices to show $H^1(D, K_D+(p^eH-(s+1)D)|_D)=0$. 
Note that $p^eH -(s+1)D$ is ample, because so are $H$ and $p^eH - p^{e_2}D (=p^{e_2}(p^{e_1}H-D)$). 
Then we get $H^1(D, K_D+(p^eH-(s+1)D)|_D)=0$ 
by the fact that $D$ is toric or a smooth ruled surface \cite[Theorem 3]{Muk13}. 
\medskip

(II) 
Assume that $R_1+R_2$ is $E+E$. 
The list of such Fano threefolds is as follows \cite[Subsection 7.2]{FanoIV}: 2-12, 2-15 2-17, 2-19, 2-21, 2-22, 2-23, 2-26, 2-28, 2-30. 
We treat the following two cases separately: 
\begin{enumerate}
\item  2-12, 2-15, 2-17, 2-19, 2-22,  2-28, 2-30. 
\item  2-21, 2-23, 2-26. 
\end{enumerate}
If (1) (resp. (2)) holds, 
then there is a blowup $f_1: X \to Y_1=\P^3$ (resp. $f_1:X \to Y_1=Q$) 
along a smooth curve $B$ on $Y_1$. 
Note that $H_1$ is the inverse image of the corresponding member $\overline{H}_1$ on $Y_1$. 
After replacing $H_1$ by a general member of $|H_1|$, 
we may assume that $H_1$ is the blowup along 
a smooth zero-dimensional scheme $B \cap \overline{H}_1$ of a smooth surface $\overline{H}_1$, and hence $H_1$ is a smooth projective surface. 
We have that $-K_X \sim H_1 + H'_2$ for some $H'_2 \sim H_2 + N$ with $N$ nef. 
Then it holds that $-K_{H_1}$ is nef and big. 
Hence $H_1$ is a smooth weak del Pezzo surface. 
The same argument as in (I) deduces that $X$ is $F$-split.

\medskip

(III) Assume that the types of the extremal rays are $C+D$. 
Since we are assuming that $X$ is not 2-2, 
$X$ is 2-18 or 2-34 \cite[Subsection 7.2]{FanoIV}. 
If $X$ is 2-34, i.e., $X \simeq \P^1 \times \P^2$), 
then $X$ is clearly $F$-split. 
The case when $X$ is 2-18 has been settled in Lemma \ref{l-2-18}.
\medskip

(IV) Assume that $R_1+R_2$ is $C+C$. 
Since we are assuming that $X$ is not 2-6, 
$X$ is 2-24 or 2-32. 
If $X$ is 2-24 (resp. 2-32), 
then $X$ is $F$-split by Lemma \ref{l-2-24-gene-sm} (resp. Lemma \ref{l-prim-toric}). 
Note that an arbitrary conic bundle $X \to S$ is generically smooth by $p>2$. 
\qedhere 

\end{proof}

%% file: section6.tex
\section{$F$-splitting and Quasi-$F$-splitting via Cartier operators}\label{sec:8}

\subsection{Quasi-$F$-splitting for 2-2, 2-6, 2-8, and 3-10}

\subsubsection{Preparation}

\begin{lem}\label{lem:SRC1}
Let $X$ be a smooth Fano threefold. 
Assume that $X$ is SRC. 
Then $H^0(X, \Omega_X^i(-D))=0$ for every $i>0$ and every pseudo-effective Cartier divisor $D$ on $X$.
\end{lem}

For the definition of SRC (separable rational connectedness), 
we refer to  \cite{Kol96}. 

\begin{proof}
The assertion follows from the essentially same proof as in \cite[Proposition 3.4]{Kaw1} 
by using the fact that the restriction of a pseudo-effective Cartier divisor $D$ to a general curve is 
pseudo-effective. 
\end{proof}

We repeatedly use the following basic lemma.

\begin{lem}\label{lem:SRC}
    Let $X$ be a smooth Fano threefold. 
Assume that there exists a conic bundle $f: X \to S$, i.e., 
$f$ is a flat morphism to a smooth projective surface $S$ such that every fibre $f^{-1}(s)$ is isomorphic to 
a conic (cf. \cite[Definition 2.3]{Tan-conic}). 
    Then the following hold. 
    \begin{enumerate}
    \item $X$ is SRC. 
    \item $H^0(X, \Omega_X^i(-D))=0$ 
    for every $i>0$ and every pseudo-effective Cartier divisor $D$ on $X$.
    \end{enumerate}
    In particular, if $X$ is one of No.~2-2, 2-6, 2-8, and 3-10, then $X$ satisfies the condition (1) in Proposition \ref{prop:criterion for qFs}.
\end{lem}

\begin{proof}
Since (1) implies (2) (Lemma \ref{lem:SRC1}), 
it is enough to show  (1). 
If $f$ is generically smooth, 
then $S$ is a smooth rational surface \cite[Proposition 3.13]{FanoIV}, and hence $X$ is SRC by \cite[Theorem 0.5]{GLP15}. 
We may assume that $f$ is wild, i.e., not generically smooth. 
Then $X$ is 2-24 or 3-10 by \cite[Corollary 8]{MS03}. 
Assume that $X$ is 2-24. 
Then the contraction of the other extremal ray is of type $C$ and its contraction gives a generically smooth conic bundle structure. We are done by the generically smooth case. 
Assume that $X$ is 3-10. 
Then $X$ is obtained as a blowup of $Q$ 
   (cf.~\cite[Corollary 8]{MS03} or \cite[Section 7]{FanoIV}). 
   Since $Q$ is rational (and hence SRC), so is $X$. 
   Thus (1) holds. 
\end{proof}

\begin{nothing}[Double covers]\label{n-double-cover}
Given smooth projective varieties $X$ and $Y$, 
we say that $f \colon X \to Y$ is a {\em double cover} if 
$f$ is a finite surjective morphism 
such that the induced field extension $K(X)/K(Y)$ is of degree $2$. 
Recall that 
we have an exact sequence 
\[
0 \to \MO_Y \to f_*\MO_X \to \MO_Y(-L) \to 0 
\]
for some Cartier divisor $L$ on $Y$ \cite[Lemma A.1]{Kaw2}. 
In particular, $\MO_Y(L) \simeq (f_*\MO_X/\MO_Y)^{-1}$. Moreover, $K_X \sim f^*(K_Y+L)$ \cite[Proposition 0.1.3]{CD89}. 
We say a double cover $f\colon X \to Y$ is {\em split} if 
the induced homomorphism 
$\MO_Y \to f_*\MO_X$ splits as an $\MO_Y$-module homomorphism (i.e., the above exact sequence splits). 
In this case, we obtain $f_*\MO_X \simeq \MO_Y \oplus \MO_Y(-L)$. 
\end{nothing}

\begin{lem}\label{l-dc-to-bdl}
Let $f: X \to Y$ be a split double cover of smooth projective varieties. 
Let $L$ be a Cartier divisor $L$ on $Y$ satisfying $\MO_Y(-L) \simeq f_*\MO_X/\MO_Y$ (cf. (\ref{n-double-cover})). 
Then there exists a closed immersion 
 $j \colon X \hookrightarrow P:=\P_Y(\MO_Y \oplus \MO_Y(-L))$ which satisfies the following properties. 
\begin{enumerate}
\item $K_P+X \sim g^*(K_Y+L)$, where $g: P=\P_Y(\MO_Y \oplus \MO_Y(-L)) \to Y$ denotes the induced $\P^1$-bundle. 
\item For the section $S := \P_Y( \MO_Y(-L))$ of 
$g\colon P=\P_Y(\MO_Y \oplus \MO_Y(-L)) \to Y$ corresponding to the second projection $\MO_Y \oplus \MO_Y(-L) \to \MO_Y(-L)$, 
it holds that $S \cap X = \emptyset$, $S|_S \sim -g^*L|_S$, $X -2S \sim 2g^*L$, and 
$\MO_P(1) \simeq \MO_P(S)$, 
\item 
 For the section $T \coloneqq \P_Y( \MO_Y)$ of  $g: P=\P_Y(\MO_Y \oplus \MO_Y(-L)) \to Y$ corresponding to the first projection $\MO_Y \oplus \MO_Y(-L) \to \MO_Y$, 
it holds that 
$S \cap T =\emptyset$, 
$T|_T \sim g^*L|_T$, $X \sim 2T$, $T-S \sim g^*L$, and $\MO_P(1) \simeq \MO_P(T-g^*L)$. 
\item 
$\Omega_{P/Y}^1 \simeq \MO_P(-g^*L -2S)$. 
\end{enumerate}
    
\end{lem}

\begin{proof}
% Recall that 
% \begin{enumerate}
%     \item $S \cap X = \emptyset$. 
% \end{enumerate}
% $X$ is a closed subscheme of the $\bA^1$-bundle $\Spec_Y\mathcal A$, 
% where $\mathcal A = \bigoplus_{d \geq 0} \MO_Y(-dL)$ {\cred ref}. 
% Since $\P_Y(\MO_Y \oplus \MO_Y(-L))$  
% Fix a fibre $\zeta$ of the $\P^1$-bundle $g:  P=\P_Y(\MO_Y \otimes \MO_Y(-L)) \to Y$. 
In what follows, we only treat the case when $p=2$, as otherwise the problem is easier. 
By \cite[the proof of Proposition 0.1.3]{CD89} (note that $f$ splits if and only if 
$f$ corresponds to a splittable admissible triple), 
there is a closed immersion $j^{\circ} : X \hookrightarrow P^{\circ}$ 
to the $\A^1$-bundle 
\[
P^{\circ} := \Spec_Y(\bigoplus_{d=0}^{\infty} \MO_Y(-dL)). 
\]
Since $P^{\circ}$ is an open subscheme of $P$, we obtain a closed immersion $j : X \to P$ over $Y$. 
By definition, we get $S \cap T =\emptyset$, 
$\MO_P(1)|_S \simeq \MO_S(-L)$, and $\MO_P(1)|_T \simeq \MO_T$, 
where $(g|_S)^*L$ denotes $L$ for $g|_S : S \xrightarrow{\simeq} Y$ by abuse of notation. 
We can write $\MO_P(1) \simeq \MO_P(S+g^*D_S) \simeq \MO_P(T+g^*D_T)$ 
for some Cartier divisors  $D_S$ and $D_T$ on $Y$. 
By $S \cap T = \emptyset$ and $\MO_P(1)|_T \simeq \MO_T$, 
we obtain $\MO_T \simeq \MO_P(1)|_T \simeq \MO_P(S+g^*D_S)|_T \simeq 
(g^* \MO_Y(D_S))|_T$, which implies $D_S \sim 0$. 
Similarly, we obtain $D_T \sim -L$ by 
$\MO_S(-L) \simeq \MO_P(1)|_S  \simeq \MO_P(T+g^*D_T)|_S \simeq (g^*\MO_Y(D_T))|_S$. Hence we get 
\[
\MO_P(1) \sim S \sim T -g^*L, 
\]
which implies 
\[
S|_S \simeq -g^*L|_S, \qquad T|_T \sim g^*L|_T. 
\]

% We have $\MO_P(1) \sim S + g^*D$ for some Cartier divisor $D$. 
% By $\MO_P(1)|_S \sim -L$ and $S|_S \sim -L$, we get $E \sim 0$, i.e., $\MO_P(1) \sim S$. 

% We can write $T-S \sim g^*D$ for some Cartier divisor $D$ on $Y$. 
% By $L \sim -S(T-S \sim g^*D$, we obtain $D \sim L$, 
% i.e., $T-S \sim g^*L$. 
% In particular, $T|_T \sim L$ and $S|_T \sim -L$. 

\begin{claim*}
The following hold. 
\begin{enumerate}
\renewcommand{\labelenumi}{(\alph{enumi})}
\item $\MO_P(K_P) \otimes \MO_P(2) \otimes g^*\MO_Y(2L) \simeq 
g^*\MO_Y(K_Y+L)$.  
\item $X \cap S = \emptyset$. 
\item $K_P+X \sim g^*(K_Y+L)$. 
\end{enumerate}
\end{claim*}

\begin{proof}[Proof of Claim] 
Let us show (c)' below, which is weaker than (c): 
\begin{enumerate}
\item[(c)'] $K_P+X \equiv g^*(K_Y+L)$, where $\equiv$ denotes the numerical equivalence. 
%\item $\MO_P(K_P) \otimes \MO_P(2) \otimes g^*\MO_Y(2L) \simeq g^*\MO_Y(K_Y+L)$.  
%\item $X \cap S = \emptyset$. 
\end{enumerate}
Since $f : X \to Y$ is a double cover, 
we can find a Cartier divisor $E$ on $Y$ such that $K_P+X \sim g^*E$. 
Then 
\[
f^*E = (g^*E)|_X \sim (K_P+X)|_X \sim K_X \sim f^*(K_Y+L), 
\]
which implies $E \equiv K_Y+L$ \cite[Corollary 1(ii) in page 304]{Kle66}. 
Thus (c)' holds. 

The assertion (a) holds by the following (cf. \cite[Proposition 7.1(2)]{FanoIII}):
\[
\MO_P(K_P) \simeq \MO_P(-2) \otimes g^*(\omega_Y \otimes \det (\MO_Y \oplus \MO_Y(-L))) \simeq \MO_P(-2) \otimes g^*\MO_Y(K_Y-L). 
\]
Let us show (b). 
By (a) and (c)', we obtain $X \equiv \MO_P(2) +2g^*L \sim 2(S+g^*L)$. 
This, together with $S|_S \sim -g^*L|_S$, implies $X|_S \equiv 0$. 
By $X \neq S$, we obtain $X \cap S = \emptyset$, as otherwise we could find a curve $C$ on $S$ which properly intersects $X$. 
Thus (b) holds. 
%Let us show (c). 
Then it holds that 
\[
g^*E|_S \sim (K_P+X)|_S \overset{{\rm (b)}}{\sim} K_P|_S \sim K_S -S|_S \sim g^*(K_Y+L)|_S, 
\]
which implies $E \sim K_Y+L$, i.e., (c) holds. 
This completes the proof of Claim. %$K_P+X \sim g^*(K_Y+L)$. 
\end{proof}

We can write $X -2S \sim g^*F$ for some Cartier divisor $F$ on $Y$. 
By $-2g^*L|_S \sim (X-2S)|_L \sim g^*F|_S$, we obtain $F\sim 2L$. 
Then $X \sim 2S +2g^*L \sim 2T$. 
This completes the proofs of (2) and (3). 

%We now show that $K_P+X \sim g^*(K_Y+L)$. 

% Therefore, 
% \[
% \MO_P(K_P) \otimes \MO_P(2) \otimes g^*\MO_Y(2L) \simeq 
% g^*\MO_Y(K_Y+L). 
% \]

Let us show (4). 
%Let us show $\Omega_{P/Y}^1 \simeq \MO_P(-g^*L -2S)$. 
%Since $g : P \to Y$ is a $\P^1$-bundle, 
We have an exact sequence 
\[
0 \to g^*\Omega_Y^1 \xrightarrow{\alpha} \Omega^1_P \to \Omega_{P/Y}^1 \to 0, 
\]
where the injectivity of $\alpha$ can be checked by taking 
the corresponding sequence of the stalks at the generic point. %exactness can be checked by taking the stalks 
Taking the wedge products, we get $\omega_P \simeq g^*\omega_Y \otimes  \Omega_{P/Y}^1$. 
By $K_P+X \sim g^*(K_Y+L)$ and $X \sim   2S +2g^*L$, we obtain 
\[
\Omega_{P/Y}^1 \simeq \MO_P(K_P -g^*K_Y) \simeq \MO_P(-X +g^*L) 
\simeq \MO_P(-g^*L-2S). 
\]
Thus (4) holds. 
\qedhere
% Since $f: X \to Y$ is a double cover, $X -2S \sim g^*E$ for some Cartier divisor $E$ on $Y$. 

% By $X \cap S = \emptyset$, 
% we obtain $2L \sim (X-2S)|_S \sim (g^*F)|_S$, which implies $F \sim 2L$. 
% Therefore, $X \sim 2S + 2g^*L \sim \MO_P(2)$. 
\end{proof}

\begin{lem}\label{l-Omega_P}
We use the same notation as Lemma \ref{l-dc-to-bdl}. 
Fix $q \in \Z$ and take a Cartier divisor $D$ on $Y$. 
Assume that 
\begin{enumerate}
\item $H^{q-1}(Y, D)=H^{q-1}(Y, D-L)=0$, and 
\item $H^q(Y, \Omega_Y^1(D))=H^q(Y, \Omega_Y^1(D-L))=0$. 
\end{enumerate}
Then $H^q(P, \Omega_P^1(g^*D))=0$. 
\end{lem}

\begin{proof}
%More generally, we get a criterion to get the vanishing of $H^k(P, \Omega_P^1(g^*D))=0$. 
We have an exact sequence 
\[
0 \to g^*(\Omega_Y^1(D)) \to \Omega_P^1(g^*D) \to \Omega^1_{P/Y}(g^*D) \to 0. 
\]
By (2), it is enough to show 
%Assuming $H^k(Y, \Omega_Y^1(D))=0$, the vanishing follows from 
 $H^q(P, \Omega^1_{P/Y}(g^*D))=0$, 
i.e., $H^q(P, g^*D-g^*L-2S)=0$. 
Using an exact sequence $0 \to \MO_P(-S) \to \MO_P \to \MO_S \to 0$ twice,
it suffices to prove that 
$H^q(P, g^*D-g^*L))=0$ and $H^{q-1}(S, g^*D-g^*L-nS)=0$ with $n \in \{0, 1\}$. 
By $-S|_S \sim g^*L|_S$, these equalities follow from (2) and (1), respectively. 
\end{proof}

\begin{prop}\label{p-2:1-QFS}
We use the same notation of Lemma \ref{l-dc-to-bdl}. 
Assume that $\dim X = \dim Y=3$, and both $L$ and $H\coloneqq -K_Y-L$ are ample. 
Consider the following conditions: 
\begin{enumerate}
\item[(0)]
\begin{enumerate}
\item[(0a)] $H^j(Y, -A))=0$ for every $j<3$ and every ample Cartier divisor $A$ on $Y$. 
\item[(0b)] $H^j(Y, \MO_Y(p^iH)) =H^j(Y, \MO_Y(p^iH -L)) =H^j(Y, \MO_Y(p^iH -2L))=0$ for every $j \in \{1, 2\}$. 
\item[(0c)] $H^3(Y, \sO_Y( p^iH -2L)) = H^3(Y, \sO_Y(p^iH-3L))=0$ for every $i>0$. 
\end{enumerate}
\item[(1)]
\begin{enumerate}
\item[(1a)] $H^1(Y, \Omega_Y^1(-mH-nL))=0$ for every $m\geq1$ and $n \geq 0$. 
\item[(1b)] $H^2(Y, \Omega_Y^1(-mH-2L))=H^2(Y, \Omega_Y^1(-mH-3L))=0$. 
\item[(1c)] $H^j(Y, \Omega^1_Y(p^iH)) =H^j(Y, \Omega^1_Y(p^iH-L))=0$ for every $j \in \{2, 3\}$ and every $i>0$. 
%\item $H^3(Y, \Omega_P^1(p^i$
\end{enumerate}
\item[(2)] 
$H^0(X, \Omega_X^2(p^iK_X))=0$ for every $i>0$. 
\end{enumerate}
Then the following hold. 
\begin{enumerate}
    \item[(I)] If (0a), (1a), and (1b) hold, then
    $H^1(X, \Omega^1_X({m}K_X))=0$ for all $m\geq 1$. 
    \item[(II)] If all the conditions above hold,
    $H^2(X, \Omega^1_X(-p^iK_X))=0$ for every $i>0$. 
    \item[(III)] If all the conditions above hold, $X$ is quasi-$F$-split. 
\end{enumerate}

\end{prop}

% \begin{enumerate}
% \item[$(*)$] $H^2(P, \Omega^{1}_{P}(p^ig^{*}H))=H^3(P, \Omega^{1}_{P}(p^ig^{*}H))=0$. 
% \end{enumerate}
% We have 
% \[
% 0 \to g^*\Omega_Y^1 \to \Omega_P^1 \to \Omega^1_{P/Y} \to 0. 
% \]
% Assuming $H^{2, 3}(Y, \Omega_Y^1(p^iH))=0$, ($*$) can be replaced by $H^{2, 3}(Y, \Omega^1_{P/Y}(p^ig^*H))=0$, 
% i.e., $H^{2, 3}(P, p^ig^*H-g^*L-2S)=0$. 
% Adding $S$ twice, enough to impose 
% $H^{2, 3}(P, p^ig^*H-g^*L)=0$ and $H^{1, 2}(S, p^ig^*H-g^*L-nS)=0$ wiht $n \in \{0, 1\}$. 
% These are $H^{2, 3}(Y, p^iH-L)=0$ and $H^{1, 2}(Y, p^iH-mL)=0$ for $m \in \{1, 2\}$ (by $-S|_S =L$). 

\begin{proof}
By Proposition \ref{prop:criterion for qFs}, (I) and (II) imply (III). 
In what follows, we shall prove (I) and (II). 
% it suffices to show that 
% \begin{enumerate}
%     \item[(I)] $H^1(X, \Omega^1_X(K_X))=0$ and 
%     \item[(II)] $H^2(X, \Omega^1_X(-p^iK_X))=0$ for every $i>0$
% \end{enumerate}
Note that we can write $-K_X \sim f^*H$  for $H := -K_Y-L$, and hence $-K_X$ is ample. 

\medskip

\noindent\textbf{Step 1: Proof of (I).}\,\,  
By the conormal exact sequence, we have an exact sequence
\[
0\to \sO_X(mK_X-X)\to \Omega^1_P|_X(mK_X) \to \Omega^1_X(mK_X)\to 0,
\]
where $\Omega^1_P|_X(mK_X) := (\Omega^1_P|_X) \otimes \MO_X(mK_X)$. 
It follows from $K_X=-g^{*}H|_X =-f^*H$, $X|_X=(X-2S)|_X=2g^{*}L|_X =2f^*L$, and $f_{*}\sO_X=\sO_Y\oplus \MO_Y(-L)$ 
that 
\begin{align*}
    H^2(X, \MO_X(mK_X-X))&=H^2(X, f^{*}\MO_Y(-(mH+2L)))\\
                 &=H^2(Y, \sO_Y(-mH-2L))\oplus H^2(Y, \sO_Y(-mH-3L))\\
                 &\overset{{\rm (0a)}}{=}0. 
\end{align*}
Thus it suffices to show $H^1(X, \Omega^1_P|_X(mK_X))=0$.

By $K_X=-g^{*}H|_X$, $S|_X=0$, and $-X=-2S-2g^{*}L$, we have the following exact sequence:
\[
0\to \Omega^{1}_{P}(-mg^{*}H-2g^{*}L)\to \Omega^{1}_{P}(-mg^{*}H+2S)\to \Omega^{1}_{P}|_X(mK_{X})\to 0.
\]
%We get $H^2(P, \Omega^{1}_{P}(-g^{*}H-2g^{*}L))=0$ by 
Thus, in order to prove (I), it is enough to show  
\begin{enumerate}
    \item[(i)] $H^1(P, \Omega^{1}_{P}(-mg^{*}H+2S))=0$ and 
    \item[(ii)] $H^2(P, \Omega^{1}_{P}(-mg^{*}H-2g^{*}L))=0$.
\end{enumerate}

\noindent\textbf{Step 1-1: Proof of (i).}\,\, 
%By the relative exact sequence, 
We have the following exact sequence:
\begin{multline*}
    0\to g^{*}(\Omega^1_Y(-mH))(2S) \to \Omega^1_P(-mg^{*}H+2S) \to \Omega^1_{P/Y}(-mg^{*}H+2S)%\simeq \sO_P(-g^{*}H-g^{*}L) 
    \to 0.
\end{multline*}
It holds that 
\[
H^1(P, \Omega^1_{P/Y}(-mg^{*}H+2S)) \simeq 
H^1(P, \sO_P(-mg^{*}H-g^{*}L)) \simeq %H^1(P, \sO_P(-g^{*}H-g^{*}L))\simeq 
H^1(Y, \sO_Y(-mH-L))\overset{{\rm (0a)}}{=}0. 
\]
%, where the last equality follows from Kodaira vanishing for $Y$.  
Then it suffices to show $H^1(P, g^{*}(\Omega^1_Y(-mH))(2S)) =0$. 
Using an exact sequence $0 \to \MO_P(-S) \to \MO_P \to \MO_S \to 0$ twice,  
the vanishing $H^1(P, g^{*}(\Omega_Y(-mH))(2S)) =0$ can be reduced, by  $S|_S=-g^{*}L|_S$, 
to those of 
\begin{itemize}
    \item $H^1(S, g^{*}(\Omega^1_Y(-mH))(nS))\overset{g|_S:\text{isom}}{\simeq} H^1(Y, \Omega^1_Y(-mH-nL))$ for $n\in\{1,2\}$ and
    \item $H^1(P, g^{*}\Omega_Y^1(-mH))\simeq H^1(Y, \Omega_Y^1(-mH))$. 
\end{itemize}
Both of them  follow from (1a). Thus (i) holds. 
%both of which follow from Bott vanishing for $Y$. 

\medskip

\noindent\textbf{Step 1-2: Proof of (ii).}\,\, 
In order to show (ii), it is enough to verify the assumptions of Lemma \ref{l-Omega_P} for the case when $q=2$ and $D = -H-2L$. 
The conditions Lemma \ref{l-Omega_P}(1) and Lemma \ref{l-Omega_P}(2) hold by (0a) and (1b), respectively. 
This completes the proofs of (ii) and (I).

% By the relative exact sequence, we have the following exact sequence:
% \begin{multline*}
% 0\to g^{*}\Omega^1_Y(-H-2L) \to \Omega^1_P(-g^{*}H-2g^{*}L)\\ \to \Omega^1_{P/Y}(-g^{*}H-2g^{*}L)\simeq \sO_P(-g^{*}H-3g^{*}L-2S) \to 0.
% \end{multline*}
% By Bott vanishing for $Y = \P^1 \times \P^2$, it holds that 
% \[
% H^2(P, g^{*}\Omega^1_Y(-H-2L)) \simeq H^2(Y, \Omega^1_Y(-H-2L))=0. 
% \]
% %by Bott vanishing for $Y$.
% It suffices to show $H^2(P, \sO_P(-g^{*}H-3g^{*}L-2S))=0$. 
% By $S|_S=-g^{*}L|_S$ and an exact sequence $0 \to \MO_P(-S) \to \MO_P \to \MO_S \to 0$,
% %the vanishing $H^2(P, \sO_P(-g^{*}H-3g^{*}L-2S))$ is  reduced to those of
% the problem is reduced to 
% \begin{itemize}
%     \item $H^1(S, \MO_P(-g^{*}H-3g^{*}L-nS)|_S)=0$ for $n\in\{0,1\}$ and 
%     \item $H^2(P, -g^{*}H-3g^{*}L)=0$.
% \end{itemize}
% These vanishings follows from Kodaira vanishing for $Y$, because  we have
% \begin{eqnarray*}
%     H^1(S, \MO_P(-g^{*}H-3g^{*}L-nS)|_S) &\simeq& H^1(S, -g^{*}H+(-3+n)g^{*}L)\\
%     &\overset{g|_S:\text{isom}}{\simeq}& H^1(Y, -H+(-3+n)L)
% \end{eqnarray*}
% and
% \[ 
% H^2(P, g^{*}(-H-3L))\simeq H^2(Y, -H-3L). 
% \]

\medskip

\noindent\textbf{Step 2: Proof of (II).}\,\,  
Fix $i \in \Z_{>0}$. 
By the conormal exact sequence, we have the following exact sequence:
\[
0 \to \sO_X(-p^iK_X-X) \to \Omega^1_P|_X(-p^iK_X) \to \Omega^1_X(-p^iK_X) \to 0.
\]
It follows from $K_X=-f^{*}H$, $X|_X=2f^{*}L$, and $f_{*}\sO_X=\sO_Y\oplus \MO_Y(-L)$ that 
\begin{align*}
H^3(X, -p^iK_X-X)&= H^3(X, p^if^{*}H-2f^{*}L)\\
&=H^3(Y, \sO_Y( p^iH -2L))\oplus H^3(Y, \sO_Y(p^iH-3L))\\
&\overset{{\rm (0c)}}{=}0. 
\end{align*}
%because $h^3(Y, p^iH-mL) = h^0(Y, )$
Thus it suffices to show $H^2(X, \Omega^1_P|_X(-p^iK_X))=0$.

We get 
\[
H^q(P,  \Omega^{1}_{P}(p^ig^{*}H))=0\qquad  \text{for} \qquad \text{every} \qquad q \in \{2, 3\}, 
\]
because Lemma \ref{l-Omega_P}, for the case when $q \in \{2, 3\}$ and $D=p^iH$, is applicable by 
(0b) and (1c). 
Since $K_X=-g^{*}H|_X$ and $-X=-2T$, we have the following exact sequence:
\[
0\to \Omega^{1}_{P}(p^ig^{*}H-2T)\to \Omega^{1}_{P}(p^ig^{*}H)\to \Omega^{1}_{P}|_X(-p^iK_{X})\to 0.
\]
Thus it suffices to show that
\[
H^3(P, \Omega^{1}_{P}(p^ig^{*}H-2T))=0.
\]
% By $2T=X=2g^{*}L+2S$, we have 
% \[
% H^3(P, \Omega^{1}_{P}(p^ig^{*}H-2g^{*}L-2S)=H^3(P, \Omega^{1}_{P}(p^ig^{*}H-2T)).
% \]
By $T|_T=g^{*}L|_T$ and an exact sequence $0 \to \MO_P(-T) \to \sO_P \to \sO_T \to 0$,
% the vanishing 
% \[
% H^2(P, \sO_P(-g^{*}H-g^{*}L-2T))\] 
% can be reduced to those of
the problem is reduced to 
\begin{itemize}
    \item $H^2(T, \Omega^1_P|_T(p^ig^{*}H-nT))=0$ for $n\in\{0,1\}$ and 
    \item $H^3(P, \Omega^{1}_{P}(p^ig^{*}H))=0$. %$H^2(T, \Omega^1_P|_T(pg^{*}H))=0$.
\end{itemize}
The second vanishing has been settled already.
Thus it suffices to show the first one. 
Fix $n \in \{0, 1\}$. 
 By the conormal exact sequence, we have the following exact sequence
    \[
    0\to \sO_T(p^ig^{*}H-(n+1)T)\to \Omega^1_{P}|_T(p^ig^{*}H-nT) \to \Omega^1_{T}(p^ig^{*}H-nT)\to 0.
    \]
    It follows  from $\sO_T(p^ig^{*}H-(n+1)T) \overset{T \simeq Y}{\simeq} \sO_Y(p^iH-(n+1)L)$ that
    \begin{align*}
        H^2(T, \sO_T(p^ig^{*}H-(n+1)T) \simeq H^2(Y, \MO_Y(p^iH -(n+1)L)) \overset{{\rm (0b)}}{=}0. 
        % &\simeq H^0(\P^1, \sO_{\P^1}(p^i-n-1))\otimes H^2(\P^2, \sO_{\P^2}(p^i-2n-2))\\
        %                             &\oplus H^1(\P^1, \sO_{\P^1}(p^i-n-1))\otimes H^1(\P^2, \sO_{\P^2}(p^i-2n-2))\\
        %                             &=0
    \end{align*}
%    for $n\in\{0,1\}$. 
Then we are done by  $H^2(T, \Omega^1_{T}(p^ig^{*}H-nT) \simeq H^2(Y, \Omega^1_Y(p^iH -nL))\overset{{\rm (1c)}}{=}0$. 
    % Since we have 
    % \begin{align*}
    % \Omega^1_{T}(p^ig^{*}H-nT)\simeq \Omega^1_{Y}(p^iH-nL)&=\Omega^1_Y(p^i-n, p^i-2n)\\
    % &\simeq \mathrm{pr}_1^{*}\Omega^1_{\P^1}(p^i-n)\otimes \mathrm{pr}_2^{*}\sO_{\P^2}(p^i-2n)\\
    % &\oplus 
    % \mathrm{pr}_1^{*}\sO_{\P^1}(p^i-n)\otimes \mathrm{pr}_2^{*}\Omega^1_{\P^2}(p^i-2n)
    % \end{align*}
    % it follows {\cred from \cite[Proposition 1.3]{Totaro(Fano)}} that
    % \begin{align*}
    %     H^2(T, \Omega^1_{T}(p^ig^{*}H-nT))&\simeq H^0(\P^1, \Omega^1_{\P^1}(p^i-n))\otimes H^2(\P^2, \sO_{\P^2}(p^i-2n))\\
    %                                      &\oplus H^1(\P^1, \Omega^1_{\P^1}(p^i-n))\otimes H^1(\P^2, \sO_{\P^2}(p^i-2n))\\
    %                                    &\oplus  H^0(\P^1, \sO_{\P^1}(p^i-n))\otimes H^2(\P^2, \Omega^1_{\P^2}(p^i-2n))\\
    %                                 &\oplus H^1(\P^1, \sO_{\P^1}(p^i-n))\otimes H^1(\P^2, \Omega^1_{\P^2}(p^i-2n))\\
    %                                 &=0
    % \end{align*}
    % for $n\in\{0,1\}$.
    % Therefore, we obtain the desired vanishing.
    \qedhere

\end{proof}

\subsubsection{2-2}

\begin{lem}\label{lem:description of (2-2)}
A smooth Fano threefold $X$ of No.~2-2 satisfies the following properties:
    \begin{enumerate}
        \item There is a split double cover $f\colon X\to Y :=\P^1\times \P^2$. 
        \item $f_*\MO_X \simeq \MO_Y \oplus \MO_Y(-L)$ for a Cartier divisor $L$ satisfying $\MO_Y(L) \simeq \MO_Y(1, 2)$. 
        \item $X$ is (isomorphic to) a divisor on $P\coloneqq \mathbb{P}_{Y}(\sO_Y\oplus \MO_Y(-L))$.
        \item $K_P+X=-g^{*}H$ and $K_X=g^{*}(K_Y+L)=-g^{*}H|_{X}$, where $H\coloneqq \sO_Y(1,1)$ and $g: P =\mathbb{P}_{Y}(\sO_Y \oplus \MO_Y(-L)) \to Y$ denotes the projection.
\item 
 There exists a section $S$ of $g$ such that  $S \cap X = \emptyset$, $S|_S \sim -g^*L|_S$, $X -2S \sim 2g^*L$, 
$\MO_P(1) \simeq \MO_P(S)$, and $\Omega_{P/Y}^1 \simeq \MO_P(-g^*L -2S)$. 
\item 
There exists a section $T$ of $g$ such that 
$S \cap T =\emptyset$, 
$T|_T \sim g^*L|_T$, $X \sim 2T$, $T-S \sim g^*L$, and $\MO_P(1) \simeq \MO_P(T-g^*L)$. 
    \end{enumerate}
\end{lem}

\begin{proof}
The assertions (1) and (2) follow from \cite[Subsection 9.2]{FanoIII}. 
Then the remaining ones hold by Lemma \ref{l-dc-to-bdl}. 
\end{proof}

\begin{lem}\label{lem:qFs for 2-2}
Let $X$ be a smooth Fano threefold of No.~2-2. 
Then the following hold. 
\begin{enumerate}
    \item $H^1(X, \Omega^1_X(K_X))=0$. 
    \item $H^2(X, \Omega^1_X(-p^iK_X))=0$ for every $i>0$. 
    \item $X$ is quasi-$F$-split. 
\end{enumerate}
\end{lem}

\begin{proof}
We use the same notation as Lemma \ref{lem:description of (2-2)}. 
It is enough to verify the conditions of Proposition \ref{p-2:1-QFS} for 
$Y  = \P^1 \times \P^2, L = \MO_Y(1, 1)$, and $H =\MO_Y(1, 2)$. 
Note that $mH-L$ is ample when $m \geq 2$. 
By Lemma \ref{lem:SRC}, Proposition \ref{p-2:1-QFS}(2) holds. 
Since $Y$ satisfies Kodaira vanishing, 
it is easy to see that Proposition \ref{p-2:1-QFS}(0) holds. 
As $Y$ satisfies Bott vanishing, 
it is obvious that  Proposition \ref{p-2:1-QFS}(1) holds. 
\end{proof}

\subsubsection{2-6-a}

\begin{definition}\label{description of 2-6 (a)}
 Let $X$ be a smooth Fano threefold of No.~2-6. 
By \cite[Section 7.2]{FanoIV}, one of the following holds up to isomorphisms. 
 \begin{enumerate}
\item[(2-6-a)] $X$ is a  hypersurface of $P:=\P^2\times \P^2$ of bidegree $(2,2)$. In this case, 
$\MO_X(-K_X) \simeq \MO_X(1, 1)$, where $\MO_X(1, 1) := \MO_P(1, 1)$. In this case, we say that $X$ is (a Fano threefold) of {\em No. 2-6-a}. 
\item[(2-6-b)] 
There is a split double cover $f: X \to W$ satisfying $f_*\MO_X \simeq \MO_W \oplus \MO_W(-L)$, 
where $L$ is a Cartier divisor on $W$ with $\MO_W(2L) \simeq \omega_W^{-1}$. 
In this case, we  say that $X$ is (a Fano threefold) of {\em No. 2-6-b}. 
 \end{enumerate}
\end{definition}

\begin{lem}\label{lem:vanishing for (2-6-1)}
    Let $X$ be a smooth Fano threefold of No.~2-6-a.
    Then the following hold.
    \begin{enumerate}
        \item $H^i(X, \sO_X(n,n))=0$ for $i\in\{1,2\}$ and $n\in\Z$.
        \item $H^3(X, \sO_X(n,n))=0$ for $n\in\Z_{\geq 0}$.
    \end{enumerate}
\end{lem}
\begin{proof}
Use an exact sequence $0 \to \MO_{\P^2 \times \P^2}(-2, -2) \to \MO_{\P^2 \times \P^2} \to \MO_X \to 0$. 
\end{proof}

\begin{lem}\label{lem:qFs for 2-6a}
Let $X$ be a smooth Fano threefold of No.~2-6-a. 
Then the following hold. 
\begin{enumerate}
    \item $H^1(X, \Omega^1_X(K_X))=0$. 
    \item $H^2(X, \Omega^1_X(-p^iK_X))=0$ for every $i>0$. 
    \item $X$ is quasi-$F$-split. 
\end{enumerate}
\end{lem}
\begin{proof}
We use the notation of Definition \ref{description of 2-6 (a)}.
By Proposition \ref{prop:criterion for qFs} and Lemma \ref{lem:SRC}, it suffices to show (1) and (2). 

\medskip

\noindent\textbf{Step 1: Proof of (1).}\,\,  By the conormal exact sequence, we have the following exact sequence:
\[
0\to \sO_X(K_X-X)\to \Omega^1_P|_X(K_X) \to \Omega^1_X(K_X)\to 0.
\]
Since we have 
\[
H^2(X, K_X-X)=H^2(X, \sO_X(-3,-3))=0
\] 
by Lemma \ref{lem:vanishing for (2-6-1)},
it suffices to show $H^1(X, \Omega^1_P|_X(K_X))=0$.

We have the following exact sequence:
\[
0\to \Omega^{1}_{P}(K_{P})\to \Omega^{1}_{P}(K_{P}+X)\to \Omega^{1}_{P}|_X(K_{X})\to 0
\]
By Bott vanishing, we have
\begin{itemize}
    \item $H^1(P, \Omega^{1}_{P}(K_{P}+X))=H^1(P, \Omega^{1}_{P}(-1,-1))=0$ and 
    \item $H^2(P, \Omega^{1}_{P}(K_{P}))=H^2(P, \Omega^{1}_{P}(-3,-3))=0$.
\end{itemize}
Therefore, (1) holds.

\medskip

\textbf{Step 2: Proof of (2).}\,\,  
Fix $i \in \Z_{>0}$. 
By the conormal exact sequence, we have the following exact sequence
\[
0 \to \sO_X(-p^iK_X-X) \to \Omega^1_P|_X(-p^iK_X) \to \Omega^1_X(-p^iK_X) \to 0.
\]
We have 
\[
H^3(X, \sO_X(-p^iK_X-X))=H^3(X, \sO_X(p^i-2,p^i-2))=0
\]
by Lemma \ref{lem:vanishing for (2-6-1)}.
Thus, it suffices to show $H^2(X, \Omega^1_P|_X(-p^iK_X))=0$.

We have the following exact sequence:
\[
0\to \Omega^{1}_{P}(-p^i(K_{P}+X)-X)\to \Omega^{1}_{P}(-p^i(K_{P}+X))\to \Omega^{1}_{P}|_X(-p^iK_{X})\to 0.
\]
Then 
\begin{itemize}
    \item $H^2(P, \Omega^{1}_{P}(-p^i(K_{P}+X)))=H^2(P, \Omega^{1}_{P}(p^i,p^i))=0$ by Bott vanishing and 
    \item $H^3(P, \Omega^{1}_{P}(-p^i(K_{P}+X)-X))=H^3(P, \Omega^{1}_{P}(p^i-2,p^i-2))=0$ by \cite[Proposition 1.3]{Totaro(Fano)}.
\end{itemize}
Therefore, the assertion holds.  
\end{proof}

\subsubsection{2-6-b}

\begin{lem}\label{lem:description of 2-6 (b)}
A smooth Fano threefold $X$ of No.~2-6-b satisfies the following properties:
    \begin{enumerate}
        \item There is a split double cover $f\colon X\to W$, where 
        $W$ is a smooth hypersurface of $\mathbb{P}^2\times \mathbb{P}^2$ of bidegree $(1,1)$.
        \item $f_*\MO_X \simeq \MO_W \oplus \MO_W(-L)$ 
        for a Cartier divisor $L=\MO_W(1, 1)$.
        \item $X$ is a divisor of $P\coloneqq \mathbb{P}_{W}(\sO_W\oplus \MO_W(-L))$.
        \item $K_P+X=-g^{*}L$, $K_X=g^{*}(K_W+L)=-g^{*}L|_{X}$.
        \item There exists a section $S$ of $g$ such that  $S \cap X = \emptyset$, $S|_S \sim -g^*L|_S$, $X -2S \sim 2g^*L$, 
$\MO_P(1) \simeq \MO_P(S)$, and $\Omega_{P/W}^1 \simeq \MO_P(-g^*L -2S)$. 
    \end{enumerate}
\end{lem}

\begin{proof}
The assertions (1) and (2) follow from \cite[Subsection 9.2]{FanoIII}. 
Then the remaining ones hold by Lemma \ref{l-dc-to-bdl}. 
\end{proof}

\begin{lem}\label{lem:vanishing of W for (2-6-2)}
    Let $Y=\P^3$ (resp.~$Y=Q$, resp. $Y=W$), where $Q$ is a smooth quadric hypersurface of $\P^4$ and $W$ is a smooth hypersurface of $\mathbb{P}^2\times \mathbb{P}^2$ of bidegree $(1,1)$.
    Let $L=\sO_{\P^3}(3)$ (resp.~$\sO_Q(2)$, resp. $\sO_{W}(1,1)$).
    Then the following hold.
    \begin{enumerate}
         \item $H^1(Y, \sO_Y(nL))=0$ for $n\in \Z$.
         \item $H^2(Y, \sO_Y(nL))=0$ for $n\in \Z$.
         \item $H^3(Y, \sO_Y(nL))=0$ for $n\in\Z_{\geq -1}$.
         \item $H^0(Y, \Omega^1_Y(nL))=0$ for $n\in \Z_{\leq 0}$.
        \item $H^1(Y, \Omega^1_Y(nL))=0$ for $n\in \Z\setminus\{0\}$.
        \item $H^2(Y, \Omega^1_Y(nL))=0$ for 
        $n\in \Z \setminus \{-1\}$. 
        \item $H^3(Y, \Omega^1_Y(nL))=0$ for $n\in\Z_{\geq 0}$.
    \end{enumerate}
\end{lem}

\begin{proof}
Since $Y$ is $F$-split (Lemma \ref{l-prim-toric}), 
(1)-(3) hold. 
The assertions (4) and (7) follow from the fact that $X$ is SRC (Lemma \ref{lem:SRC}). 
Let us  prove (5) and (6). 
If $Y=\P^3$, then these follow from the Bott vanishing theorem and 
\cite[Proposition 1.3]{Totaro(Fano)}. 
In what follows, we assume $Y \in \{Q, W\}$. 
If $Y = Q$ (resp. $Y=W$), then 
\begin{itemize}
    \item we have an embedding $Y \subset P$ 
for $P:= \P^4$ (resp. $P=\P^2 \times \P^2$), 
\item we set $H := \MO_{\P^4}(1)$ (resp. $H :=\MO_{\P^2 \times \P^2}(1, 1)$), and 
\item we get $L \sim sH|_Y$ for $s:= 2$ (resp. $s:=1$). It holds that $Y \sim sH$. 
\end{itemize}
We have the following exact sequence:
\[
0\to \Omega_P^1((n-s)H)\to \Omega_{P}^1(nH) \to \Omega^1_{P}(nH)|_Y \to 0.
\]
By Bott vanishing and \cite[Proposition 1.3]{Totaro(Fano)}, we have 
\begin{itemize}
         \item $H^1(P, \Omega^1_P(nH))=0$ for $n\in \Z\setminus\{0\}$, 
        \item $H^2(P, \Omega^1_P(nH))=0$ for $n\in \Z$, and 
        \item $H^3(P, \Omega^1_P(nH))=0$ for $n \in \Z$. 
\end{itemize}
We then get  
\begin{itemize}
         \item $H^1(Y, \Omega^1_P(nH)|_Y)=0$ for $n\in \Z\setminus\{0\}$ and
        \item $H^2(Y, \Omega^1_P(nH)|_Y)=0$ for $n \in \Z$.
\end{itemize}
By $Y \sim sH$ and the conormal exact sequence, we have the following exact sequence:
\[
0\to \sO_Y(nH - sH)\to \Omega_{P}^1(nH)|_Y \to \Omega^1_{Y}(nH)\to 0.
\]
Recall that 
\begin{itemize}
         \item $H^2(Y, \sO_Y(nH))=0$ for $n\in \Z$ and
        \item $H^3(Y, \sO_Y(nH))=0$ for $n\geq -s$. 
\end{itemize}
Hence we get 
\begin{enumerate}
         \item[\textup{(5)'}] $H^1(Y, \Omega^1_Y(nH))=0$ for $n\neq 0$ and
        \item[\textup{(6-a)}] $H^2(Y, \Omega^1_Y(nH))=0$ for $n\geq 0$.
\end{enumerate}
In particular, (5) holds. %is obtained. 
Comparing (6) with (6-a), it suffices to show  (6-b) below by Serre duality. 
\begin{enumerate}
        \item[\textup{(6-b)}] $H^1(Y, \Omega^2_Y(nH))=0$ for $n >  s$.  
\end{enumerate}
Taking the wedge product $\bigwedge^2$ to the conormal exact sequence, 
we get the following exact sequence:
\[
0\to \Omega^1_Y(nH-sH) \to \Omega^2_{P}(nH)|_Y \to \Omega^2_{Y}(nH)\to 0.
\]
Since we have $H^2(Y, \Omega_Y^1(nH-sH))=0$ for $ n >%{\cred \geq} 
s$ (6-a), 
it is enough to prove $H^1(Y, \Omega^2_{P}(nH)|_Y)=0$ for $n 
> %{\cred \geq}  
s$. 
This holds by an exact sequence 
\[
0 \to \Omega^2_P(nH-sH) \to \Omega^2_P(nH) \to \Omega^2_P(nH)|_Y \to 0,
\]
because Bott vanishing implies 
$H^i(P, \Omega_P^2(mH))=0$ for $i>0$ and $m>0$. 
\end{proof}

\begin{lem}\label{lem:qFs for 2-6b}
Let $X$ be a smooth Fano threefold of No.~2-6-b. 
Then the following hold. 
\begin{enumerate}
    \item $H^1(X, \Omega^1_X(K_X))=0$. 
    \item $H^2(X, \Omega^1_X(-p^iK_X))=0$ for every $i>0$. 
    \item $X$ is quasi-$F$-split. 
\end{enumerate}
\end{lem}

\begin{proof}
It is enough to verify the conditions in Proposition \ref{p-2:1-QFS}. 
Note that we have $H = -K_Y -L =L$. 
Then the conditions in Proposition \ref{p-2:1-QFS} follow from Lemma \ref{lem:SRC} and Lemma \ref{lem:vanishing of W for (2-6-2)}. 
\end{proof}

The above argument can be applied for some hyperelliptic Fano threefolds. 
Let us start by recalling the definition.

\begin{dfn}
We say that a smooth Fano threefold $X$ is {\em hyperelliptic} if $X$ is of index one, 
$|-K_X|$ is base point free, and 
the induced morphism $f\colon X \to Y \coloneqq \varphi_{|-K_X|}(X)$ is a double cover. 
\end{dfn}

It is known that if $\rho(X)=1$, then $Y$ is isomorphic to $\P^3$ or $Q$ in the above notation \cite[Theorem 6.5]{FanoI}. 
The assumption $p>5$ in Proposition \ref{p-hyperell-QFS}(i) is sharp as we shall see later (Example \ref{e-p=5-nonQFS}).

\begin{prop}\label{p-hyperell-QFS}
Let $X$ be a hyperelliptic smooth Fano threefold such that $\rho(X)=1$. 
Let $f\colon X \to Y \coloneqq \varphi_{|-K_X|}(X)$ be the double cover induced by $\varphi_{|-K_X|}$. 
Assume the following. 
\begin{enumerate}
\item[(i)] If $Y \simeq  \P^3$, then $p>5$. 
\item[(ii)] If $Y \simeq Q$, then $p>3$. 
\end{enumerate}
Then the following hold. 
\begin{enumerate}
    \item $H^1(X, \Omega^1_X(K_X))=0$. 
    \item $H^2(X, \Omega^1_X(-p^iK_X))=0$ for every $i>0$. 
    \item $X$ is quasi-$F$-split. 
\end{enumerate}
\end{prop}

\begin{proof}
We have $f_*\MO_X \simeq \MO_Y \oplus \MO_Y(-L)$ for a Cartier divisor $L=\MO_Y(r-1)$, 
where $r$ denotes the index of $Y$ (i.e., if $Y = \P^3$ (resp.~$Y=Q$), then $r=4$ (resp.~$r=3$)). 
It is enough to verify the conditions in Proposition \ref{p-2:1-QFS}. 
Note that we have $L = (r-1)H$ and $\MO_Y(1) = \MO_Y(H)$. 
Then Lemma \ref{lem:SRC} and Lemma \ref{lem:vanishing of W for (2-6-2)} imply all the conditions in Proposition \ref{p-2:1-QFS} except for Proposition \ref{p-2:1-QFS}(0c). 
By our assumptions (i) and (ii), Proposition \ref{p-2:1-QFS}(0c) directly follows from Serre duality, e.g., 
\[
h^3(Y, \MO_Y(p^iH-3L)) =h^0(Y, \MO_Y(K_Y+3L-p^iH)) = h^0(Y, \MO_Y((2r-3-p^i)H))=0
\]
for every $i>0$.
\qedhere

\end{proof}

\subsubsection{2-8}

\begin{lem}\label{lem:description of (2-8)}
A smooth Fano threefold $X$ of No.~2-8 satisfies the following properties:
    \begin{enumerate}
        \item $X$ is (isomorphic to) a divisor on $P\coloneqq \mathbb{P}_{\mathbb{P}^2}(\sO_{\P^2}\oplus \sO_{\P^2}(1) \oplus \sO_{\P^2}(2))$.
        \item $\sO_P(K_P)=\MO_P(-3)$, $\sO_P(X)=\MO_P(2)$, $K_X=\sO_X(-1),$ and $\sO_X(X)=\sO_X(2)=\sO_X(-2K_X)$.
        \item $|\sO_P(1)|$ is base point free and $\MO_P(1)$ is big. 
        \item Let $\varphi : P \to P'$ be the birational morphism to a normal projective variety $P'$ 
        such that $\varphi_*\MO_P = \MO_{P'}$ 
        and $\MO_P(1) \simeq \varphi^*\MO_{P'}(1)$ 
        for some ample invertible sheaf $\MO_{P'}(1)$ on $P'$. 
        Then $\varphi : P \to P'$ is a small birational morphism which is an isomorphism around $X$.   
    \end{enumerate}
\end{lem}

We say that a birational morphism  $\varphi : P \to P'$ is {\em small} if $\dim \Ex(\varphi) \leq \dim P -2$. 

\begin{proof}
Note that a Fano threefold of No. 2-8 is characterised by the following properties (i) and (ii) \cite[Theorem 5.34]{FanoIII}:
\begin{enumerate}
\item[(i)] $X$ is a Fano threefold with $\rho(X)=2$. 
\item[(ii)] One of the extremal rays is of type $C_1$, and the the other extremal ray is of type $E_3$ or $E_4$. 
\end{enumerate}
Then we may apply \cite[Proposition 5.29 and Lemma 5.30]{FanoIII}. 
By \cite[Proposition 5.29(3), Lemma 5.30]{FanoIII}, 
$X$ is a divisor on $P\coloneqq \mathbb{P}_{\mathbb{P}^2}(\sO_{\P^2}\oplus \sO_{\P^2}(1) \oplus \sO_{\P^2}(2))$ satisfying $X \sim \MO_P(2)$. Thus (1) holds. 
Moreover, \cite[the proof of Lemma 5.30]{FanoIII} implies that 
$-K_P \sim \MO_P(3)$. 
Then the adunction formula implies $K_X \sim (K_P+X)|_X \sim \MO_P(-3+2)|_X = \MO_X(-1)$. 
Thus (2) holds. 

Let us show (3) and (4). 
We have three sections $\Gamma_0, \Gamma_1, \Gamma_2$ of 
the induced $\P^2$-bundle $g: P = \mathbb{P}_{\mathbb{P}^2}(E) \to \P^2$, where 
$E :=\sO_{\P^2}\oplus \sO_{\P^2}(1) \oplus \sO_{\P^2}(2)$. 
corresponding to the projections of $E$ 
to the factors $\MO_{\P^2}, \MO_{\P^2}(1), \MO_{\P^2}(2)$, respectively. 
Similarly, we have the following three prime divisors which are $\P^1$-bundles over $\P^2$: 
\begin{itemize}
\item $D_0 := \P_{\P^2}(\MO_{\P^2}(1) \oplus \MO_{\P^2}(2))$, corresponding to $E \to \MO_{\P^2}(1) \oplus \MO_{\P^2}(2)$. 
\item $D_1 := \P_{\P^2}(\MO_{\P^2} \oplus \MO_{\P^2}(2))$, corresponding to $E \to \MO_{\P^2} \oplus \MO_{\P^2}(2)$. 
\item $D_2 := \P_{\P^2}(\MO_{\P^2} \oplus \MO_{\P^2}(1))$, corresponding to $E \to \MO_{\P^2} \oplus \MO_{\P^2}(1)$. 
\end{itemize}
By construction, we have $D_i \cap \Gamma_i = \emptyset$ for every $i \in \{ 0, 1, 2\}$. 
Fix a line $L$ on $\P^2$ and set $F := g^*L$, which is a prime divisor on $P$.  
For each $i \in \{0, 1, 2\}$, we can write $\MO_P(1) \sim D_i + bF$ for some $b \in \Z$. 
It holds that 
\[
\MO_{\P^2}(i) = \MO_P(1)|_{\Gamma_i} = (D_i +bF)|_{\Gamma_i} = bF|_{\Gamma_i} \simeq \MO_{\P^2}(b), 
\]
i.e., $b=i$. Then 
%$\MO_P(1) \sim D_i +iF$ for each $i \in$
we have that 
\[
\MO_P(1) \sim D_0  \sim D_1 +F \sim D_2 + 2 F. 
\]
By $D_0 \cap D_1 \cap D_2 = \emptyset$, 
$|\MO_P(1)|$ is base point free. 
% We have $D_0 \cdot D_1 = \Gamma_2, D_1 \cdot D_2 = \Gamma_0, D_2 \cdot D_0 = \Gamma_1$, which imply  
% \[
% \MO_P(1)^4 = D_0 \cdot (D_1+F) \cdot (D_2+F)
% \]
Since $D_1$ is relatively ample over $\P^2$ and $F$ is the pullback of an ample divisor, 
the divisor $\MO_P(1) \sim D_1 +F$ is big. 
Thus (3) holds. 
Let $\varphi : P \to P'$ be as in the statement of (4). 

We now show that $\Ex(\varphi) = \Gamma_0$. 
By $\MO_P(1)|_{\Gamma_0} \simeq \MO_{\Gamma_0}$, 
$\varphi(\Gamma_0)$ is a point. 
In particular, $\Ex(\varphi) \supset \Gamma_0$. 
Pick a curve $C$ such that $\varphi(C)$ is a point. 
It suffices to show $C \subset \Gamma_0$. 
We have $\MO_P(1) \cdot C=0$. 
Note that $g(C)$ is not a point, because $\MO_P(1)$ is $g$-ample. 
In particular, $F \cdot C >0$. 
Therefore, 
\[
0 = \MO_P(1) \cdot C = (D_1+F) \cdot C =(D_2 + 2F) \cdot C. 
\]
By $F \cdot C >0$, we obtain $D_1 \cdot C <0$ and $D_2 \cdot C <0$. 
Hence $C \subset D_1 \cap D_2 = \Gamma_0$, as required.

It is enough to prove $X \cap \Gamma_0 =\emptyset$. 
Suppose $X \cap \Gamma_0 \neq \emptyset$. 
By $\MO_P(X)|_{\Gamma_0} \simeq \MO_P(2)|_{\Gamma_0} \simeq \MO_{\Gamma_0}$, 
we obtain $\P^2 \simeq \Gamma_0 \subset X$. 
Then the Stein factorisation $\psi : X \to X'$ of the composite morphism 
$\varphi|_X : X \hookrightarrow P \to P'$ is a birational morphism which contracts $\Gamma_0 \simeq \P^2$ to a point. 
This is absurd, because $X$ has no extremal ray of $E_2$ or $E_5$. Thus (4) holds. 
\end{proof}

\begin{lem}\label{lem:qFs for (2-8)}
    A smooth Fano threefold of No.~2-8 is quasi-F-split. 
\end{lem}
\begin{proof}
We use the notation of Lemma \ref{lem:description of (2-8)}.
By Proposition \ref{prop:criterion for qFs} and Lemma \ref{lem:SRC}, it suffices to show that 
\begin{enumerate}
    \item[(1)] $H^1(X, \Omega^1_X(K_X))=0$.
    \item[(2)] $H^2(X, \Omega^1_X(-p^iK_X))=0$ for every $i >0$. 
\end{enumerate}
\noindent\textbf{Step 1: Proof of (1).}\,\,
By the conormal exact sequence, we have the following exact sequence
\[
0\to \sO_X(K_X-X) \to \Omega^1_P|_X(K_X) \to \Omega^1_X(K_X)\to 0.
\]
Since $\sO_X(K_X)=\sO_X(-1)$ and $\sO_X(X)=\sO_X(2)$, we have an exact sequence
\[
H^2(P, \MO_P(-3))\to H^2(X, \sO_X(-3)) (= H^2(X, K_X-X))\to H^3(P, \MO_P(-5)).
\]
Since $P$ is toric and $\MO_P(1)$ is nef, 
it follows from $K_P \sim \MO_P(-3)$ and \cite[Proposition 1.3]{Totaro(Fano)} that
\begin{itemize}
    \item $H^2(P, \sO_P(-3))\simeq H^2(P, \sO_P)=0$ and 
    \item $H^3(P, \sO_P(-5))\simeq H^1(P, \sO_P(2))=0$.
\end{itemize}
Thus $H^2(X, K_X-X)=0$, and it suffices to show $H^1(X, \Omega^1_P|_X(K_X))=0$.

Since we have a closed embedding $X\subset P'$ around which $P'$ is smooth (Proposition \ref{lem:description of (2-8)}(4)), 
we have the following exact sequence:
\[
0\to \Omega^{[1]}_{P'}(K_{P'})\to \Omega^{[1]}_{P'}(K_{P'}+X)\to \Omega^{[1]}_{P'}|_X(K_{X})\simeq \Omega^{1}_{P}|_X(K_{X})\to 0,
\]
By Bott vanishing \cite[Theorem 1.1 or Corollary 1.3]{Fuj07}, we have
\begin{itemize}
    \item $H^1(P', \Omega^{[1]}_{P'}(K_{P'}+X))=H^1(P', \Omega^{[1]}_{P'}(-1))=0$ and 
    \item $H^2(P', \Omega^{[1]}_{P'}(K_{P'}))=H^2(P', \Omega^{[1]}_{P'}(-3))=0$.
\end{itemize}
This completes the proof of (1). 
%Therefore, we obtain the assertion.\\

\medskip

\noindent\textbf{Step 2: Proof of (2).}\,\,
By the conormal exact sequence, we have the following exact sequence
\[
0\to \sO_X(-p^iK_X-X) \to \Omega^1_P|_X(-p^iK_X) \to \Omega^1_X(-p^iK_X)\to 0.
\]
We have 
\[
H^3(X, \sO_X(-p^iK_X-X))=H^3(X, \sO_X(p^i-2))=0.
\]
Thus it suffices to show that $H^2(X, \Omega^1_P|_X(-p^iK_X))=0$.

We have the following exact sequence:
\[
0\to \Omega^{1}_{P}(-p^i(K_{P}+X)-X)\to \Omega^{1}_{P}(-p^i(K_{P}+X))\to \Omega^{1}_{P}|_X(-p^iK_{X})\to 0.
\]
Since $\sO_P(1)$ is nef and $P$ is a smooth toric variety, we have
\begin{enumerate}
    \item $H^2(P, \Omega^{1}_{P}(-p^i(K_{P}+X)))=H^2(P, \Omega^{1}_{P}(p^i))=0$ and
    \item $H^3(P, \Omega^{1}_{P}(-p^i(K_{P}+X)-X))=H^3(P, \Omega^{1}_{P}(p^i-2))=0$.
\end{enumerate}
by \cite[Proposition 1.3]{Totaro(Fano)}. Thus (2) holds. 
\end{proof}

\subsubsection{3-10}

\begin{lem}\label{description of (3.10)}
Let $X$ be a smooth Fano threefold $X$ of No.~3-10 such that there is a wild conic bundle structure 
$f: X \to \P^1 \times \P^1$. 
Then  the following properties hold:
    \begin{enumerate}
        \item $X$ is (isomorphic to) a divisor on $P \coloneqq\P_{\P^1 \times \P^1}(\MO \oplus \MO(1, 0) \oplus \MO(0, 1))$ satisfying $\MO_P(X) \sim \MO_P(2)$. 
        \item Each of $X$ and $\sO_P(1)$ is nef and big. 
        \item $-K_P$ and $-(K_P+X)$ are ample.
        \item $\MO_X(K_X) \simeq \sO_X(-1)\otimes f^{*}\sO_{\P^1\times \P^1}(-1,-1)$, where $\MO_X(-1) \coloneqq \MO_P(-1)|_X$.
        \item $-2^i(K_{P}+X)-X$ is nef for every $i>0$.
    \end{enumerate}
\end{lem}

\begin{proof}
The assertion (1) follows from \cite[Corollary 8]{MS03}. 
Let us show (2). 
By $X \sim \MO_P(2)$, it suffices to show that $|\MO_P(1)|$ is base point free and $\MO_P(1)$ is big. 
Set $L_0 := \MO_{\P^1 \times \P^1}, L_1 :=  \MO_{\P^1 \times \P^1}(1, 0), L_2 := \MO_{\P^1 \times \P^1}(0, 1)$, and 
\[
E \coloneqq L_0 \oplus L_1 \oplus L_2 = \MO \oplus \MO(1, 0) \oplus \MO(0, 1). 
\]
We have three sections of $\pi \colon P= \P_{\P^1 \times \P^1}(\MO \oplus \MO(1, 0) \oplus \MO(0, 1)) \to \P^1 \times \P^1$: 
\begin{itemize}
\item Set $\Gamma_0 := \P(\MO)$, which is corresponding to the projection   $E = \MO \oplus \MO(1, 0) \oplus \MO(0, 1) \to \MO$. 
We get $\MO_P(1)|_{\Gamma_0} = \MO = L_0$. 
\item Set $\Gamma_1 := \P(\MO(1, 0))$, which is corresponding to the projection   $E = \MO \oplus \MO(1, 0) \oplus \MO(0, 1) \to \MO(1, 0)$. 
We get $\MO_P(1)|_{\Gamma_1} = \MO(1, 0) = L_1$. 
\item Set $\Gamma_2 := \P(\MO(0, 1))$, which is corresponding to the projection   $E = \MO \oplus \MO(1, 0) \oplus \MO(0, 1) \to \MO(0, 2)$. 
We get $\MO_P(1)|_{\Gamma_2} = \MO(0, 1) = L_2$. 
\end{itemize}
Similarly, we have three prime divisors on $P$ which are $\P^1$-bundles over $\P^1 \times \P^1$: 
\begin{itemize}
\item Set $D_0 := \P(L_1 \oplus L_2)$, which is corresponding to the projection 
$E \to L_1 \oplus L_2$. 
\item Set $D_1 := \P(L_0 \oplus L_2)$, which is corresponding to the projection 
$E \to L_0 \oplus L_2$. 
\item Set $D_2 := \P(L_0 \oplus L_1)$, which is corresponding to the projection 
$E \to L_0 \oplus L_1$. 
\end{itemize}
By construction, we get $\Gamma_i \cap D_i = \emptyset$ for every $i \in \{0, 1, 2\}$. 

We now show that 
\[
\MO_P(1) \sim D_0 \sim D_1 +\pi^*\MO(1, 0) \sim D_2 + \pi^*\MO(0, 1). 
\]
%We have $K_X +D_0+D_1 + D_2 =\pi^*K_{\P^1 \times \P^1}$ (proof: take the restriction to $S_0$, and apply adjunction). 
Fix $i  \in \{0, 1, 2\}$. 
Note that we have 
$\MO_P(1) \sim D_i + \pi^*M_i$ for some $M_i$. 
By restricting this to $\Gamma_i$, we obtain 
\[
L_i = \MO_P(1)|_{\Gamma_i} =(D_i +\pi^*M_i)|_{\Gamma_i} = M_i. 
\]
Hence we get $\MO_P(1) \sim D_i +L_i$, as required.

It follows from $D_0 \cap D_1 \cap D_2 = \emptyset$ that $|\MO_P(1)|$ is base point free. 
We have 
\[
\MO_P(2) \sim  D_1 +D_2+ \pi^*\MO(1, 1). 
\]
Since $\MO_{\P^1 \times \P^1}(1, 1)$ is ample and $D_1+D_2$ is an effective $\pi$-ample divisor, 
$\MO_P(2)$ is big. Thus (2) holds.

Let us show (3). 
The following holds (cf. \cite[Proposition 7.1(2)]{FanoIII}): 
\[
K_P \simeq \MO_P(-3) \otimes \pi^*(K_{\P^1 \times \P^1} \otimes \det E) 
\simeq \MO_P(-3) \otimes \pi^*\MO(-1, -1). 
\]
Since $\MO_P(1)$ is nef and $\pi$-ample, $-K_P$ is ample. 
Similarly, $-(K_P+X)$ is ample by  
\[
K_P+X \sim \MO_P(-1) \otimes \pi^*\MO(-1, -1). 
\] 
Thus (3) holds. This linear equivalence implies (4). 
Finally, (5) follows from 
\[
-2^i(K_{P}+X)-X \sim \MO_P(2^i-2) \otimes \pi^*\MO(2^i, 2^i). 
\]
\end{proof}

\begin{lem}\label{lem:qFs for (3-10)}
    A smooth Fano threefold $X$ of No.~3-10 is quasi-F-split. 
\end{lem}

\begin{proof}
By Proposition \ref{p-rho3}, we may assume that $p=2$ and $X$ has a wild conic bundle structure.
In what follows, we use the notation of Lemma \ref{description of (3.10)}.
By Proposition \ref{prop:criterion for qFs} and Lemma \ref{lem:SRC}, it suffices to show that 
\begin{enumerate}
    \item[(1)] $H^1(X, \Omega^1_X(K_X))=0$.
    \item[(2)] $H^2(X, \Omega^1_X(-2^iK_X))=0$ for every $i>0$.
\end{enumerate}

\medskip

\noindent\textbf{Step 1: Proof of (1).}\,\,
By the conormal exact sequence, we have the following exact sequence:
\[
0\to \sO_X(K_X-X) \to \Omega^1_P|_X(K_X) \to \Omega^1_X(K_X)\to 0.
\]
Considering the restriction $\sO_P\to \sO_X$, we have an exact sequence
\[
H^2(P, K_P)\to H^2(X, K_X-X) \to H^3(P, K_P-X).
\]
Since $X$ is nef, we have 
$H^2(P, K_P)\simeq H^2(P, \sO_P)=0$ and $ H^3(P, K_P-X)\simeq H^1(P, X)=0$ by \cite[Proposition 1.3]{Totaro(Fano)}.
Thus, we have $H^2(X, K_X-X)=0$. 
Then it suffices to show $H^1(X, \Omega^1_P|_X(K_X))=0$. 

We have the following exact sequence:
\[
0\to \Omega^{1}_{P}(K_{P})\to \Omega^{1}_{P}(K_{P}+X)\to \Omega^{1}_{P}|_X(K_{X})\to 0
\]
Since $-(K_P+X)$ and $-K_P$ are ample, we get 
\[
H^1(P, \Omega^{1}_{P}(K_{P}+X))=0\,\,\,\text{and}\,\,\,H^2(P, \Omega^{1}_{P}(K_{P}))=0
\]
by Bott vanishing. Thus (1) holds. % we obtain the assertion.\\

\medskip

\noindent\textbf{Step 2: Proof of (2).}\,\,
By the conormal exact sequence, we have an exact sequence
\[
0\to \sO_X(-2^iK_X-X) \to \Omega^1_P|_X(-2^iK_X) \to \Omega^1_X(-2^iK_X)\to 0.
\]
Since $K_X=\sO_X(-1)\otimes \pi^{*}\sO_{\P^1\times \P^1}(-1,-1)$ and $X=\sO_P(2)$, we have
\begin{align*}
    H^3(X, -2^iK_X-X)&\simeq H^0(X, (2^i+1)K_X+X)\\
                     &=H^0(X, \sO_P(-2^i+1) \otimes \pi^*\sO_{\P^1\times \P^1}(-2^i-1, -2^i-1))\\
                     &=0.
\end{align*}
Thus it suffices to show $H^2(X, \Omega^1_P|_X(-2^iK_X))=0$.

We have the following exact sequence:
\[
0\to \Omega^{1}_{P}(-2^i(K_{P}+X)-X)\to \Omega^{1}_{P}(-2^i(K_{P}+X))\to \Omega^{1}_{P}|_X(-2^iK_{X})\to 0
\]
Since $-2^i(K_{P}+X)$ and $-2^i(K_{P}+X)-X$ is nef,
we have
\[
H^2(P, \Omega^{1}_{P}(-2^i(K_{P}+X))=0\,\,\,\text{and}\,\,\,H^3(P, \Omega^{1}_{P}(-2^i(K_{P}+X)-X))=0
\]
by \cite[Proposition 1.3]{Totaro(Fano)}. 
Thus (2) holds.
\end{proof}

\subsection{$F$-splitting for 2-2 and 2-6}\label{ss-Cartier-Fsplit}

\begin{prop}\label{p-F-split-meta}
We use the same notation of Lemma \ref{l-dc-to-bdl}. 
Assume that $\dim X = \dim Y=3$. 
Moreover, suppose that the following hold. 
\begin{enumerate}
\item[(0)]
\begin{enumerate}
\item[(0a)] $H^1(Y, pH-mL)=0$ for $m \in \{2, 3\}$. 
\item[(0b)] $H^2(Y, pH-mL)=0$ for $m \in \{1, 2, 3, 4\}$. 
\item[(0c)] $H^3(Y, pH - mL)=0$ for $m \in \{1, 2, 3, 4, 5\}$. 
\end{enumerate}
\item[(1)]
\begin{enumerate}
\item[(1a)] $H^1(Y, \Omega_Y^1(pH-mL))=0$ for $m \in \{1, 2, 3\}$. 
\item[(1b)]  $H^2(Y, \Omega_Y^1(pH-mL))=0$ for $m \in \{2, 3, 4\}$. 
\item[(1c)]$H^3(Y, \Omega_Y^1(pH-4L))=0$. 
\end{enumerate}
\item[(2)]
\begin{enumerate}
\item[(2a)] $H^1(Y, \Omega_Y^2(pH-nL))=0$ for $n \in \{0, 1, 2\}$. 
\item[(2b)] $H^2(Y, \Omega_Y^2(pH-2L))=0$. 
\end{enumerate}
\end{enumerate}
Then the following hold. 
\begin{enumerate}
    \item[(A)] $H^2(X, \Omega_X^1(-pK_X-X))=0$.
    \item[(B)] $H^1(X, \Omega^2_P|_X(-pK_X))=0$.
    \item[(C)] $H^1(X, \Omega_X^2(-pK_X))=0$. 
\end{enumerate}
\end{prop}

\begin{proof}
By taking the wedge product $\bigwedge^2$ of the conormal exact sequence, we have an exact sequence 
\[
0 \to \Omega^1_X(-pK_X-X) \to \Omega^2_P|_X(-pK_X) \to \Omega^2_X(-pK_X) \to 0.
\]
Therefore, we get the following implication: 
\[
{\rm (A)} + {\rm (B)} \Rightarrow {\rm (C)}. 
\]
In what follows, we shall prove (A) and (B).

\noindent\textbf{Step 1: Proof of (A).}\,\,   
By the conormal exact sequence, we have an exact sequence
\[
0 \to \sO_X(-pK_X-2X) \to \Omega^1_P|_X(-pK_X-X) \to \Omega^1_X(-pK_X-X) \to 0.
\]
We recall that $K_X=-g^{*}H|_X =-f^*H$, $X-2S\sim 2g^*L =2f^*L$, $S|_X=0$, and $f_{*}\sO_X=\sO_Y\oplus \MO_Y(-L)$.
We then get 
\begin{align*}
    H^3(X, -pK_X-2X)&=H^3(X, pf^{*}H-4f^{*}L)\\
                    &=H^3(Y, \sO_Y(pH -4L))\oplus H^3(Y, \sO_Y(pH-5L))\\
                    &\overset{{\rm (0c)}}{=} 0. 
\end{align*}
%since $p\geq 7$.
Thus it suffices to show $H^2(X, \Omega^1_P|_X(-pK_X-X))=0$.

Since $K_X=-g^{*}H|_X$, $-X=-2S-2g^{*}L$, and $S|_X=0$, we have the following exact sequence:
\[
0\to \Omega^{1}_{P}(pg^{*}L-4g^{*}L-2S)\to \Omega^{1}_{P}(pg^{*}H-2g^{*}L)\to \Omega^{1}_{P}|_X(-pK_{X}-X)\to 0.
\]
By applying Lemma \ref{l-Omega_P} for $q=2$ and 
$D = pH-2L$, 
(0a) and (1b) imply $H^2(P, \Omega^{1}_{P}(pg^{*}H-2g^{*}L)) =0$. 
Then %it is enough to prove 
the problem is reduced to 
\[
H^3(P, \Omega^{1}_{P}(pg^{*}L-4g^{*}L-2S))=0.
\]
We have the following exact sequence:
\begin{multline*}
0\to g^{*}(\Omega^1_Y(pH-4L))(-2S) \to \Omega^1_P(pg^{*}H-4g^{*}L-2S)\\ \to \Omega^1_{P/Y}(pg^{*}H-4g^{*}L-2S)\simeq \sO_P(pg^{*}H-5g^{*}L-4S) \to 0
\end{multline*}
Thus it is enough to prove that 
\begin{enumerate}
    \item[(I)] $H^3(P, g^{*}(\Omega^1_Y(pH-4L))(-2S))=0$ and 
    \item[(II)] $H^3(P, \sO_P(pg^{*}H-5g^{*}L-4S))=0$.
\end{enumerate}

\medskip

\noindent\textbf{Step 1-1: Proof of (I).}\,\, 
%By $S|_S=-g^{*}L|_S$, 
By using an exact sequence $0 \to \MO_P(-S) \to \MO_P \to \MO_S \to 0$ twice, 
%using the restriction $\sO_P \to \sO_S$ repeatedly,
%$H^3(P, g^{*}(\Omega^1_Y(pH-4L))(-2S))$ can be reduced to those of
the problem is reduced to 
\begin{enumerate}
    \item[(Ia)] $H^2(S, g^{*}(\Omega^1_Y(pH+(n-4)L)))=0$ for $n\in\{0,1\}$ and 
    \item[(Ib)] $H^3(P, g^{*}(\Omega^1_Y(pH-4L)))=0$.
\end{enumerate}
%Let us show (Ia). 
By 
\[
H^2(P, g^{*}(\Omega^1_Y(pH+(n-4)L))) = H^2(Y, \Omega_Y^1(pH+(n-4)L)), 
\]
(Ia) follows from (1b). 
We have 
\[
H^3(P, g^{*}(\Omega^1_Y(pH-4L))) \simeq H^3(Y, \Omega^1_Y(pH-4L)),
\]
and hence (1c) implies (Ib). 
This completes thep proof of (I).

\medskip

\noindent\textbf{Step 1-2: Proof of (II).}\,\, 
By $T \sim S +g^*L$, 
We have $H^3(P, \sO_P(pg^{*}H-5g^{*}L-4S)) \simeq H^3(P, \sO_P(pg^{*}H-g^{*}L-4T))$. 
By using an exact sequence $0 \to \MO_P(-T) \to \MO_P \to \MO_T \to 0$ four times, 
%By $pg^{*}H-5g^{*}L-4S \sim pg^*H -g^*L-4T$, in order to show $H^3(P, \sO_P(pg^{*}H-5g^{*}L-4S))=0$, 
it is enough to prove (IIa) and (IIb) below. 
%can be reduced to those of
\begin{enumerate}
    \item[(IIa)] $H^2(T, pg^{*}H-g^{*}L -nT) =0$ for $n\in\{0,1,2,3\}$.  
    %$H^2(S, pg^{*}H+(-5+n)g^{*}L)\simeq H^2(Y, \sO_Y(p-5+n,p-10+2n))$ for $n\in\{0,1,2,3\}$ and 
    \item[(IIb)] $H^3(P, pg^{*}H-g^{*}L) =0$. 
    %$H^3(S, pg^{*}H-5g^{*}L)\simeq H^3(Y, \sO_Y(p-5,p-10))$.
\end{enumerate}
%\commentbox{{\cred (IIa) does not vanish when $n=0$ and $p=7$?? $H^2(\MO_{\P^2}(p-10+n)) =H^2(\MO_{\P^2}(-3))\neq 0$.}}
By 
\[
H^2(T, pg^{*}H-g^{*}L -nT)  \simeq H^2(Y, pH -(n+1)L), 
\]
(IIa) follows from (0b). 
We have 
\[
H^3(P, pg^{*}H-g^{*}L) \simeq 
H^3(Y, pH-gL), 
\]
and hence (0c) implies (IIb). 
Thus (II) holds.

\medskip

\noindent\textbf{Step 2: Proof of (B).}\,\,  
Since $K_X=-g^{*}H|_X$, $S|_X=0$, and $-X=-2S-2g^{*}L$, we have the following exact sequence:
\[
0\to \Omega^{2}_{P}(pg^{*}H-2g^{*}L)\to \Omega^{2}_{P}(pg^{*}H+2S)\to \Omega^{2}_{P}|_X(-pK_{X})\to 0
\]
Thus it suffices to show that
\begin{enumerate}
    \item[(III)] $H^1(P, \Omega^{2}_{P}(pg^{*}H+2S))=0$ and 
    \item[(IV)] $H^2(P, \Omega^{2}_{P}(pg^{*}H-2g^{*}L))=0$.
\end{enumerate}

\medskip

\noindent\textbf{Step 2-1: Proof of (III).}\,\, 
Taking the wedge product $\bigwedge^2$ 
of the relative exact sequence $0 \to g^*\Omega^1_{Y} \to \Omega^1_P \to \Omega^1_{P/Y} \to 0$, we get
\[
0\to g^{*}\Omega^2_Y \to \Omega^2_P \to g^{*}\Omega^1_Y\otimes \Omega^1_{P/Y} \simeq g^{*}\Omega^1_Y(-g^{*}L-2S) \to 0.
\]
Thus we have
\[
0\to g^{*}(\Omega^2_Y(pH))(2S) \to \Omega^2_P(pg^{*}H+2S)\to g^{*}(\Omega^1_Y(pH-L)) \to 0.
\] 
It holds that 
\[
H^1(P, g^{*}(\Omega^1_Y(pH-L)))=H^1(Y, \Omega^1_Y(pH-L)) 
\overset{{\rm (1a)}}{=}0. 
\]
Then it suffices to show $H^1(P, g^{*}(\Omega^2_Y(pH))(2S)) =0$.  
By  $S|_S=-g^{*}L|_{S}$ and an exact sequence  $0 \to \MO_P(-S) \to \sO_P \to \sO_S \to 0$, 
the problem is reduced to the vanishings of the following:  %the vanishing $H^1(P, g^{*}(\Omega^2_Y(pL))(2S))$ can be reduced to those of
\begin{itemize}
    \item $H^1(S, g^{*}(\Omega^2_Y(pH-nL))\simeq 
    H^1(Y, \Omega^2_Y(pH-nL))$ for $n\in\{1,2\}$. 
    \item $H^1(P, g^{*}(\Omega^2_Y(pH))) \simeq H^1(Y, \Omega^2_Y(pH))$.
\end{itemize}
Both  of them follow from (2a). 
Thus (III) holds.

\medskip

\noindent\textbf{Step 2-2: Proof of (IV).}\,\, 
By the relative exact sequence $0 \to g^*\Omega^1_Y \to \Omega^1_P \to \Omega^1_{P/Y} \simeq \MO_P(-g^*L-2S) \to 0$, we get the following exact sequence:
\[
0\to g^{*}(\Omega^2_Y(pH-2L)) \to \Omega^2_P(pg^{*}H-2g^{*}L) \to g^{*}(\Omega_Y^1(pH-3L))(-2S) \to 0.
\]
We have $H^2(P, g^{*}(\Omega^1_Y(pH-2L)))=H^2(Y, \Omega^2_Y(pH-2L)) 
\overset{{\rm (2b)}}{=}0$. % by Bott vanishing.
By  $S|_S=-g^{*}L|_{S}$ and an exact sequence $0 \to \MO_P(-S) \to \sO_P \to \sO_S \to 0$,
the vanishing of $H^2(P, g^{*}(\Omega_Y(pg^{*}H-3g^{*}L)(-2S))$ can be reduced to those of
\begin{itemize}
    \item $H^1(S, g^{*}(\Omega^1_Y(pH+(-3+n)L)\simeq H^1(Y, \Omega_Y^1(pH+(-3+n)L))$ for $n\in\{0,1\}$ and 
    \item $H^2(P, g^{*}(\Omega^1_Y(pH-3L)))\simeq H^2(Y, \Omega^1_Y(pH-3L))$.
\end{itemize}
These follow from (1a) and (1b). 
Thus (IV) holds. 
\end{proof}

\subsubsection{2-2}

\begin{lem}
    A smooth Fano threefold $X$ of No.~2-2 is $F$-split if $p\geq 7$.
\end{lem}

\begin{proof}
We follow the notation of Lemma \ref{lem:description of (2-2)}. 
It is enough to verify the conditions (1)-(4) in Proposition \ref{prop:criterion for F-split}. 
Proposition \ref{prop:criterion for F-split}(1) holds by 
Lemma \ref{lem:SRC}. 
Lemma \ref{lem:qFs for 2-2} implies 
Proposition \ref{prop:criterion for F-split}(2) and 
Proposition \ref{prop:criterion for F-split}(3). 

It suffices to show 
Proposition \ref{prop:criterion for F-split}(4).  
It is enough to verify the conditions of Proposition \ref{p-F-split-meta}. 
Recall that $Y = \P^1 \times \P^2$, $H = \MO_Y(1, 1)$, and $L =\MO_Y(1, 2)$. 
Since $Y$ is toric, $Y$ satisfies Bott vanishing. 
Then it is enough to check the following (concerning Proposition \ref{p-F-split-meta}(0), use the fact that $pH-4L-K_Y$ is ample): 
\begin{enumerate}
\item[(0)] $H^3(Y, pH-5L)=0$. 
\item[(1)] $H^2(Y, \Omega^1_Y(pH-4L))= H^3(Y, \Omega^1_Y(pH-4L))=0$. 
\end{enumerate}
The assertion (0) follows from 
\[
H^3(Y, pH-5L) =H^3(Y, \MO_Y(p-5, p-10)) \simeq H^1(\P^1, \MO_{\P^1}(p-5)) \otimes H^2(\P^2, \MO_{\P^2}(p-10)) =0. 
\]
Let us show (1). We have 
    \begin{align*}
   \Omega^1_Y(pH-4L) &= \Omega^1_Y(p-4, p-8)\\
    &\simeq \mathrm{pr}_1^{*}\Omega^1_{\P^1}(p-4)\otimes \mathrm{pr}_2^{*}\sO_{\P^2}(p-8)\\
    &\oplus 
    \mathrm{pr}_1^{*}\sO_{\P^1}(p-4)\otimes \mathrm{pr}_2^{*}\Omega^1_{\P^2}(p-8)
    \end{align*}
    Then it holds that 
\begin{align*}
&H^3(Y, \Omega^1_Y(pH-4L)\\
=&H^1(\Omega^1_{\P^1}(p-4))\otimes H^2(\sO_{\P^2}(p-8)) \oplus H^1(\sO_{\P^1}(p-4))\otimes H^2(\Omega^1_{\P^2}(p-8))\\
=& 0  
\end{align*}
and 
\begin{align*}
&H^2(Y, \Omega^1_Y(pH-4L)\\
=&H^0(\Omega^1_{\P^1}(p-4))\otimes H^2(\sO_{\P^2}(p-8)) \oplus 
H^1(\Omega^1_{\P^1}(p-4))\otimes H^1(\sO_{\P^2}(p-8))\\
\oplus &
H^0(\sO_{\P^1}(p-4)\otimes H^2(\Omega^1_{\P^2}(p-8)) \oplus H^1(\sO_{\P^1}(p-4))\otimes H^1(\Omega^1_{\P^2}(p-8))\\
=&0, 
\end{align*}
where we have  $H^2(\Omega^1_{\P^2}(p-8))\simeq H^0(\Omega^1_{\P^2}(8-p))^*=0$ by  $p\geq 7$ and the Euler exact sequence \cite[Ch. II, Example 8.20.1]{Har77}, where $(-)^* := \Hom_k(-, k)$.
Indeed, if $p>7$, then this immediately follows from Bott vanishing. 
When $p=7$, use the Euler exact sequence 
$0 \to \Omega_{\P^2}^1(1) \to \MO_{\P^2}^{\oplus 3} \to \MO_{\P^2}(1) \to 0$, which is still exact even after applying $H^0(\P^2, -)$ by Bott vanishing, and hence $h^0(\Omega_{\P^2}^1(1)) = h^0( \MO_{\P^2}^{\oplus 3}) -h^0(\MO_{\P^2}(1))=0$. 
\end{proof}

\subsubsection{2-6-a}

\begin{lem}
    A smooth Fano threefold of No.~2-6-a is $F$-split if $p\geq 5$.
\end{lem}

\begin{proof}
We follow the notation of Definition \ref{description of 2-6 (a)}. 
It is enough to verify the conditions (1)-(4) in Proposition \ref{prop:criterion for F-split}. 
Proposition \ref{prop:criterion for F-split}(1) holds by 
Lemma \ref{lem:SRC}. 
Lemma \ref{lem:qFs for 2-6b} implies 
Proposition \ref{prop:criterion for F-split}(2) and 
Proposition \ref{prop:criterion for F-split}(3). 

It suffices to show 
Proposition \ref{prop:criterion for F-split}(4).  
By the first paragraph of the proof of Proposition \ref{p-F-split-meta}, 
it is enough to prove  that
\begin{enumerate}
    \item[(1)] $H^2(X, \Omega_X^1(-pK_X-X))=0$.
    \item[(2)] $H^1(X, \Omega_P^2|_X(-pK_X))=0$.
\end{enumerate}

\medskip

\noindent\textbf{Step 1: Proof of (1).}\,\, 
By the conormal exact sequence, we have an exact sequence
\[
0 \to \sO_X(-pK_X-2X) \to \Omega^1_P|_X(-pK_X-X) \to \Omega^1_X(-pK_X-X) \to 0.
\]
We have 
\[
H^3(X, \sO_X(-pK_X-2X))=H^3(X, \sO_X(p-4,p-4))=0
\]
by $p \geq 5$ and Lemma \ref{lem:vanishing for (2-6-1)}.
Thus it suffices to show $H^2(X, \Omega^1_P|_X(-pK_X-X)=0$.

We have the following exact sequence:
\[
0\to \Omega^{1}_{P}(-p(K_{P}+X)-2X)\to \Omega^{1}_{P}(-p(K_{P}+X)-X)\to \Omega^{1}_{P}|_X(-pK_X-X)\to 0.
\]
Then the required vanishing $H^2(X, \Omega^1_P|_X(-pK_X-X)=0$ follows from 
\[
H^2(P, \Omega^{1}_{P}(-p(K_{P}+X)-X))=H^2(P, \Omega^{1}_{P}(p-2,p-2))=0
\] 
and 
\[
H^3(P, \Omega^{1}_{P}(-p(K_{P}+X)-2X))=H^3(P, \Omega^{1}_{P}(p-4,p-4))=0, 
\]
where each vanishing follows from Bott vanishing. Thus (1) holds. 
%Thus, we conclude (1).

\medskip

\noindent\textbf{Step 2: Proof of (2).}\,\,  
We have the following exact sequence:
\[
0\to \Omega^{2}_{P}(-p(K_{P}+X)-X)\to \Omega^{2}_{P}(-p(K_{P}+X))\to \Omega^{2}_{P}|_X(-pK_X)\to 0.
\]
By Bott vanishing, we get 
\[
H^1(P, \Omega^{2}_{P}(-p(K_{P}+X)))=H^1(P, \Omega^{2}_{P}(p,p))=0
\] 
and 
\[
H^2(P, \Omega^{2}_{P}(-p(K_{P}+X)-X))=H^2(P, \Omega^{2}_{P}(p-2,p-2))=0.
\]
Therefore, (2) holds. 
\end{proof}

\subsubsection{2-6-b}

\begin{lem}\label{l-2-6-b-F-split}
    A smooth Fano threefold $X$ of No.~2-6-b is $F$-split if $p\geq 5$.
\end{lem}

\begin{proof}
We follow the notation of Lemma \ref{lem:description of 2-6 (b)}. 
It is enough to verify the conditions (1)-(4) in Proposition \ref{prop:criterion for F-split}. 
Proposition \ref{prop:criterion for F-split}(1) holds by 
Lemma \ref{lem:SRC}. 
Lemma \ref{lem:qFs for 2-6b} implies 
Proposition \ref{prop:criterion for F-split}(2) and 
Proposition \ref{prop:criterion for F-split}(3). 

It suffices to show 
Proposition \ref{prop:criterion for F-split}(4).  
It is enough to verify the conditions of Proposition \ref{p-F-split-meta}. 
Recall that $W \in \MO_{\P^2 \times \P^2}(1, 1)$ and $L =H =\MO_Y(1, 1)$. 
Then all the conditions of Proposition \ref{p-F-split-meta} hold by Lemma \ref{lem:vanishing of W for (2-6-2)} and Serre duality. 
\end{proof}

By a similar argument, we obtain an analogous result for the hyperelliptic case. 
We shall later prove  that the assumption on $p$ in (1) is optimal 
(Example \ref{e-p=11}).

\begin{prop}\label{p-hyperell-F-split}
Let $X$ be a smooth Fano threefold such that $\rho(X)=r_X =1$ and $|-K_X|$ is not very ample, 
where $r_X$ denotes the index of $X$. 
Let $f: X \to Y$ be the double cover induced by $|-K_X|$, where $Y \in \{ \P^3, Q\}$ 
(cf.\ \cite[Theorem 6.5]{FanoI}). 
Then the following hold. 
\begin{enumerate}
\item If $p \geq 13$ and $Y = \P^3$, then $X$ is $F$-split. 
\item If $p \geq 11$ and $Y = Q$, then $X$ is $F$-split. 
\end{enumerate}
\end{prop}

\begin{proof}
We use the same notation of the proof of Proposition \ref{p-hyperell-QFS}. 
By the same argument as in 2-6-b (cf. the proof of Lemma \ref{l-2-6-b-F-split}), 
it is enough to verify Proposition \ref{prop:criterion for F-split}(4), 
which follows from Lemma \ref{lem:vanishing of W for (2-6-2)}. 
\end{proof}

\subsection{Hodge numbers}

\begin{lem}\label{lem:Hodgenumber1(hyperelliptic)}
    In the notation of Section \ref{n-double-cover}, suppose that $(Y,L)=(Q, \sO_{Q}(2)), (W, \sO_W(1,1))$, or $Y$ is toric and $L$ is ample.
    Then $H^0(X,\Omega^2_X)=0$.
\end{lem}
\begin{proof}
    By the exact sequence
    \[
    0\to \Omega^1_X(-X)\to \Omega^2_P|_X \to \Omega^2_X \to 0,
    \]
    it suffices to show 
    \begin{enumerate}
        \item[(a)] $H^0(X,\Omega^2_P|_X)=0$ and 
        \item[(b)] $H^1(X,\Omega^1_X(-X))=0$.
    \end{enumerate}  

    \textbf{(a)}\,\,
    Consider the exact sequence
    \[
    0\to \Omega^2_P(-X)\to \Omega^2_P\to \Omega^2_P|_X\to 0.
    \]
    Since $Y$ is rational, so is $P$, which shows $H^0(P,\Omega^2_P)=0$.
    We show $H^1(P, \Omega^2_P(-X))=0$.
We have an exact sequence
    \[
    0\to g^{*}\Omega^2_Y \to \Omega^2_P \to \Omega^1_{P/Y}\otimes g^{*}\Omega^1_Y\to 0.
    \]
    Recall that $X=2g^{*}L+2S$ and
    $\Omega^1_{P/Y}=\sO_P(-g^{*}L-2S)$.
    Thus, we obtain an exact sequence
    \[
    0\to g^{*}(\Omega^2_Y(-2L))(-2S) \to \Omega^2_P(-X) \to g^{*}(\Omega^1_{Y}(-3L))(-4S)\to 0.
    \]

    Recall $S|_S=-g^{*}L$.
    First, we show $H^1(P, g^{*}(\Omega^2_Y(-2L))(-2S))=0$.
    Using
    \begin{multline*}
        0\to g^{*}(\Omega^2_Y(-2L))(-nS)\to g^{*}(\Omega^2_Y(-2L))(-(n-1)S)\to\\ g^{*}(\Omega^2_Y(-2L))(-(n-1)S)|_S=g^{*}(\Omega^2_Y((n-3)L))\to 0,
    \end{multline*}
    it suffices to show that
    \begin{enumerate}
        \item $H^0(P, g^{*}\Omega^2_Y(-L))=H^0(P, g^{*}\Omega^2_Y(-2L))=0$ and 
        \item $H^1(P, g^{*}\Omega^2_Y(-2L))=0$,
    \end{enumerate}
    which follows from Lemma \ref{lem:vanishing of W for (2-6-2)} or Bott vanishing.

    Next, we show $H^1(P, g^{*}(\Omega^1_{Y}(-3L))(-4S))=0$.
    Using
    \[
    0\to g^{*}(\Omega^1_Y(-3L))(-nS)\to g^{*}(\Omega^1_Y(-3L))(-(n-1)S)\to g^{*}(\Omega^1_Y((n-4)L))\to 0
    \]
    it suffices to show that
    \begin{enumerate}
        \item $H^0(P, g^{*}\Omega^1_Y(-nL))=0$ for $n\in\{0,1,2,3\}$ and 
        \item $H^1(P, g^{*}\Omega^1_Y(-3L))=0$,
    \end{enumerate}
    which follows from Lemma \ref{lem:vanishing of W for (2-6-2)} or Bott vanishing.

    \textbf{(b)}\,\,
    Consider the exact sequence
    \[
    0\to \sO_X(-2X)\to \Omega^1_P|_X(-X)\to \Omega^1_X(-X)\to 0.
    \]
    Since $H^2(X, \sO_X(-2X))=H^2(X, \sO_X(-4g^{*}L))=H^2(Y, \sO_Y(-4L))\oplus H^2(Y, \sO_Y(-5L))=0$,
    it suffices to show that $H^1(X, \Omega^1_P|_X(-X))=0$.
    By
    \[
    0\to \Omega_P^1(-2X)\to \Omega_P^1(-X)\to \Omega_P^1(-X)|_X\to 0
    \]
    It suffices to show that
    \begin{enumerate}
        \item[(b-1)] $H^0(X, \Omega_P^1(-X)|_X)=0$, 
        \item[(b-2)] $H^1(P, \Omega_P^1(-2X))=1$, 
        \item[(b-3)]  $H^1(P, \Omega_P^1(-X))=1$, and
        \item[(b-4)] $H^2(P, \Omega_P^1(-2X))=0$.
    \end{enumerate}

    \textbf{(b-1)}: We omit the proof as it is easy.\\

    \textbf{(b-2)}: We show that $H^1(P, \Omega_P^1(-2X))=1$.\\
    We have
    \[
    0\to g^{*}(\Omega^1_Y(-4L))(-4S)\to \Omega^1_P(-2X) \to \sO_P(-5g^{*}L-6S)\to 0
    \]
    By using the usual restriction exact sequence repeatedly, we have $H^1(P, \sO_P(-5g^{*}L-6S))=H^0(S,\sO_S)=1$.
    We prove $H^{i}(P, g^{*}(\Omega^1_Y(-4L))(-4S))=0$ for $i\in\{1,2\}$. This is reduced to 
    \begin{enumerate}
        \item $H^{i-1}(P, g^{*}(\Omega^1_Y(-nL)))=0$ for $n\in\{1,2,3,4\}$ and 
        \item $H^{i}(P, g^{*}(\Omega^1_Y(-4L)))=0$
    \end{enumerate}
    which follows from Lemma \ref{lem:vanishing of W for (2-6-2)} or Bott vanishing.

    \textbf{(b-3)}: We show that $H^1(P, \Omega_P^1(-X))=1$.\\
    By 
    \[
    0\to g^{*}\Omega^1_Y\to \Omega^1_P \to \Omega^1_{P/Y}=\sO_P(-g^{*}L-2S)\to 0
    \]
    we have
    \[
    0\to g^{*}(\Omega^1_Y(-2L))(-2S)\to \Omega^1_P(-X) \to \sO_P(-3g^{*}L-4S)\to 0.
    \]
    It is easy to see that $H^1(P, \sO_P(-3g^{*}L-4S))=1$. 
    We show
    $H^i(P, g^{*}(\Omega^1_Y(-2L))(-2S))=0$ for $i\in\{1,2\}$. Considering restriction to $S$, the vanishing is reduced to 
    \begin{enumerate}
        \item $H^{i-1}(P, g^{*}\Omega_Y^1(-L))=H^i(P, g^{*}\Omega_Y^1(-2L))=0$ and 
        \item $H^i(P, g^{*}\Omega_Y^1(-2L))=0$,
    \end{enumerate}
    which follows from Lemma \ref{lem:vanishing of W for (2-6-2)} or Bott vanishing.

    \textbf{(b-4)}: We show $H^2(P, \Omega_P^1(-2X))=0$.
   We have
    \[
    0\to g^{*}(\Omega^1_Y(-4L))(-4S)\to \Omega^1_P(-2X) \to \sO_P(-5g^{*}L-6S)\to 0
    \]
    It is easy to see $H^2(P, \sO_P(-5g^{*}L-6S))=0$.
    The vanishing $H^2(P, g^{*}(\Omega^1_Y(-4L))(-4S))=0$ has been proven in (b-2). Thus, we conclude.
    % \begin{enumerate}
    %     \item $H^1(P, g^{*}(\Omega^1_Y(-nL)))=0$ for $n\in\{1,2,3,4\}$ and 
    %     \item $H^2(P, g^{*}(\Omega^1_Y(-4L)))=0$
    % \end{enumerate}
    % which follows from \cite[Lemma 6.13]{Kawakami-Tanaka3}.
    
\end{proof}

\begin{lem}\label{lem:Hodgenumber2(hyperelliptic)}
    In the notation of Section \ref{n-double-cover}, suppose that $(Y,L)=(Q, \sO_{Q}(2)), (W, \sO_W(1,1))$, or $Y$ is toric and $L$ is ample.
    Then $h^{1,1}(X)=\rho(X)$.
\end{lem}
\begin{proof}
    Consider the short exact sequence
    \[
    0\to \sO_X(-X) \to \Omega^1_P|_X \to \Omega^1_X\to 0.
    \]
    Since $H^i(\sO_X(-X))=H^i(Y, g_{*}\sO_X(-2g*L))=H^i(Y, \sO_Y(-3L))\oplus H^i(Y, \sO_Y(-4L))=0$ for $i\in\{1,2\}$, we have $H^1(\Omega^1_X)=H^1(\Omega^1_{P}|_X)$.
    We also have $H^0(\Omega^1_P|_X)=0$.

    Consider the short exact sequence
    \[
    0\to \Omega^1_P(-X) \to \Omega^1_P \to \Omega^1_P|_X \to 0.
    \]

    We show $H^2(P, \Omega^1_P(-X))=0$.
    By the short exact sequence
    \[
    0\to g^{*}\Omega^1_Y\to \Omega^1_P \to \Omega^1_{P/Y}=\sO_P(-g^{*}L-2S)\to 0
    \]
    we have the short exact sequence
    \[
    0\to g^{*}(\Omega^1_Y(-2L))(-2S)\to \Omega^1_P(-X) \to \sO_P(-3g^{*}L-4S)\to 0
    \]
    Since $H^2(P, \sO_P(-3g^{*}L-4S))=0$, it suffices to show
    $H^2(P, g^{*}(\Omega^1_Y(-2L))(-2S))=0$. Considering restriction to $S$, the vanishing is reduced to 
    \begin{enumerate}
        \item $H^1(P, g^{*}\Omega_Y^1(-L))=H^1(P, g^{*}\Omega_Y^1(-2L))=0$ and 
        \item $H^2(P, g^{*}\Omega_Y^1(-2L))=0$,
    \end{enumerate}
    which follows from Lemma \ref{lem:vanishing of W for (2-6-2)} or Bott vanishing.\\

    We show $H^1(P, \Omega^1_P(-X))=1$.
    We have the short exact sequence
    \[
    0\to g^{*}(\Omega^1_Y(-2L))(-2S)\to \Omega^1_P(-X) \to \sO_P(-3g^{*}L-4S)\to 0
    \]
    Then considering restriction to $S$, we have $H^1(P, \sO_P(-3g^{*}L-4S))=1$.
We also have $H^i(P, g^{*}(\Omega^1_Y(-2L))(-2S))=0$ for $i\in\{1,2\}$. 
    Thus, $H^2(P, \Omega^1_P(-X))=1$.

    We have $h^1(\Omega^1_P)=h^0(\sO_Y)+h^1(\Omega^1_Y)=1+\rho(Y)=\rho(P)$, where the first equality follows from \cite[Chapter1, Corollaire 4.2.13]{Gros}. Thus, we obtain $h^1(X, \Omega^1_P|_X)=h^1(\Omega^1_P)-1 =\rho(Y) =\rho(X)$ 
   (see \cite[Section 7]{FanoIV} for the last equality). 
\end{proof}

\subsection{Akizuki-Nakano vanishing}
\begin{lem}\label{lem:ANV(hyeprelliptic)}
    Let $X$ be a smooth Fano threefold such that $\mathrm{Pic}(X)=\mathbb{Z} K_X$ and $-K_X$ is not very ample.
    Then $H^1(X, \Omega^1_X(nK_X))=0$ for all $n>0$.
\end{lem}
\begin{proof}
   If $g=2$ (resp.~$g=3$), then $Y\simeq \P^3$ (resp.~$\Q^3$) and $L=\sO_{\P^3}(3)$ (resp.~$\sO_{\Q^3}(2)$) and $H=\sO_{\P^3}(1)$ (resp.~$\sO_{\Q^3}(1)$).
   Now, the assertion follows from Theorem \ref{p-2:1-QFS} (I) and
   Lemma \ref{lem:vanishing of W for (2-6-2)}.
\end{proof}

%% file: section7.tex
\section{Proofs of the main theorems}

In this section, we prove the main theorems in the Introduction.

\begin{thm}[\textup{Theorem \ref{Introthm:quasi-F-split}}]\label{thm:quasi-F-split}
    Let $X$ be a smooth Fano threefold such that $\rho(X)>1$ or $r_X>1$.
    Then $X$ is quasi-$F$-split.
\end{thm}
\begin{proof}
If $r_X \geq 3$, then $X \simeq \P^3$ or $X$ is isomorphic to a smooth quadric threefold, and hence $X$ is $F$-split. 
If $r_X =2$, then $X$ is quasi-$F$-split by \cite[Theorem A and Remark 2.8]{Kawakami-Tanaka(dPvar)}. 
Therefore, we may assume that $\rho(X) \geq 2$. 
In this case, 
the assertion  holds by former parts as follows: 
\begin{itemize}
\item $\rho(X) \geq 6$: Proposition \ref{p-rho6}. 
\item $\rho(X) = 5$:  Proposition \ref{p-rho5}. 
\item $\rho(X) = 4$:  Proposition \ref{p-rho4}. 
\item $\rho(X) = 3$:  Proposition \ref{p-rho3} and Lemma \ref{lem:qFs for (3-10)}. 
\item $\rho(X) = 2$:  Proposition \ref{p-rho2}, Lemma \ref{lem:qFs for 2-2}, 
Lemma \ref{lem:qFs for 2-6a}, Lemma \ref{lem:qFs for 2-6b}, and Lemma \ref{lem:qFs for (2-8)}. 
\end{itemize}

% We now show (2), i.e., Kodaira vanishing for $X$.
% When $\rho(X)\geq 2$ or $r_X\geq 2$, this follows from (1) and Kodaira vanishing for quasi-$F$-split varieties \cite[Theorem 3.15]{KTTWYY1}. 
% When $\rho(X)=r_X=1$, then Kodaira vanishing holds by \cite[Corollary 4.5]{FanoI}.
\end{proof}

\begin{thm}[\textup{Theorem \ref{Introthm:ANV}(=\cite[Theorem B]{Kawakami-Tanaka(Lift1)})}]\label{thm:ANV}
    Let $X$ be a smooth Fano threefold over an algebraically closed field $k$ of characteristic $p>0$. 
Then Akizuki-Nakano vanishing holds on $X$, that is, 
if $A$ is an ample  Cartier divisor $A$ on $X$, then we have
\[
H^j(X,\Omega^i_X(-A))=0
\]
for all integers $i,j\geq 0$ satisfying  $i+j<3$.
\end{thm}
\begin{proof}
    If $\rho(X)>1$ or $r_X>1$, then $X$ is quasi-$F$-split.
    Thus the assertion follows from \cite[Corollary 4.10]{Petrov}.
    If $\mathrm{Pic}(X)=\mathbb{Z} K_X$ and $-K_X$ is not very ample, then the assertion %follows from 
    holds by Lemma \ref{lem:ANV(hyeprelliptic)}.
    Finally, if $\mathrm{Pic}(X)=\mathbb{Z} K_X$ and $-K_X$ is very ample, then the assertion follows from \cite[Theorem 6.3]{Kawakami-Tanaka(Lift1)}.
\end{proof}

\begin{thm}[\textup{Theorem \ref{Introthm:W(k)-lift}(=\cite[Theorem A]{Kawakami-Tanaka(Lift1)})}]\label{thm:lift}
    Let $X$ be a smooth Fano threefold over an algebraically closed field $k$ of positive characteristic.
Then $X$ lifts to $W(k)$.
\end{thm}
\begin{proof}
    Since $H^2(X,T_X)\simeq H^1(X,\Omega^1_X(K_X))=0$ and 
    an ample invertible sheaf $\omega_X^{-1}$ %Cartier divisor $-K_X$ 
    lifts to $W(k)$, we conclude the assertion by Theorem \ref{thm:ANV} 
    %, we conclude by 
    and \cite[Theorem 8.5.19]{FAG}.
\end{proof}

% \begin{thm}\label{mainthm:ANV}
%     Let $X$ be a smooth Fano threefold such that 
%     Suppose that $g\neq 6,9$ (Here, we assume $\rho(X)=r_X=1$).
%     Then $X$ satisfies Akizuki-Nakano vanishing.
% \end{thm}
% \begin{proof}
%     Fix ample Cartier divisor $A$ on $X$.
%     We show $H^j(X, \Omega^i_X(-A))=0$ for $i+j <3$. 
% If $i=0$ (resp. $(i, j) = (1, 0)$,  resp. $(i, j) = (2, 0)$), then this follows from Theorem \ref{thm:quasi-F-split}(2) 
% \end{proof}

% \begin{thm}\label{mainthm:lift}
%     Let $X$ be a smooth Fano threefold.
%     Then $X$ lifts to $W(k)$.
% \end{thm}
% \begin{proof}
%     If $g\neq 6,9$, then $H^2(X,T_X)\simeq H^0(X,\Omega^1_X(K_X))=0$ by Theorem \ref{mainthm:ANV}.
%     When $g=6,9$, we have already seen that $X$ lifts to $W(k)$ in \cite{Kawakami-Tanaka3}
% \end{proof}

\begin{thm}[\textup{Theorem \ref{Introthm:hoge number}(=\cite[Theorem D]{Kawakami-Tanaka(Lift1)})}]\label{thm:hodge}
Let $X$ be a smooth Fano threefold over an algebraically closed field $k$ of positive characteristic. 
Take a lift $f\colon \mathcal{X}\to W(k)$ of $X$ to $W(k)$, whose existence is ensured by Theorem \ref{thm:lift}. 
Let $X_{\overline{K}}$ be the geometric generic fibre over $W(k)$. 
Then all the Hodge numbers $h^{i,j}(X)\coloneqq \dim_{k} H^j(X,\Omega^i_X)$ of $X$ coincide with those of $X_{\overline{K}}$, that is, 
\[
h^j(X,\Omega^i_X)=h^j(X_{\overline{K}},\Omega^i_{X_{\overline{K}}})
\]
holds for all $i,j\geq 0$.
\end{thm}
\begin{proof}
    As in the proof of \cite[Theorem 6.2]{Kawakami-Tanaka(Lift1)}, it suffices to show that $H^0(X,\Omega^2_X)=0$ and $h^1(X, \Omega^1_X)=\rho(X)$.
    By \cite[equation (1.5)]{Hodgecohomology}, we may assume that $X$ is 
    a primitive Fano threefold (for the definition, see \cite[Subsection 2.1(6)]{FanoIII}). 
    %primitive. 
    If $\mathrm{Pic}(X)=\mathbb{Z} K_X$ and $-K_X$ is not very ample, then the assertion follows from Lemma \ref{lem:Hodgenumber1(hyperelliptic)} and Lemma  \ref{lem:Hodgenumber2(hyperelliptic)}.
    If $\mathrm{Pic}(X)=\mathbb{Z} K_X$ and $-K_X$ is very ample, then the assertion follows from \cite[Theorem 6.2]{Kawakami-Tanaka(Lift1)}. 
    Hence we may assume that $\mathrm{Pic}(X) \neq \mathbb{Z} K_X$, i.e., 
    $r_X \geq 2$ or $\rho(X) \geq 2$.

 Assume that $r_X\geq 2$. Then the assertion follows from Lemma \ref{lem:ANV for hypersurfaces} for $V_1,V_2$,
 from \cite[Theorem 2.9]{Kawakami-Tanaka(Lift1)} for $V_3,V_4$, and from \cite[Lemma 8.2 and Theorem 8.13]{Kawakami-Tanaka(Lift1)} for the other 
 cases. % $(-K_X/2)^3 \geq 5$. %ano threefolds.

    \medskip
    
    Assume that $\rho(X)\geq2$. 
    By \cite[Theorem 1.1]{FanoIII},
    $X$ is one of \[\{\text{2-2, 2-6, 2-8, 2-18, 2-24, 2-32, 2-34, 2-35, 2-36, 3-1, 3-2, 3-27, 3-31}\}.\]
    If $X$ is one of $\{\text{2-6-a, 2-24, 2-32, 2-34, 2-35, 2-36, 3-27, 3-31}\}$, then the assertion follows from \cite[Theorem 2.9]{Kawakami-Tanaka(Lift1)}.
    If $X$ is one of \[\{\text{2-2, 2-6-b, 2-8, 2-18, 3-1}\},\]
    then the assertion follows from Lemmas \ref{lem:Hodgenumber1(hyperelliptic)} and \ref{lem:Hodgenumber2(hyperelliptic)}.
    
    Finally, if $X$ is 3-2, then it is rational.
    In fact, by \cite[Proposition 4.32 and its proof]{FanoIV}, there exists a morphism $\pi : X \to \P^1$ 
 such that $\pi_*\MO_X = \MO_{\P^1}$ 
 and $(-K_X)^2 \cdot F =6$ for a fibre $F$ of $\pi$. 
Then $X$ is rational by \cite[Theorem 1.4]{BT24}.
Now, the assertion follows from \cite[Lemma 8.2 and Theorem 8.13]{Kawakami-Tanaka(Lift1)}.
\end{proof}

\begin{thm}[\textup{Theorem \ref{Introthm:E_1-degeneration}(=\cite[Theorem C]{Kawakami-Tanaka(Lift1)})}]\label{thm:degenerate}
Let $X$ be a smooth Fano threefold over an algebraically closed field $k$ of characteristic $p>0$. 
Then the following hold. 
\begin{enumerate}
    \item The Hodge to de Rham spectral sequence 
\[
E_1^{i,j} 
= H^j(X, \Omega_X^i) 
\Rightarrow 
H^{i+j}(X, \Omega_X^{\bullet}) =E^{i+j}
\]
degenerates at $E_1$.
    \item Crystalline cohomology $H^i_{\mathrm{cris}}(X/W(k))$ is torsion-free for every $i\geq 0$.
\end{enumerate}
\end{thm}
\begin{proof}
    The assertion follows from Theorems \ref{thm:lift} \ref{thm:hodge}, and \cite[Proposition 6.5]{Kawakami-Tanaka(Lift1)}.
\end{proof}

\begin{proof}[Proof of Theorem \ref{Introthm:log ANV}]
    The assertion follows from Theorem \ref{Introthm:quasi-F-split} and \cite[Corollary 7.6]{KTTWYY1}.
\end{proof}

\begin{proof}[Proof of Theorem \ref{intro-F-split}]
The assertion follows from Theorem \ref{t-GFS} and Section \ref{ss-Cartier-Fsplit}. 
\end{proof}

%% file: section8.tex
\section{Examples}

In this section, we gather examples of non-$F$-split or non-quasi-$F$-split smooth Fano threefolds.

\begin{example}[$X = S \times \P^1$]
Let $S$ be a smooth del Pezzo surface which is not $F$-split. 
Then $X \coloneqq S \times \P^1$ is a Fano threefold which is not $F$-split. 
Therefore, 
if the characteristic $p$ of the base field $k$ and 
$\rho$ satisfies one of (1)-(3) below, 
then there exists a non-$F$-split smooth Fano threefold $X$ over $k$ 
satisfying  $\rho(X)=\rho$.
\begin{enumerate}
\item $p=2$ and $\rho=7$
\item $p \in \{2, 3\}$ and $\rho=8$. 
\item $p \in \{2, 3, 5\}$ and $\rho=9$. 
\end{enumerate}
\end{example}

\begin{example}[Wild conic bundles]
%$p=2, \rho \in \{2, 3\}$, non-$F$-split]
Wild conic bundles are not $F$-split. 
Indeed, if $X$ is $F$-split and $f: X \to S$ is a conic bundle, 
then $f$ is generically reduced \cite[Lemma 2.4]{GLP15} and hence not wild.  
Therefore, if $p=2$ and $X$ is a smooth Fano threefold which is 2-24 or 3-10, 
then $X$ is not necessarily $F$-split. 
\end{example}

\begin{example}[$p=7$, non-$F$-split]\label{e-p=7}
%[No. 2-1, $p=5$, non-$F$-split]
Assume $p=7$. 
Then 
\[
X := \{ x_0^4 + x_1^4+x_2^4 +x_3^4 + x_4^4 =0\} \subset \P^4
\]
is not $F$-split. Indeed, we have 
\[
(x_0^4 + x_1^4+x_2^4 +x_3^4 + x_4^4)^6 \not\in (x_0^7, x_1^7, x_2^7, x_3^7, x_4^7), 
\]
and hence Fedder's criterion \cite[Proposition 2.1]{Fed83} implies that $X$ is not $F$-split (cf.\ \cite[Proposition A.8]{KTY}). 
\end{example}

\begin{example}[$p=11$, non-$F$-split]\label{e-p=11}
%[No. 2-1, $p=5$, non-$F$-split]
Assume $p=11$. 
Take 
\[
X := \{ x_0^6 + x_1^6+x_2^6 +x_3^6 + y^2 =0\} \subset \P(1, 1, 1, 1, 3),
\]
i.e., $\P(1, 1, 1, 1, 3) = \Proj\,k[x_0, x_1, x_2, x_3, y]$ with $\deg x_i = 1$ and $\deg y =3$ for every $0 \leq i \leq 3$, and 
\[
X := \Proj\,\frac{k[x_0, x_1, x_2, x_3, y]}{(x_0^6 + x_1^6+x_2^6 +x_3^6 + y^2)}. 
\]
Let us prove that 
\begin{enumerate}
\item $X$ is not $F$-split, and 
\item $X$ is a smooth Fano threefold. 
\end{enumerate}
The assertion (1) follows from Fedder's criterion \cite[Proposition 2.1]{Fed83} (cf.\ \cite[Proposition A.8]{KTY}). 

Let us show (2). 
It is enough to prove that $X$ is smooth, as the other assertions in (2) follow from the adjunction formula and the fact that $X$ is an ample $\Q$-Cartier effective Weil divisor on $\P(1, 1, 1, 1, 3)$ (which implies the connectedness of $X$). 
Suppose that $[a_0:a_1:a_2:a_3:b]$ is a singular point of $X$, where 
$a_0, a_1, a_2, a_3, b \in k$. 
Recall that we have 
$[a_0:a_1:a_2:a_3:b] = [\lambda a_0:\lambda a_1:\lambda a_2:\lambda a_3:\lambda^3b]$  for every $\lambda \in k \setminus \{0\}$. 
For $X_1 := x_1/x_0, X_2 :=x_2/x_0, X_3 :=x_3/x_0, Y := y/x_0^3$, 
we have 
\[
D_+(x_0) = \Spec k\left[X_1, X_2, X_3, Y\right] (\simeq \A^4) \subset \P(1, 1, 1, 1, 3). 
\]
and 
\[
X \cap D_+(x_0) = \{ 1+X_1^6 +X_2^6 + X_3^6 +Y^2=0\}, 
\]
which is smooth. 
Then it holds that $a_0=0$. 
By symmetry, we get $a_0=a_1 = a_2= a_3=0$, which implies $[a_0:a_1:a_2:a_3:b] = [0:0:0:0:b] = [0:0:0:0:1]$. 
However, $X$ does not pass through $[0:0:0:0:1]$, which is absurd. 
Thus (2) holds. 
\end{example}

\begin{example}[No. 2-3, $p=3$, non-$F$-split]\label{e-2-3-p=3}
Assume $p=3$. 
We construct a Fano threefold $X$ 
which is 2-3 and not $F$-split. 
Let $V_2$ be a Fano threefold of index $2$ such that 
$(-K_{V_2})^3 =16$ and $V_2$ is not $F$-split 
(e.g., $V_2 :=\{ x_0^4 + x_1^4 + x_2^4+x_3^4 + y^2 =0\}$ in $\P(1, 1, 1, 1, 2) = \Proj\,k[x_0, x_1, x_2, x_3, y]$ for $\deg x_0 =\deg x_1 = \deg x_2 = \deg x_3$ and $\deg y=2$, cf. the proof of Example \ref{e-p=11}). 
Take a Cartier divisor $H$ on $V_2$ such that $2H \sim -K_{V_2}$. 
Then $|H|$ is base point free and it induces a finite double cover. % $V_2 \to Z$. 
By a Bertini theorem \cite[Corollary 4.3]{Spr98}, we may assume that $H$ is a smooth prime divisor on $V_2$, 
which is a smooth del Pezzo surface with $K_H^2=2$. 
Pick a general member $C$ of $|-K_H|$, which is a smooth elliptic curve 
\cite[Theorem 1.4]{KN}. 
By an exact sequence 
\[
H^0(V_2, \MO_{V_2}(H)) \to H^0(H, \MO_{V_2}(H)|_H) \to H^1(V_2, \MO_{V_2}(H -H)) = H^1(V_2, \MO_{V_2}) =0, 
\]
there exists a member $H' \in |H|$ such that $H \cap H' = H'|_H = C$. 
Let $\sigma : X \to V_2$ be the blowup along $C$. 
Since $\sigma$ coincides with the resolution of the indeterminacies of the pencil generated by $H$ and $H'$, 
there is a contraction $\pi : X \to \P^1$ of type $D$ such that the proper transforms of $H$ and $H'$ are fibres of $\pi$.  
By Kleimann's criterion, $X$ is a smooth Fano threefold, which is of No.~2-3. 
\end{example}

\begin{example}[No. 2-1, $p=5$, non-$F$-split]\label{e-2-1-p=5}
Assume $p=5$. 
We construct a Fano threefold $X$  
which is 2-1 and $X$ is not $F$-split. 
Let $V_1$ be a Fano threefold of index $2$ such that 
$(-K_{V_1})^3 = 8$ and $V_1$ is not $F$-split. 
We can find such an example by setting 
\[
V_1 :=\{ x_0^6 + x_1^6 + x_2^6 + y^3+z^2 =0\} \subset \P(1, 1, 1, 2, 3) = \Proj\,k[x_0, x_1, x_2, y, z],\]
where $\deg x_0 =\deg x_1 = \deg x_2 =1, \deg y=2, \deg z=3$ (cf.\ \cite[Section 3.1]{Oka21}). 
Take a Cartier divisor $H$ on $V_1$ such that $2H \sim -K_{V_1}$. 
%Set $H:= -(1/2)K_{V_2}$. 
Then we have a scheme-theoretic equality ${\rm Bs}\,|H| =P$  for some  closed point $P$ on $X$. 
We take generic members $H^{\gen}_1$ and $H^{\gen}_2$ of $|H|$ twice. 
Then 
%is regular and smooth around $P$ {\cred ref, Fujita?}. 
%Again by taking the generic member $C$ on $H^{\gen}$, 
$C = H^{\gen}_1 \cap H^{\gen}_2$ is a regular curve of genus one. 
By $p=5>3$, $C$ is a smooth elliptic curve \cite[Corollary 1.8]{PW22}. 
Therefore, the intersection $C := H_1 \cap H_2$ of two general members $H_1$ and $H_2$ of $|H|$ is a smooth elliptic curve. 
Take the blowup $X \to V_1$ along $C$. 
Then we can apply the same argument as in Example \ref{e-2-3-p=3}. 
\end{example}

\begin{example}[$p=5$, non-quasi-$F$-split]\label{e-p=5-nonQFS}
Assume $p=5$. 
Take 
\[
X := \{ x_0^6 + x_1^6+x_2^6 +x_3^6 + y^2 =0\} \subset \P(1, 1, 1, 1, 3),
\]
i.e., $\P(1, 1, 1, 1, 3) = \Proj\,k[x_0, x_1, x_2, x_3, y]$ with $\deg x_i = 1$ and $\deg y =3$ for every $0 \leq i \leq 3$, and 
\[
X := \Proj\,\frac{k[x_0, x_1, x_2, x_3, y]}{(x_0^6 + x_1^6+x_2^6 +x_3^6 + y^2)}. 
\]
Then $X$ is a smooth Fano threefold by the same proof as in Example \ref{e-p=11}. 
Moreover, $X$ is not quasi-$F$-split by \cite[Corollary 4.19(i), Proposition A.8]{KTY}. 
\end{example}